\documentclass[12pt]{article}
\usepackage{amsmath, amssymb}
\usepackage{amsfonts}
\usepackage[dvips]{graphicx}

\newtheorem{theorem}{Theorem}[section]
\newtheorem{lemma}{Lemma}[section]
\newtheorem{definition}{Definition}[section]

\def\QED{\mbox{\rule[0pt]{1.5ex}{1.5ex}}}

\def\endproof{\hspace*{\fill}~\QED\par\endtrivlist\unskip}
\def\keywords{\vspace{-.3em}
    \if@twocolumn
      \small\it Keywords\/\bf---$\!$%
    \else
      \begin{center}\small\bf Keywords\end{center}\quotation\small
    \fi}
\def\endkeywords{\vspace{0.6em}\par\if@twocolumn\else\endquotation\fi
    \normalsize\rm}
\newtheorem{rem}{Remark}[section]

\newenvironment{indention}[1]{\addtolength{\leftskip}{#1}\addtolength{\rightskip}{#1}\begingroup}{\endgroup\par}

\def\Label#1{\label{#1}\ [\ #1 \ ]\ }
\def\Label{\label}
\title{\Large Practical implementation and error bound \\ of  integer-type algorithm \\ for higher-order differential equations}
\author{Fuminori SAKAGUCHI\thanks{Faculty of Engineering, University of Fukui, 3-9-1 Bunkyo, Fukui 910-8507, Japan ({\tt fsaka@u-fukui.ac.jp})} \and Masahito HAYASHI\thanks{Graduate School of Information Sciences, Tohoku University, Sendai 980-8579, Japan ({\tt hayashi@math.is.tohoku.ac.jp}) Centre for Quantum Technologies, National University of Singapore, 3 Science Drive 2, Singapore 117542}}

\date{}
\begin{document}
% \documentclass{siamltex}
% \usepackage{amsfonts,graphicx,color}
% \usepackage{graphics}

% \pagestyle{plain} \thispagestyle{plain}

% \newcommand{\R}[1]{{\color{red} #1}}

% \begin{document}

% \newtheorem{rem}{Remark}[section]

\maketitle
\renewcommand{\thesection}{\arabic{section}}
\begin{abstract}
In our preceding paper, we have proposed an algorithm for obtaining  
finite-norm solutions of higher-order linear ordinary differential equations of the Fuchsian type $\displaystyle\Bigl( \sum_{m=0}^M p_m (x) \bigl(\frac{d}{dx}\bigr)^m \Bigr) f(x) = 0$
(where $p_m$ is a rational function with rational-number-valued coefficients), by using only the four arithmetical operations on integers, 
and we proved its validity. 
For any nonnegative integer $\displaystyle\, k$, it is guaranteed mathematically that this method can produce all the solutions satisfying $\,\int |f(x)|^2 \, (x^2+1)^k \, dx < \infty \,$, under some conditions.  
%at least when $p_M(x)$ has no real root.
We materialize this algorithm in practical procedures. An integer-type quasi-orthogonalization used there can suppress the explosion of calculations.
% This method is based on a kind of orthonormal system of localized `quasi-sinusoidally' oscillating wavepackets with a spindle-shaped envelope.  We can extract almost only the true converging components for the solutions by using integer-valued quasi-orthogonalization processes based on an idea similar to a combination of the Gram-Schmidt orthogonalization and the Euclidean algorithm. Moreover, this integer-valued quasi-orthogonalization method can suppress the explosion of calculations which may appear in the usual optimization methods and in the exact orthogonalization.  In this paper, we explain how to realize this method in a practical algorithm and the reason that we can extract almost only true solution components. 
Moreover, we give an 
%an estimation of the 
upper limit of the errors.  We also give some results of numerical experiments and compare them with the corresponding exact analytical solutions, which show that the proposed algorithm is successful in yielding solutions with an extraordinarily high accuracy (using only arithmetical operations on integers).
\end{abstract}

AMS: 65L99, 42C15, 65L70, 65L60, 34A45 
\pagestyle{myheadings}
\thispagestyle{plain}
\markboth{F. SAKAGUCHI AND M. HAYASHI}{Practical integer-type algorithm for higher-order differential equations}

\maketitle
\vspace{8mm}

\thispagestyle{plain}

\section{Introduction} 
Since linear higher-order ordinary equations 
with general coefficient functions cannot be solved analytically except in very special simple cases
~\cite{CoLe}, numerical methods are useful in many applications, such as the eigenfunction problem for a differential operator  on a complete  function space ${\cal H}$. Many approaches with numerical methods have been proposed for this problem. One method is based on the approximation of solutions in a finite-dimensional subspace of the original function space ${\cal H}$
~\cite{KrVa}~\cite{Bre}.

Using an idea somewhat similar to this approach, our preceding paper~\cite{paper1} 
proved the validity of our newly-proposed 
%proposed an 
integer-type numerical method 
using smooth (i.e. analytical) basis functions 
which was able to obtain approximations of all the true solutions in 
${\cal H}$ 
%$C^M(\mathbb{R} )\cap {\cal H}$ 
of the 
Fuchsian-type 
differential equation   $\displaystyle\Bigl( \sum_{m=0}^M r_m (x) \bigl(\frac{d}{dx}\bigr)^m \Bigr) f(x) = 0$ with rational-function-type 
coefficient functions $r_m (x)$ 
using only the four arithmetical operations on integers, for a class of complete function spaces 
${\cal H}$ 
containing $L^2(\mathbb{R} )$ as a special case. 
This algorithm can be applied even for the 
non-Fuchsian cases, 
%cases where $p_M (x)$ and the denominators of $p_m(x)$ $(m=0,1,\ldots , M)$ has zero points, 
under some restrictions. 
However, the method proposed in that paper is considerably different from the usual methods based on approximation in a finite-dimensional subspace, such as the Ritz and Galerkin methods~\cite{KrVa}~\cite{Bre}. 
The main differences from usual Galerkin methods were explained in the introduction of~\cite{paper1}. 
The method proposed in~\cite{paper1} can be regarded as a kind of `semi-analytical method' rather than a purely numerical method, in that it is closely related to the functional analysis, Fourier series and the Laurent expansion of complex functions. In addition, as a remarkable characteristic of the proposed method, all the basis functions used there are rational functions of the coordinate which are related to a power series of the Cayley transform of the coordinate. This characteristic enables us to close all procedures in the method only within four arithmetic operations between rational numbers and hence between integers. These facts 
imply 
%implies 
that the proposed method can be discussed from some viewpoints of mathematical analysis. 

This method is perfectly free from round-off error because it consists only of integer-type 
operations, 
when we choose function spaces and their basis systems appropriately.  
%algorithms. 
It has only two types of errors 
instead of round-off errors. 
One is the `pure' truncation error due to the components outside the subspace 
(contained even in {\it true exact} solutions), 
and the other is the `mixture error' due to the (slight) mixture of extra solutions not corresponding to true solutions in ${\cal H}$.  
However, as is proved in ~\cite{paper1}, the latter  mixture error converges to zero as the dimension of the subspace tends to infinity, with this method. 
Moreover, the `pure' truncation error decays very rapidly in the proposed method, due to the relationship between the Fourier series and the expansion used in the proposed method.

Moreover, this method requires a small amount of calculations for obtaining high-accuracy solutions. For example, when the coefficients in the expansion of a true solution by the basis functions decay exponentially, the amount of calculations required by this method is almost proportional asymptotically to the cube of the number of required significant digits. 
This is a strong advantage in comparison with the Runge-Kutta methods and the finite elements methods which require the amount of calculations with exponential orders of the required significant digits.  
From this advantage, as is shown later in numerical results, by the proposed method, we can easily attain the accuracy with several hundred or several thousand significant digits by an ordinary personal computer in many cases. 

Since the structure and the procedures of this method are somewhat complicated, in our preceding paper ~\cite{paper1}, we explained only its abstract structure, its mathematical validity and a rough sketch of the procedures in the proposed algorithms. In this paper, we will explain the concrete materialization of this method, and will provide some theoretical analysis of the accuracy of this method.

In the previous paper~\cite{paper1}, we have not explained how to calculate the matrix elements of the band-diagonal matrix efficiently in a concrete discussion. In this paper, we will show that we can calculate them from a small number of coefficients which can be derived easily from the recursion formulae among the basis functions only by four arithmetic operations among integers when the coefficients in the rational functions $p_m(x)$ $(m=0,1,\dots ,M)$ are rational complex numbers.  Moreover, we will explain the detailed procedure of 
a quasi-minimization of the ratio between two quadratic forms which is effective
in removing extra solutions of the system of simultaneous
linear equations, while it has been omitted in~\cite{paper1}. 
We will propose and illustrate in detail an integer-type quasi-orthogonalization, which is effective in the above quasi-minimization and requires a relatively small amount of calculations.

The proof of the convergence to true solutions of the numerical solutions obtained by means of this quasi-orthogonalization is one of the main result of this paper. 
Moreover, we will prove the halting of this process. Another main result of this paper is the proof of an inequality which gives an upper bound of errors contained in numerical solutions obtained by the proposed method.

The contents of this paper are as follows: In Section \ref{sec:sv}, we survey the basic structures of this method proposed 
and proved 
in~\cite{paper1}. In section \ref{sec:fb}, the function spaces and the basis systems used in this method are surveyed. In section \ref{sec:me}, we explain how to calculate the `matrix elements' of the matrix representation of a differential operator  using only the coefficients of the differential operator. In Section \ref{sec:cp}, we specify the concrete procedures  for the algorithm proposed only abstractly in~\cite{paper1}, for the cases where the space of 
true solutions in ${\cal H}$ 
%true solutions in $C^M(\mathbb{R} )\cap {\cal H}$ 
is one-dimensional. In Section \ref{sec:suboptimality}, we prove that the realization proposed in Section \ref{sec:cp}  satisfies the conditions required for the convergence of the mixture error to zero. In Section \ref{sec:multidim}, we extend the realization proposed in Section \ref{sec:cp} to general cases where  the space of 
obtainable true solutions in ${\cal H}$ 
%true solutions in $C^M(\mathbb{R} )\cap {\cal H}$ 
is multi-dimensional, and we prove the convergence of the mixture error to zero even for this extension in Section \ref{sec:conv_multidim}. Since the procedures proposed in Section \ref{sec:cp} contain iterations under given conditions, we prove that they 
%terminate 
halt in a finite number of steps  
in Section \ref{sec:halt}. In Section \ref{sec:ul}, we will state a theoretical upper bound for the errors. In Section \ref{sec:nm}, 
we give numerical results which clarifies 
how accurate results our method gives. 
Moreover, in that section, we compare theoretically the order of the amount of required calculations between some typical existing methods and the proposed method. 
Finally, we present our conclusions.

\section{Survey of basic structure}
\label{sec:sv}
As was introduced in ~\cite{paper1}, we propose an integer-type algorithm for finding the $C^M$-solutions in a function space ${\cal H}$ of a higher-order linear ordinary  differential equation $R(x,\frac{d}{dx})f=0$ of Fuchsian type given by a  differential operator $\displaystyle R(x,{\textstyle \frac{d}{dx}})= \sum_{m=0}^M r_m(x) \, ({\textstyle {\textstyle \frac{d}{dx}}})^m$ with rational functions $r_m(x)$ $(m=0,1,...,M)$ of $x$ with rational-(complex-)valued coefficients such that $r_M(\pm i)\ne 0$ and the coefficient functions $r_m(z)$ $(m=1,2,\ldots ,M-1)$ have no singularity at $z=\pm i$. 
Even for non-Fuchsian cases, 
the proposed method can be applied with some restrictions~\cite{paper1}~\cite{paper2}. 

% Even when the condition $ ^\forall x\in \mathbb{R}, \,\, p_M(x)\ne 0$ is not satisfied, 
% if the ODE belongs to Fuchsian class, we can apply this algorithm, only under a slight modification.

% For non-Fuchsian cases, 
% the proposed method can be applied with some restrictions~\cite{paper1}~\cite{paper2}. However, in these cases, the mathematical structure of the method are much more complicated than the cases where this condition is satisfied. Therefore, for simplicity, we limit the survey only to the cases where 
% $p_M(x)$ has no zero points, 
%this condition is satisfied, 
% and later we explain briefly only the difference between such cases and the cases which are not so. 
By multiplying the least common multiple of the denominators $r_m(x)$ 
$(m=0,1,...,M)$, this problem can be simplified, without loss of generality, into the problem to solve the ODE 
$\displaystyle P(x,{\textstyle \frac{d}{dx}})= \sum_{m=0}^M 
p_m(x) \, ({\textstyle {\textstyle \frac{d}{dx}}})^m$ with 
polynomials $p_m(x)$ $(m=0,1,...,M)$ of $x$ with rational-(complex-)valued coefficients such that 
$p_M(\pm i)\ne 0$. 

Moreover, for Fuchsian cases, as is shown in another paper~\cite{paper1}, by multiplying an appropriately chosen polynomial, the above problem can be modified into the problem to solve the ODE 
$\displaystyle Q(x,{\textstyle \frac{d}{dx}})= \sum_{m=0}^M 
q_m(x) \, ({\textstyle {\textstyle \frac{d}{dx}}})^m$ with 
polynomials $q_m(x)$ $(m=0,1,...,M)$ of $x$ with rational-(complex-)valued coefficients such that 
$q_M(\pm i)\ne 0$ and the multiplicities of the zero points of $q_M(x)$ are not smaller than $M$.
When $p_M(x)$ has no zero point, without any conflict with the above condition about the multiplicities of zero points, we may set $Q(x,\frac{d}{dx})=P(x,\frac{d}{dx})$. 

This method is based on the band-diagonal matrix representation of the operator $B_Q$ which is defined by the closure of the operator $\tilde{B}_Q:{\cal H}\to {{{\cal H}^\Diamond}}$ ($\tilde{B}f=Q(x,\frac{d}{dx})f$)  with domain $D(\tilde{B}_Q):=
\{f\in C^M(\mathbb{R} ) \cap {\cal H} | Q(x,{\textstyle \frac{d}{dx}}) f \in {{{\cal H}^\Diamond}} \},$ 
where ${{{\cal H}^\Diamond}}$ is a function space which contains (as a set) ${\cal H}$, but which has a different norm, $\langle\cdot ,\cdot\rangle _{{\cal H}^\Diamond}$, from the norm for ${\cal H}$, $\langle\cdot ,\cdot\rangle_{{\cal H}}$. Under the conditions {\bf C1}-{\bf C4}, {\bf C1$^+$} and {\bf C2$^+$} below, we can show the validity of this band-diagonal matrix representation $b_{m}^n:=\langle B_Q e_n, e^\Diamond_m \rangle _{{\cal H}^\Diamond}$ with respect to orthonormal basis systems $\{e_n\}_{n\in \mathbb{Z}^+} $ and $\{e^\Diamond_n\}_{n\in \mathbb{Z}^+} $ for ${\cal H}$ and ${{{\cal H}^\Diamond}}$, respectively, \,  i.e. a one-to-one correspondence is guaranteed between a solution $f$ in $C^M(\mathbb{R}\setminus S)\cap {\cal H}$ ($S$: set of singular points of the ODE) of the differential equation $Q(x,\frac{d}{dx})f=0$ and a vector $\vec{f}$ in $V\cap \ell ^2(\mathbb{Z}^+)$ with
\begin{eqnarray}
V:= \bigl \{ \vec{f} \, \bigl| \!\!\!\!\!\!\!\!\!\!\sum_{n=\max(0,m-\ell _0 )}^{m+\ell _0} \!\!\!\!\!\!\!\!\!\!  b_m^n f_n = 0 \,\, (m\in\mathbb{Z}^+)\, \bigr\} \, , 
\label{eqn:def_sp_sol_infty}\end{eqnarray}
where $\displaystyle f=\sum_{n=0}^\infty f_ne_n$~\cite{paper1}. The conditions required are: 
\begin{description}
\item[C1]
There exists a CONS $\{e_n \, |\, n\in\mathbb{Z}^+\}$ of ${\cal H}$ such that
$e_n \in D(\tilde{B}_Q)$.
\item[C2]
There exist an integer $\ell _0 $ and a CONS $\{e_n^\Diamond |\, n\in\mathbb{Z}^+\}$ of ${{\cal H}^\Diamond }$ 
such that \par\noindent $b_{m}^n:=\langle B_Q e_n, e_m^\Diamond  \rangle _{{{\cal H}^\Diamond }}=0$ 
%and $c_{m}^n:=\langle e_n, e_m^\Diamond  \rangle _{{{\cal H}^\Diamond }}=0$ 
when $|n-m|> \ell _0 $.
% \,\, (Here note that $e_n \in {\cal H} \subset {{\cal H}^\Diamond }$.)
\item[C3]
%There exists a linear operator $C$ from a densly subset of ${{\cal H}^\Diamond }$ to a densely subset of ${\cal H}$
There exists a linear operator $C_Q$ with domain $D(C_Q)$ from a dense subspace of ${{\cal H}^\Diamond }$ to ${\cal H}$ 
such that $e_m^\Diamond  \in D(C_Q)$ and 
%$\langle B f, e_m^\Diamond  \rangle _{{{\cal H}^\Diamond }}=\langle f, C e_m^\Diamond  \rangle _{{{\cal H}^\Diamond }}$
$\langle B_Q f, e_m^\Diamond  \rangle _{{{\cal H}^\Diamond }}=\langle f, C_Q e_m^\Diamond  \rangle _{{\cal H}}$
for $f \in D(\tilde{B}_Q)$.
\item[C4]
%Any function $f$ in ${\cal H}$ satisfying 
For any sequence $\{f_n\}_{n=0}^{\infty} \in \ell ^2$ satisfying $\displaystyle\sum_{n=\max(0,m-\ell _0 )}^{m+\ell _0} \!\!\!\!\!\!\!\!\!\!  b_m^n f_n = 0 \,\, (m\in\mathbb{Z}^+)$,
the sum $\displaystyle \sum_{n=0}^N f_n e_n$ converges (with respect to the ${\cal H}$-norm) to a solution $f\in C^M(\mathbb{R} )\cap {\cal H}$ of  $Q(x,\frac{d}{dx})f=0$ 
as $N\to\infty$.
\item[C1$^+$]
There exists a positive function $\upsilon$ in $C^M(\mathbb{R} )$ s.t. $\displaystyle \langle f,g \rangle _{\cal H}=\int_{-\infty }^\infty \!\!\! f(x)\overline{g(x)}\upsilon(x)dx$.
\item[C2$^+$]
There exists a positive function ${\upsilon^\Diamond}$ in $C^M(\mathbb{R} )$ s.t. $\displaystyle \langle f,g \rangle _{{\cal H}^\Diamond}
=\int_{-\infty }^\infty \!\!\! f(x)\overline{g(x)}{\upsilon^\Diamond}(x)dx$.
\end{description}

Especially when the differential equation has no singular points, Conditions {\bf C1$^+$} and {\bf C2$^+$} can be exempted, and then we can show the one-to-one correspondence between a solution $f$ in 
$C^M(\mathbb{R})\cap {\cal H}$ of the differential equation $Q(x,\frac{d}{dx})f=0$ (the same as $P(x,\frac{d}{dx})f=0$) and a vector $\vec{f}$ in $V\cap \ell ^2(\mathbb{Z}^+)$ only from Conditions {\bf C1}-{\bf C4}.

From the above discussions, in order to obtain true solutions of the ODE, we should extract only square-summable vector solutions of the system of simultaneous linear equations $\displaystyle\sum_{n=\max(0,m-\ell _0 )}^{m+\ell _0} \!\!\!\!\!\!\!\!\!\!  b_m^n f_n = 0 \,\, (m\in\mathbb{Z}^+)$. This extraction can be made approximately by the method explained below when the following condition {\bf C5} is satisfied. 
\begin{description}
\item[C5]
There exists an integer $j_0\in \mathbb{Z}^+$ such that $b_m^{m+\ell _0 } \ne 0 $
for any integer $m \ge j_0 $ $ (m\in \mathbb{Z} ^+)$.
\end{description}
In this case, the dimension $D$ of $V$ is equal to that of 
\begin{eqnarray}
\Pi_{p_0} V=\bigl \{ \{f_n\}_{n=0}^{p_0} \, \bigl| \, 
\sum_{n=0}^{p_0} b_m^n f_n = 0 \,\, (m=0,1,...,j_0-1)\, \bigr\},\Label{eqn:Haya2}
\end{eqnarray}
where 
\begin{eqnarray}\label{eqn:def_p0}
p_0:=j_0+\ell _0 -1
\end{eqnarray}
and the truncation operator $\Pi_m$ is defined by 
\begin{eqnarray}
(\Pi_m \vec{f})_n = \left\{\begin{array}{@{\,}ll} f_n & \,\,\,\, 
(n\le m) \\ 0 & \,\,\,\, (n>m)\,\, . \end{array}\right.
\Label{eqn:1}
\end{eqnarray}
%The range $\Pi_m$ can be regarded as $\ell^2(\{0,\ldots, m\})$.
In the following, for simplicity, we sometimes identify $\Pi_m\vec{f}$ with the corresponding  $(m+1)$-dimensional vector.  

For the extraction only of square-summable vector solutions, 
we choose a bounded bilinear form $\Omega(\vec{f}, \, \vec{g})$ on 
$\ell ^2(\mathbb{Z}^+)\times \ell ^2(\mathbb{Z}^+)$ 
(and the corresponding quadratic form $\Omega(\vec{f}):=\Omega(\vec{f}, \, \vec{f})$ on 
$\ell ^2(\mathbb{Z}^+)$ ) and the integers $K$ and $N$ 
satisfying
\begin{eqnarray}
& ^\forall \vec{f} \in \ell ^2(\mathbb{Z}^+), \quad 
\Omega(\vec{f})\ge \|\vec{f}\|_{\ell ^2}^2 :=
\displaystyle \sum_{n=0}^{\infty} |f_n|^2,\quad
% \forall \vec{f} \in \ell ^2(\mathbb{Z}^+).
\Label{OC6} \\
& N \ge K \ge j_0+\ell _0 -1,
\Label{OC7}
\end{eqnarray} 
and define the ratio and its minimum:
\begin{eqnarray}
\sigma_{K,N}^{(\Omega )}(\vec{f}) 
&:=&
\frac{\Omega (\Pi_N \vec{f})}
{\|\vec{f}\|_{\ell ^2,K}^2} \quad \hbox{ for } \vec{f}  \in V\setminus\{0\}\\
\underline{\sigma_{K,N}^{(\Omega )}}
&:=&
\min_{\vec{f} \in V\setminus\{0\}} \sigma_{K,N}^{(\Omega )}(\vec{f}) .
\Label{eqn:def_ratio}
\end{eqnarray}
% This definition is well defined because 
% Conditions {\bf C2}, {\bf C5} and (\ref{OC7}) guarantee the relation
% \begin{eqnarray}
% \| \vec{f}\|_{\ell ^2,K}:=\|\Pi_K\vec{f}\|_{\ell ^2}
% >0
% \quad \hbox{ for } \vec{f}  \in V\setminus\{0\}.
% \Label{lemma:K_norm}
% \end{eqnarray}
Similarly, 
we define 
\begin{eqnarray}
\sigma_{K,\infty}^{(\Omega )}(\vec{f}) 
&:=& 
\frac{\Omega (\vec{f})}
{\|\vec{f}\|_{\ell ^2,K}^2} \quad \hbox{ for } \vec{f}  \in V\setminus\{0\}\\
\underline{\sigma_{K,\infty}^{(\Omega )}}
&:=&
\min_{\vec{f} \in V\setminus\{0\}} \sigma_{K,\infty}^{(\Omega )}(\vec{f}) .
\end{eqnarray}
% Hence, 
% our solution space 
% $V_{K}$ is given with the following condition
The proposed method yields all the vectors in a linear space $V_{K}$ satisfying the following condition:
\begin{eqnarray}
V_{K}\subset 
\Pi_{K}
\left((\sigma_{K,N}^{(\Omega )})^{-1}[0,c \underline{\sigma_{K,N}^{(\Omega )}} ]\right) 
\cup \{0\}
\Label{4-3-1}
\end{eqnarray}

The proposed method is materialized by means of a practical algorithm or finding a basis system of $V_K$ satisfying (\ref{4-3-1}). 
% The problem to find  $V_K$ satisfying (\ref{4-3-1}) is a kind of quasi-minization problem of the ratio between two quadratic forms. An usual easy method for this type of problem is based on the calculation of an eigenvector of the matrix $A^{-\frac{1}{2}}BA^{-\frac{1}{2}}$ associated with the minimum eigenvalue with the matrices $A$ and $B$ defined by 
% $(A)_{i\, j}:=(\Pi_K\vec{F}^{(i)},\, \Pi_K\vec{F}^{(j)})_{\ell ^2}$ and $(B)_{i\, j}:=\Omega(\Pi_N\vec{F}^{(i)}, \Pi_N\vec{F}^{(j)})$. However, this normal method is quite difficult to apply, 
% because the matrices $A$ and $B$ are usually very close to a singular matrix with rank $1$ due to the most diverging components in $V$, 
% and hence this usual method is particularly subject to 
% the `canceling' due to round-off errors.  
%In order to overcome  this difficulty, 
This algorithm is based on the intermediate idea between the Gram-Schmidt orthogonalization and the Euclidean algorithm, which requires a relatively small amount of calculations. 

When the purpose is to
 calculate the truncated elements $\ell ^2$ solution $\Pi_K(V \cap \ell ^2(\mathbb{Z}^+))$,
the error is evaluated by the norm concerning inner product 
$\langle \vec{x},\vec{y}\rangle_{\ell ^2,K}:=
\langle \Pi_K\vec{x},\Pi_K\vec{y}\rangle_{\ell ^2}$.
Denoting the 
the projection to 
$W$ concerning this inner product by $P_{W,K}$,
we can evaluate the accuracy of our result $V_{K}$ for the worst case by
% \begin{eqnarray}
% \sup_{\vec{x} \in V_{K}\setminus\{0\}}
% \frac{ \|P_{V\cap \ell ^2(\mathbb{Z}^+),K} \vec{x} - \vec{x}\|_{\ell ^2,K} }{\|\vec{x}\|_{\ell ^2,K} }.
% \end{eqnarray}
% However, 
% the subspace $V_{K}$ is not uniquely defined, and is chosen with the condition
% (\ref{4-3-1}).
% The accuracy of our algorithm should be evaluated with the worst case as follows:
\begin{eqnarray}
&& \sup_{V_{K}\subset 
\Pi_{K}
((\sigma_{K,N}^{(\Omega )})^{-1}[0,c \underline{\sigma_{K,N}^{(\Omega )}} ]) \cup \{0\}} \,\, 
\sup_{\vec{x} \in V_{K}\setminus\{0\}}
\frac{ \|P_{V,K} \vec{x} -\vec{x}\|_{\ell ^2,K} }{\|\vec{x}\|_{\ell ^2,K} }.
\end{eqnarray}
It can be shown that this value goes to $0$ and all the true solutions of the ODE can be approximated by our solution space, from the following theorems~\cite{paper1}. 

\begin{theorem}{\rm (Theorem 2.5 of~\cite{paper1})}\Label{thm:conv1}
For fixed $K$, when $N$ goes to infinity,
the convergence
\begin{eqnarray}
\sup_{
\vec{x} \in 
(\sigma_{K,N}^{(\Omega )})^{-1} [0,c \underline{\sigma_{K,N}^{(\Omega )}} ] }
\frac{ \|P_{V,K} \vec{x} -\vec{x}\|_{\ell ^2,K} }{\|\vec{x}\|_{\ell ^2,K} }
\to 0
\end{eqnarray}
holds.
\end{theorem}

\begin{theorem}{\rm (Theorem 2.6 of~\cite{paper1})}\Label{thm:conv2}
Assume that $\Omega(\vec{x})=\|\vec{x}\|_{\ell ^2}$.
When we choose sufficiently large numbers $N_0$ and $c_0$,
then for any $N \ge N_0$ and $c\ge c_0$, we have 
\begin{eqnarray}
P_{V,K}V_K = \Pi_K( V \cap \ell ^2(\mathbb{Z}^+))\Label{eq:2}
\end{eqnarray}
for any choice of $V_K$.
\end{theorem}
Their proofs have been given in Section 4 of~\cite{paper1}. 
As a practical choice of the quadratic form $\Omega(\vec{f},\vec{g})$, we can use $\displaystyle \Omega(\vec{f},\vec{g})=\sum_{n=0}^\infty  f_n \,\overline{g_n}\, w_n$
with a non-decreasing `weight number sequence' satisfying $w_n=1$ for $n\le K$ and $w_n=R$ for $n\ge J$ with  integers $K$ and $J$ such that $j+\ell _0 -1\le K\le J\le N$. To enable discussion of the upper error bound, we limit the bilinear form to this class. 
In particular, the weight number sequence $w_n$ used in numerical experiments of this study is 
\begin{eqnarray} \label{eqn:w_n} 
w_n :=\left\{\begin{array}{@{\,}ll} 1 & (n\le K) \\ e^{r(\mu_n-\mu_K)} & \!\!\! (K<n<J) \\ R:=e^{r(\mu_J-\mu_K)} & (n\ge N) \end{array}\right. 
 \\   \mbox{ with } \,\,\,\,\, 
\mu_n:=\left|\ddot{n}_{k_0,n}-\frac{k_0+1}{2}\right |-\frac{k_0+1}{2}. 
\nonumber \end{eqnarray}
Empirically, the choice with $r=10^8$, $K=2\lfloor \frac{3N}{8} \rfloor + k_0$ and 
%$J=2\lfloor \frac{7N}{16} \rfloor + k$ 
$J=2\lfloor \frac{7N}{16} \rfloor + k_0$ or 
$K=2\lfloor \frac{7N}{16} \rfloor + k_0$ and 
$J=2\lfloor \frac{15N}{32} \rfloor + k_0$   
often gives good results, for example. This choice preserves the symmetry property of the basis system introduced below. The specification of the class of bilinear forms given in this paragraph is not related to the following sections, except for Sections \ref{sec:ul} and \ref{sec:nm}.

\section{Function spaces and basis systems used in the algorithm}
\label{sec:fb}
In this section, we will introduce the function spaces ${\cal H}$, ${\cal H}^\Diamond $ and their basis systems \par\noindent $\{e_n |\, n\in\mathbb{Z}^+\}$, $\{e_n^\Diamond |\, n\in\mathbb{Z}^+\}$ which satisfy the conditions {\bf C1}-{\bf C4} and {\bf C5}.

First, define the inner product and norm parametrized by $k\in\mathbb{Z}$ as 
\begin{eqnarray*}
\langle f, \,g\rangle_{(k)} := \int_{-\infty}^\infty f(x) \,\overline{g(x)}\ \,(x^2 +1)^k \,dx \, \,\,\,\,\,\, {\rm and} \,\,\,\,\,\, 
\| f\|_{(k)} := \int_{-\infty}^\infty |f(x)|^2  \,(x^2 +1)^k \,dx \, . 
\end{eqnarray*}
Now we can introduce the Hilbert space of functions with inner product $\langle f, \,g\rangle_{(k)}$
\begin{eqnarray*}
L_{(k)}^2 (\mathbb{R}) :=\bigl\{ \, f : {\rm measurable }\, \bigl| \, \|f\|_{(k)}  < \infty \bigr\}.  
\end{eqnarray*}
Then, $L_{(k)}^2 (\mathbb{R}) \subset L_{({\kappa })}^2 (\mathbb{R} )$ if $k\ge {\kappa }$.  
Note that $L_{(0)}^2 (\mathbb{R}) = L^2 (\mathbb{R}) $. 
For the spaces ${\cal H}$ and ${{{\cal H}^\Diamond}}$, 
we choose 
\begin{eqnarray}
{\cal H} = L_{k_0}^2(\mathbb{R} ) ,\,\,\, 
{{{\cal H}^\Diamond}} = L_{{k_0^\Diamond}}^2(\mathbb{R} ) 
\label{eqn:def_sp}\end{eqnarray}
with integers $k_0$ and ${k_0^\Diamond}$ satisfying $k_0\ge 0$ and $\displaystyle {k_0^\Diamond}\le k_0-\max_{m\in \{0,1,\ldots ,M\}} (\deg p_m-m)$.

Next, we will introduce the basis function systems.
Define the wavepacket functions 
\begin{eqnarray}
\psi_{k,\, \ddot{n}}(x) := \frac{1}{(x+i)^{k+1}} \left( \frac{x-i}{x+i} \right)^{\ddot{n}} \,\,\,\,\,\, (k,\ddot{n}\in \mathbb{Z}) . 
\label{eqn:def_psi}\end{eqnarray}
The indices of functions in $\left\{\,\psi_{k_0,\, \ddot{n}}\, \bigl| \, \ddot{n}\in\mathbb{Z} \right\}$ are bilaterally expressed, whereas the indices of basis functions in  $\{ e_n \}_{n=0}^{\infty}$ are unilaterally expressed, and they are `matched' to one another by the one-to-one mapping defined by ${\ddot{n}}_{k,n}$ in (\ref{eqn:def_bs}) below. In order to avoid confusion between them, in this paper, the integer indices with double dots\, $\ddot{}$\, denote the bilateral ones in $\mathbb{Z}$, in contrast with the unilateral ones (without double dots) in $\mathbb{Z}^+$. 
The functions introduced in the above satisfy the symmetry property   $\overline{\psi_{k_0,n}}=\psi_{k_0,-n-k_0-1}$. 

These are sinusoidal-like wavepackets with spindle-shaped envelopes. An example of the shapes of these wavepackets is illustrated in 
Figures \ref{fig:1} and \ref{fig:2}.  
As is explained in Appendix \ref{app:shape} in detail, the wavepackets defined by (\ref{eqn:def_psi}) are `almost-sinusoidally' oscillating wavepackets with spindle-shaped envelopes $|\psi_{k,n}(x)| = (x^2 +1)^{-\frac{k+1}{2}} \,$,\ when $k\ge 0$, and their approximation to a sinusoidal wavepacket with Gaussian envelope 
\begin{eqnarray*}
\psi_{k,n} {\textstyle \left(\frac{x}{\sqrt{k}}\right)} \approx ({\rm const.}) \, \cdot \,\, e^{ i\frac{k+2n+1}{\sqrt{k}} \, x} \,\,\, e^{-\frac{1}{2}x^2} \,\,\, .
\end{eqnarray*}
holds for sufficiently large $k$, in the sense that we can show the convergence \par\noindent 
$\displaystyle \, \lim_{k\to\infty } \parallel \, \Xi_{k,n} \,\, e^{- i\frac{k+2n+1}{\sqrt{k}} \, x} \, \psi_{k,n}{\textstyle \left(\frac{x}{\sqrt{k}}\right)} \, - \, e^{-\frac{1}{2}x^2} \parallel \, = 0 \,\,\,\,\,  (\, \Xi_{k,n}:\,{\rm const.})\, $ with respect to  the $L^2$-norm. This property shows the suitability of the wavepackets for expanding 
almost-localized smooth solutions. 

Moreover, as is shown in our preceding paper~\cite{paper2}, they are related to the basis functions of Fourier series, under a change of variable $x\to \theta:=2\arctan x$ used also in another field ~\cite{Qia}~\cite{Che}. 

\begin{figure}[t]
\begin{center}
\begin{minipage}{0.45\linewidth}
\includegraphics{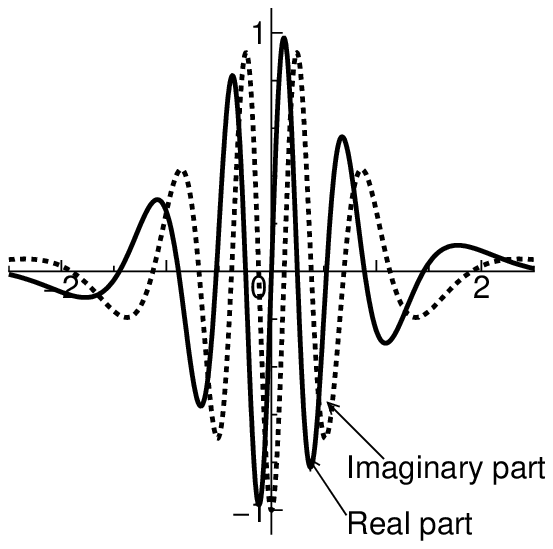}
\caption{ Graph of $\psi_{2,5}(x)$} 
\label{fig:1}
\end{minipage}
\begin{minipage}{0.45\linewidth}
\includegraphics{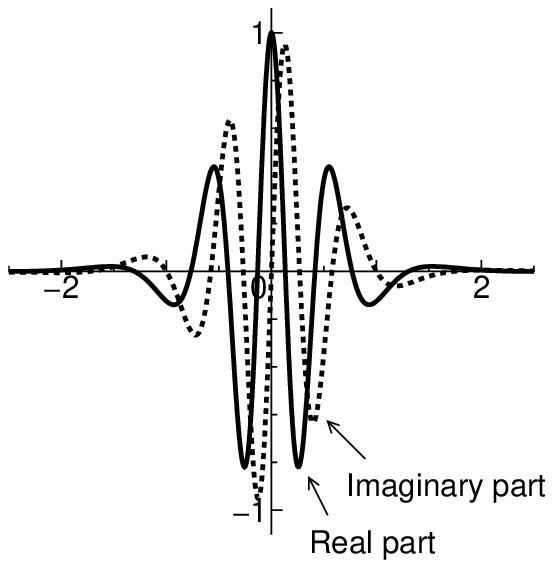}
\caption{ Graph of $\psi_{5,3}(x)$} 
\label{fig:2}
\end{minipage}
\end{center}
\end{figure}

These functions are used for the orthonormal basis systems $\{e_n\}_{n=0}^\infty $ and $\{e^\Diamond_n\}_{n=0}^\infty $ of ${\cal H}$ and ${{{\cal H}^\Diamond}}$, respectively, as follows: 

Define 
\begin{eqnarray}\label{eqn:def_bs}
e_n= \sqrt{\textstyle\frac{1}{\pi}} \,\psi_{k_0,\, \ddot{n}_{k_0,n}} , 
\,\,\, 
e^\Diamond_n= \sqrt{\textstyle\frac{1}{\pi}} \,\psi_{{k_0^\Diamond},\, \ddot{n}_{{k_0^\Diamond},n}} 
\\ \nonumber 
{\rm with}\,\,\,\,\,\,  \ddot{n}_{k,n} := \left\lfloor {\textstyle -\frac{k+1}{2}} \right\rfloor + (-1)^{n+k+1}\, \left\lfloor {\textstyle \frac{n+1}{2}} \right\rfloor \,\, .
\end{eqnarray}
Under the choices (\ref{eqn:def_sp}) and (\ref{eqn:def_bs}), we can show that the assumptions {\bf C1}-{\bf C4} are satisfied~\cite{paper1}. Moreover, if $q_M(\pm i)\ne 0$, the assumption {\bf C5} is also satisfied~\cite{paper1}. We can show that 
\begin{eqnarray}
\ell _0=2M+k_0-k_0^\Diamond \,\,\,\,\,\, {\rm and } \,\,\,\,\,\, j_0=\max({k_0^\Diamond},0)
\label{eqn:lj}\end{eqnarray}
for Conditions {\bf C2} and  {\bf C5}~\cite{paper1}~\cite{paper2}. In addition, if all the coefficients of the polynomials $q_m(x)$ $(0\le m\le M)$ belong to $\mathbb{Q} +\mathbb{Q} i$, then {\bf C7} is always satisfied~\cite{paper1}. 

Here we point out some properties of $\psi_{k,\, {\ddot{n}}}$ defined in (\ref{eqn:def_psi}), which will be important later.
\begin{theorem}
\label{thm:ixd_psi}
For any integer $n$, 
\begin{eqnarray}
\psi_{k,\, {\ddot{n}}}(x) = - \frac{i}{2} \left( \psi_{k-1,\, {\ddot{n}}} (x) - \psi_{k-1,\, {\ddot{n}}+1} (x) \right) \,\,\,\,\,\,\,\, %(k\ge 1) 
\label{eqn:id_CSI}\end{eqnarray}
\begin{eqnarray} 
x\, \psi_{k,\, {\ddot{n}}} (x) = \frac{1}{2} \left( \psi_{k-1,\, {\ddot{n}}} (x) + \psi_{k-1,\, {\ddot{n}}+1} (x) \right)\,\,\,\,\,\,\,\, %(k\ge 1) 
\label{eqn:mult_CSI}\end{eqnarray}
\begin{eqnarray}
{\textstyle \frac{d}{dx}}\, \psi_{k ,\, {\ddot{n}}} (x)   
= {\ddot{n}} \,  \psi_{k+1 ,\, {\ddot{n}}-1} (x) - ({\ddot{n}}+k+1) \, \psi_{k+1 ,\, {\ddot{n}}} (x) \,\, . %\,\,\,\,\,\, (k\ge 0) \,\,\,\, . 
\label{eqn:diff_CSI}\end{eqnarray}
\end{theorem} \par\noindent 
This theorem can be derived directly from the definition of $\psi_{\ddot{n},k}$. It is essential for the conditions {\bf C2} and {\bf C5} as is shown in ~\cite{paper1}, and it will be useful for showing how to calculate the matrix elements $b_m^n$ by integer-type programs in Section \ref{sec:me}.

\section{Matrix elements calculated using only the coefficients of the polynomials in the differential equation}
\label{sec:me}
In order to find solutions of the simultaneous linear equations $\sum_n b_m^n f_n=0$ $(m\in\mathbb{Z}^+ )$, we have to determine the `matrix elements' $b_m^n$ $(m,n\in \mathbb{Z}^+; \, |m-n|\le \ell _0 )$ from the differential operator $Q(x,\frac{d}{dx})$.  
In this section, we explain how we can calculate these matrix elements from just the coefficients $q_{m,j}$ $(m=0,1, \ldots M;\, j=0,1,\ldots \deg q_m )$ of the polynomials $\displaystyle q_m(x)=\sum_{j=0}^{\deg q_m} q_{m,j}x^j$ in the differential operator $\displaystyle \sum_{m=0}^M q_m(x)({\textstyle \frac{d}{dx}})^m$ by means of a recursive use of the properties (\ref{eqn:id_CSI}), (\ref{eqn:mult_CSI}) and (\ref{eqn:diff_CSI}) of $\psi_{k,\ddot{n}}$.

In ~\cite{paper1} and in Section \ref{sec:sv} of this paper, the matrix element $b_m^n$ is defined in a general framework only as the inner product $\langle B_Qe_n, e^\Diamond_m \rangle _{\widetilde{\cal_H}}$. However, because of these properties of $\psi_{k,\ddot{n}}$, we can determine $b_m^n$ from just the coefficients $q_{m,j}$ $(m=0,1, \ldots M;\, j=0,1,\ldots \deg q_m)$, without any calculations of inner products. 

In order to explain this, we introduce the matrix elements $\ddot{b}_{\ddot{m}}^{\ddot{n}}$ $(\ddot{m},\ddot{n}\in \mathbb{Z} )$ in the `bilateral expression'. Define 
\begin{eqnarray}
 \ddot{b}_{\ddot{m}}^{\ddot{n}}:=\frac{1}{\pi }\, \langle Q(x,{\textstyle \frac{d}{dx}})\psi_{k_0,\ddot{n}}, \psi_{k_0^\Diamond ,\ddot{m}}\rangle _{{\cal H}^\Diamond } = b_{\ddot{n}_{k_0,m}}^{\ddot{n}_{k_0,n}}
\label{eqn:def_ddorb}\end{eqnarray}
where the last equality is derived from the definitions in (\ref{eqn:def_bs}) together with the `matching' between the unilateral and bilateral expressions.

The function $Q(x,{\textstyle \frac{d}{dx}})\psi_{k_0,\ddot{n}}$ can be expressed as a linear combination of $\psi_{{k_0^\Diamond},\ddot{m}}$ with $\ddot{m}=\ddot{n}-M,\,\ddot{n}-M+1,\ldots ,\ddot{n}+M+k_0-k_0^\Diamond $, by the linear combination of the results of the recursive use of (\ref{eqn:diff_CSI}) $m$ times, that of (\ref{eqn:mult_CSI}) $j$ times and that of (\ref{eqn:id_CSI}) $k_0-{k_0^\Diamond}+m-j$ times. 
From this fact, it is easily shown that the matrix elements  $\ddot{b}_{\ddot{m}}^{\ddot{n}}$ $(\ddot{m},\ddot{n}\in \mathbb{Z} ; \, \ddot{n}-M\le \ddot{m}\le \ddot{n}+M+k_0-k_0^\Diamond )$ can be written as a polynomial in $\ddot{n}$ whose order is not greater than $M$, under fixed $k$. For $\ddot{m}\le \ddot{n}-M-1$ or $\ddot{m}\ge \ddot{n}+M+k_0-{k_0^\Diamond}+1$,\,\,  $\ddot{b}_{\ddot{m}}^{\ddot{n}}=0$.

Since these polynomials are similar to one another for the matrix elements with common difference $\ddot{m}-\ddot{n}$, under the representation of  $\ddot{b}_{\ddot{m}}^{\ddot{n}}$ by a function of $\ddot{m}-\ddot{n}$ and $\ddot{n}$, the power series expansion 
\begin{eqnarray}
\ddot{b}_{\ddot{m}}^{\ddot{n}} = \sum_{{s}=0}^M \,\ddot{\beta }_{\ddot{m}-\ddot{n},\,{s}} \,\, \ddot{n}^{s}
\label{eqn:expansion_matrix_element}\end{eqnarray}
with coefficients $\ddot{\beta }_{\ddot{r},\,{s}} $ $(s=0,1,\dots ,M)$ is convenient. Here note that  $\ddot{\beta }_{\ddot{r},\,{s}}=0 $ for $\ddot{r}\le -M-1$ or $\ddot{r}\ge M+k_0-k_0^\Diamond+1$. Hence, we can calculate all the matrix elements $\ddot{b}_{\ddot{m}}^{\ddot{n}}$ by the expansion (\ref{eqn:expansion_matrix_element}) from just the  $(2M+k_0-k_0^\Diamond+1)(M+1)$ coefficients $\ddot{\beta }_{\ddot{r},\,{s}} $ with  $\ddot{r}=-M,-M+1,\ldots ,M+k_0-k_0^\Diamond$ and $s=0,1,\dots ,M$.

By these relations, the calculations of the coefficients $\ddot{\beta }_{\ddot{r},\,{s}} $ can be performed by the procedures in Table \ref{tbl:1}, where the modules {\bf (I), (X)} and {\bf (D)} just correspond to (\ref{eqn:exp_id}), (\ref{eqn:exp_mult}) and (\ref{eqn:exp_diff}), respectively.

It is easily shown that all the coefficients $\ddot{\beta }_{\ddot{r},\,{s}} $ are rational-(complex-)valued when all the coefficients $q_{m,j}$ of the polynomial in the differential operator are. Hence, if all the coefficients $q_{m,j}$ are rational-(complex-)valued,  so are all the matrix elements $\ddot{b}_{\ddot{m}}^{\ddot{n}}$ (and hence $b_m^n$). This fact shows that the condition {\bf C7} is satisfied if all the coefficients $q_{m,j}$ are rational-(complex-)valued.

In terms of the unilateral expression, the matrix elements $b_m^n$ can be calculated by 
\begin{eqnarray}  \label{eqn:matrix_elements}
b_m^n  &=& \ddot{b}_{\ddot{n}_{k_0,m}}^{\ddot{m}_{k_0,n}}=  \sum_{{s}=0}^M \ddot{\beta }_{(\ddot{n}_{k_0,m}-\ddot{n}_{k_0,n}),\, {s}} \,\, (\ddot{m}_{k_0,n})^{\, {s}} \\ 
&=& \left\{\begin{array}{@{\,}ll} \displaystyle 
\sum_{{s}=0}^M\ddot{\beta }_{\bigl((-1)^{m+k+1}(m-n)/2\bigr),\,{s}} \,\, (\ddot{n}_{k_0,n})^{\, {s}} & (\mbox {if }\, m+n\mbox{ is even}) \\ \\ \displaystyle 
\sum_{{s}=0}^M\ddot{\beta }_{\bigl((-1)^{m+k+1}(m+n+1)/2\bigr),\,{s}}\,\, (\ddot{n}_{k_0,n})^{\, {s}} & (\mbox {if }\, m+n\mbox{ is odd}) \, , 
\end{array}\right.
\nonumber \end{eqnarray}
from just the coefficients $\ddot{\beta }_{\ddot{r},\,{s}} $ with  $\ddot{r}=-M,-M+1,\ldots ,M+k_0-k_0^\Diamond$ and \par\noindent $s=0,1,\dots ,M$. 
Here note that a matrix element $b_m^n$ with $n$ and $m$ such that $m+n$ is an odd integer greater than $\ell _0 +k_0-k_0^\Diamond$ $(=2M+2k_0-2{k_0^\Diamond})$ vanishes even when $|m-n|\le \ell _0 $ because $\ddot{\beta }_{\ddot{r},s}=0$ for $\ddot{r}\le -M-1$ or $\ddot{r}\ge M+k_0-k_0^\Diamond+1$. This fact implies the `sparseness' of the `band' in the band-diagonal representation of $B_Q$. Moreover, this fact implies that the calculations of the matrix elements $b_m^n$ can be carried out using only integer-type programs.

The unilateral expression is somewhat less convenient for practical programs than the bilateral one, and we use the bilateral expression in actual calculations. However, in this paper, we follow the unilateral expression for consistency with the mathematical framework introduced in~\cite{paper1}.

\begin{table}[p]
\caption{Calculation of $\ddot{\beta }_{{\ddot{r}},{s}}$ for $-M\le \ddot{r}\le M+k_0-k_0^\Diamond$\,\,\,\,\,\, $\displaystyle (\, \sum_{{s}=0}^M \,\ddot{\beta }_{{\ddot{r}},{s}} \, {\ddot{n}}^{s}\, = \ddot{b}_{{\ddot{n}}+{\ddot{r}}}^{\ddot{n}} ) $}
\begin{center}
\scalebox{1.}{
\begin{tabular}{|l|}
\hline 
$\begin{array}{@{\,}ll}  \\ 
\,\,\, \mbox{{\bf for} ${\ddot{r}}=-M\, $ to $\, M+k_0-k_0^\Diamond$} \\ 
\,\,\, \left\lfloor\begin{array}{@{\,}ll} 
\,\,\, \mbox{{\bf for} ${s}=0\, $ to $\, M$} \\  
\,\,\, \left\lfloor\begin{array}{@{\,}ll} 
\,\,\, \mbox{$\ddot{\beta }_{{\ddot{r}},{s}}\leftarrow 0$} \\ 
\end{array}\right. \\  
\end{array}\right. \\ 
\,\,\, \mbox{{\bf for} $m=0$ to $M$} \\  
\,\,\, \left\lfloor\begin{array}{@{\,}ll} 
\,\,\, \mbox{{\bf for} $j=0\,$ to $\,\deg q_m$} \\ 
\,\,\, \left\lfloor\begin{array}{@{\,}ll} 
\,\,\, \mbox{$\kappa \leftarrow k$} \\ 
\,\,\, \mbox{{\bf for} ${\ddot{r}}=-m\, $ to $\, m+k_0-k_0^\Diamond$} \\  
\,\,\, \left\lfloor\begin{array}{@{\,}ll} 
\,\,\, \mbox{{\bf for} ${s}=0\, $ to $\, M$} \\  
\,\,\, \left\lfloor\begin{array}{@{\,}ll} 
\,\,\, \mbox{$\alpha _{{\ddot{r}},{s}} \leftarrow \delta _{{\ddot{r}}0}\,\delta_{{s}0}$} \\ 
\end{array}\right. \\  
\end{array}\right. \\  
\,\,\, \mbox{{\bf Iterate (D)} below $m$ times.} \\ 
\,\,\, \left\lfloor\begin{array}{@{\,}ll} 
\,\,\, \mbox{\bf (D)}  \left\{\begin{array}{@{\,}ll} 
\,\,\, \mbox{{\bf for} ${\ddot{r}}=-m\, $ to $\, m+k_0-k_0^\Diamond-1$} \\  
\,\,\, \left\lfloor\begin{array}{@{\,}ll} 
\,\,\, \mbox{{\bf for} ${s}=0\, $ to $\, m-1$} \\  
\,\,\, \left\lfloor\begin{array}{@{\,}ll} 
\,\,\, \mbox{$\widehat{\alpha }_{{\ddot{r}},\, {s}+1} \leftarrow - \alpha _{{\ddot{r}},{s}} + \alpha _{{\ddot{r}}+1,\, {s}} $} \\ 
\,\,\, \mbox{$\widehat{\alpha }_{{\ddot{r}},{s}} \leftarrow \widehat{\alpha }_{{\ddot{r}},{s}} - (\ddot{r}+\kappa +1)\alpha _{{\ddot{r}},{s}} + \ddot{r}\,\alpha _{{\ddot{r}}+1,\, {s}} $} \\ 
\end{array}\right. \\  
\end{array}\right. \\  
\,\,\, \mbox{{\bf for} ${\ddot{r}}=-m\, $ to $\, m+k_0-k_0^\Diamond-1$} \\  
\,\,\, \left\lfloor\begin{array}{@{\,}ll} 
\,\,\, \mbox{{\bf for} ${s}=0\, $ to $\, m$} \\  
\,\,\, \left\lfloor\begin{array}{@{\,}ll} 
\,\,\, \mbox{$\alpha _{{\ddot{r}},{s}}\leftarrow \widehat{\alpha }_{{\ddot{r}},{s}} $} \\ 
\end{array}\right. \\ 
\end{array}\right. \\ 
\,\,\, \mbox{$\kappa \leftarrow \kappa +1$} \\ 
\end{array}\right. \\ 
\end{array}\right. \\ 
\,\,\, \mbox{{\bf Iterate (X)} below $j$ times.} \\ 
\,\,\, \left\lfloor\begin{array}{@{\,}ll} 
\,\,\, \mbox{\bf (X)}  \left\{\begin{array}{@{\,}ll} 
\,\,\, \mbox{{\bf for} ${\ddot{r}}=-m+1\, $ to $\, m+k_0-k_0^\Diamond$} \\  
\,\,\, \left\lfloor\begin{array}{@{\,}ll} 
\,\,\, \mbox{{\bf for} ${s}=0\, $ to $\, m$} \\  
\,\,\, \left\lfloor\begin{array}{@{\,}ll} 
\,\,\, \mbox{$\widehat{\alpha }_{{\ddot{r}}, {s}} \leftarrow \frac{1}{2}\alpha _{{\ddot{r}}-1, {s}}+ \frac{1}{2}\alpha _{{\ddot{r}}, {s}}$} \\ 
\end{array}\right. \\ 
\end{array}\right. \\ 
\,\,\, \mbox{{\bf for} ${\ddot{r}}=-m+1\, $ to $\, m+k_0-k_0^\Diamond$} \\  
\,\,\, \left\lfloor\begin{array}{@{\,}ll} 
\,\,\, \mbox{{\bf for} ${s}=0\, $ to $\, m$} \\  
\,\,\, \left\lfloor\begin{array}{@{\,}ll} 
\,\,\, \mbox{$\alpha _{{\ddot{r}},{s}}\leftarrow \widehat{\alpha }_{{\ddot{r}},{s}} $} \\ 
\end{array}\right. \\ 
\end{array}\right. \\ 
%\,\,\, \mbox{$\kappa \leftarrow \kappa -1$} \\ 
\end{array}\right. \\ 
\end{array}\right. \\ 
\,\,\, \mbox{{\bf Iterate (I)} below $k_0-k_0^\Diamond-j+m$ times.} \\ 
\,\,\, \left\lfloor\begin{array}{@{\,}ll} 
\,\,\, \mbox{\bf (I)}\,\,\,   \left\{\begin{array}{@{\,}ll} 
\,\,\, \mbox{{\bf for} ${\ddot{r}}=-m+1\, $ to $\, m+k_0-k_0^\Diamond$}  \\  
\,\,\, \left\lfloor\begin{array}{@{\,}ll} 
\,\,\, \mbox{{\bf for} ${s}=0\, $ to $\, m$} \\  
\,\,\, \left\lfloor\begin{array}{@{\,}ll} 
\,\,\, \mbox{$\widehat{\alpha }_{{\ddot{r}},{s}} \leftarrow \frac{i}{2} \alpha _{{\ddot{r}}-1,{s}} - \frac{i}{2} \alpha _{{\ddot{r}},{s}}$} \\ 
\end{array}\right. \\ 
\end{array}\right. \\ 
\,\,\, \mbox{{\bf for} ${\ddot{r}}=-m+1\, $ to $\, m+k_0-k_0^\Diamond$}  \\  
\,\,\, \left\lfloor\begin{array}{@{\,}ll} 
\,\,\, \mbox{{\bf for} ${s}=0\, $ to $\, m$} \\  
\,\,\, \left\lfloor\begin{array}{@{\,}ll} 
\,\,\, \mbox{$\alpha _{{\ddot{r}},{s}}\leftarrow \widehat{\alpha }_{{\ddot{r}},{s}} $} \\ 
\end{array}\right. \\ 
\end{array}\right. \\ 
%\,\,\, \mbox{$\kappa \leftarrow \kappa -1$} \\ 
\end{array}\right. \\  
\end{array}\right. \\ 
\,\,\, \mbox{{\bf for} ${\ddot{r}}=-k_0+k_0^\Diamond-m\, $ to $\, m$} \\  
\,\,\, \left\lfloor\begin{array}{@{\,}ll} 
\,\,\, \mbox{{\bf for} ${s}=0\, $ to $\, M$} \\  
\,\,\, \left\lfloor\begin{array}{@{\,}ll} 
\,\,\, \mbox{$\ddot{\beta }_{{\ddot{r}},{s}} \leftarrow \ddot{\beta }_{{\ddot{r}},{s}} + q_{m,j}\,\alpha _{{\ddot{r}},{s}}$} \\
\end{array}\right. \\ 
\end{array}\right. \\  
\end{array}\right. \\  
\end{array}\right. \\  \\ 
\end{array} $ \\ 
\hline
\end{tabular}
}
\end{center}
\label{tbl:1}
\end{table}

The coefficients $\ddot{\beta }_{\ddot{r},\,{s}} $ with  $\ddot{r}=-M,-M+1,\ldots ,M+k_0-k_0^\Diamond$ and $s=0,1,\dots ,M$ can be calculated by the recursive use of the relations (\ref{eqn:id_CSI}), (\ref{eqn:mult_CSI}) and (\ref{eqn:diff_CSI}), which can be realized by the procedures given in Table \ref{tbl:1} which are programmable on computers, as follows: From these relations, the 
%innovation 
renewal  
of the coefficients $\ddot{\alpha }_{\ddot{r},\,{s}}$ in Table \ref{tbl:1} can be carried out by means of the relations 
\begin{eqnarray}
\sum_{\ddot{r}} \sum_{{s}} \,\ddot{\alpha }_{\ddot{r},\,{s}} \,\, \ddot{n}^{s} \psi_{k,\ddot{n}+\ddot{r}}(x)=
\sum_{\ddot{r}} \sum_{{s}} \,\left(\frac{i}{2}\ddot{\alpha }_{\ddot{r}-1,\,{s}} -\frac{i}{2} \ddot{\alpha }_{\ddot{r},\,{s}} \right) \,\, \ddot{n}^{s} \psi_{k-1,\ddot{n}+\ddot{r}}(x)
\label{eqn:exp_id}\end{eqnarray}
\begin{eqnarray}
\hspace{.8cm}  x \left(\sum_{\ddot{r}} \sum_{{s}} \,\ddot{\alpha }_{\ddot{r},\,{s}} \,\, \ddot{n}^{s} \psi_{k,\ddot{n}+\ddot{r}}(x)\right)=
\sum_{\ddot{r}} \sum_{{s}} \,\left(\frac{1}{2}\ddot{\alpha }_{\ddot{r}-1,\,{s}} +\frac{1}{2} \ddot{\alpha }_{\ddot{r},\,{s}} \right) \,\, \ddot{n}^{s} \psi_{k-1,\ddot{n}+\ddot{r}}(x)
\label{eqn:exp_mult}\end{eqnarray}
\begin{eqnarray} \label{eqn:exp_diff} 
\frac{d}{dx} \left(\sum_{\ddot{r}} \sum_{{s}} \,\ddot{\alpha }_{\ddot{r},\,{s}} \,\, \ddot{n}^{s} \psi_{k,\ddot{n}+\ddot{r}}(x)\right) \hspace{7cm} \\  = 
\sum_{\ddot{r}} \sum_{{s}} \,\Bigr(-(\ddot{r}+k+1)\ddot{\alpha }_{\ddot{r},\,{s}} + (\ddot{r}+1)\ddot{\alpha }_{\ddot{r}+1,\,{s}} -\ddot{\alpha }_{\ddot{r},\,{s}-1}+ \ddot{\alpha }_{\ddot{r}+1,\,{s}-1}\Bigr) \,\, \ddot{n}^{s} \psi_{k+1,\ddot{n}+\ddot{r}}(x)\, . 
\nonumber \end{eqnarray}

\section{Concrete procedures for the algorithm}
\label{sec:cp}
%\subsection{necassary procedures for the algorithm}

In this section, we explain the concrete procedures of our algorithm in detail, though its basic idea was sketched in the {\bf Algorithm} in Section 2 of~\cite{paper1}. As a first explanation, for simplicity, we will explain it for the cases with $\dim V\cap \ell ^2(\mathbb{Z}^+) =1$, though a similar method is possible even when $\dim V\cap \ell ^2(\mathbb{Z}^+) \ge 2$ as will be shown in Sections \ref{sec:multidim} and \ref{sec:conv_multidim}. 

%We begin with some definitions. 
In the following, let $\Pi_m$ $(m\in\mathbb{Z}^+ )$ be the projector defined in (\ref{eqn:1}) of Section \ref{sec:sv}. 
%With this projector and the space $V$ defined (\ref{eqn:def_sp_sol_infty}), define the spaces 
%$\Pi _m V:=\{ \Pi_m\vec{f}\, |\, \vec{f}\in V\} \,\,\,\,\, (m\in \mathbb{Z}^+ )$ 
%In particular, $\Pi _{p_0}$  with $p_0$ defined in (\ref{eqn:def_p0}) is important. 
% \begin{eqnarray}
% \tilde{V}:=V_{p_0} \,\,\,\,\,\, {\rm with } \,\,\,\,\,\, p_0:=j_0+\ell _0 -1 \,\, , 
% \label{eqn:def_tildeU}\end{eqnarray}
% where the integers $\ell _0$ and $j_0$ are the same as the integers in Conditions {\bf C2} and {\bf C5} in Section \ref{sec:sv}, respectively. 
Under the choice of function spaces and basis systems introduced in Section \ref{sec:fb}, from (\ref{eqn:lj}), 
\begin{eqnarray*}
p_0=\max(2M+k_0,\, 2M+k_0-k_0^\Diamond )\,\, .
\end{eqnarray*}
% In addition, define the spaces 
% \begin{eqnarray}
% V_m^{\,\ell ^2}:=\bigl\{\Pi _m\vec{f}\, \bigl|\, \vec{f} \in V\cap \ell ^2(\mathbb{Z}^+)\bigr\}\,\,\,\,\, (m\in \mathbb{Z}^+ ) \, . 
% \label{eqn:def_Ul2}\end{eqnarray}
Moreover, we introduce two inner products and their corresponding norms in $V$ and $\Pi _mV$ under the inequality (\ref{OC7}), 
\begin{eqnarray}\label{eqn:def_ip}
\langle \vec{f}, \vec{g}\rangle _{\ell ^2,K}&:=& \langle \Pi _K\vec{f}, \Pi _K\vec{g}\rangle _{\ell ^2}, 
\,\,\,\,\,\,\,\,\,\, 
\langle \vec{f}, \vec{g}\rangle _{Q,N}:=\Omega (\Pi _N\vec{f},\Pi _N\vec{g}), \\ 
\nonumber 
\|\vec{f}\|_{\ell ^2,K}^2 &:=& \langle \vec{f}, \vec{f}\rangle _{\ell ^2,K}, 
\,\,\,\,\,\,\,\,\,\, \,\,\,\,\,\,\,\,\,\, \,\,\,\,\,\,  
\|\vec{f}\|_{Q,N}^2 := \langle \vec{f}, \vec{f}\rangle _{Q,N}=\Omega (\Pi _N\vec{f}),
\end{eqnarray}
where the statement $\|\vec{f}\|_{\ell ^2,K}\!\!\! =0 \, \Longrightarrow \, \vec{f}=0$ is guaranteed, for $\vec{f}\in V$ or $\vec{f}\in \Pi _m V$,   because of {\bf C2} and (\ref{OC7}).

The algorithm consists mainly of three parts, where the first ({\bf Step 1} of {\bf Algorithm} of~\cite{paper1}) is the calculation of the basis vectors of $\Pi _{p_0} V$ and the second ({\bf Step 2}) is  the calculation of the basis vectors of $\Pi _N V$ and the third ({\bf Step 3}) is the removal of the components which do not belong to $\Pi _N (V\cap \ell ^2(\mathbb{Z}^+))$ from linear combinations of these basis vectors. 

For {\bf Step 1} and {\bf Step 2}, the matrix elements $b_m^n$ $(m,n\in \mathbb{Z}^+; \, |m-n|\le \ell _0,\, n\le N)$ are required. As has been explained in Section \ref{sec:me}, these elements can be calculated from just the coefficients $\ddot{\beta }_{\ddot{r},s}$  
$(\ddot{m},\ddot{n}\in \mathbb{Z} ; \, \ddot{n}-M\le \ddot{m}\le \ddot{n}+M+k_0-k_0^\Diamond )$, by a simple substitution of $m$ and $n$ into the power expansion  (\ref{eqn:expansion_matrix_element}). Hence the calculation of these coefficients suffices, and it should be carried out before {\bf Step 1}; we regard it as a preliminary step {\bf Step 0}.

The step {\bf Step 3} consists of the three sub-steps {\bf Step 3.1.a}-{\bf Step 3.1.c} given below. With these remarks, the concrete procedures of the algorithm can be stated, as follows: 

\vspace{3mm}
\noindent {\bf Algorithm}
\vspace{3mm}
\begin{description}
\item[Step 0]
Calculation of coefficients $\ddot{\beta }_{\ddot{r},s}$ 
\item[] 
\hspace{5mm}
Calculate $\ddot{\beta }_{\ddot{r},s}$ $(\ddot{m},\ddot{n}\in \mathbb{Z} ; \, \ddot{n}-M\le \ddot{m}\le \ddot{n}+M+k_0-k_0^\Diamond )$ by the procedure  given in Table \ref{tbl:1}, as is explained in  Section \ref{sec:me}.  
\item[Step 1]
Calculation of basis vectors of $\Pi _{p_0} V$:
\item[] 
\hspace{5mm} Find a basis system $\{F_n^{(1)}\}_{n=0}^{p_0} ,\ldots, \{F_n^{(D)}\}_{n=0}^{p_0}$
 of $\Pi _{p_0} V$  by Gaussian elimination, with the matrix elements $b_m^n$ calculated by (\ref{eqn:expansion_matrix_element}) and the result of {\bf Step 0}, where the dimension $D$ is determined automatically by Gaussian elimination. This is easy because $p_0$ is small.% and it can be done by integer-type programs because of the rationality of variables if the matrix elements $b_m^n$ with $0\le m\le j_0-1$ and $0\le n \le p_0$ are rational-(complex-)valued.
\item[Step 2] Recursive calculation of the basis vectors of $\Pi _nV$ $(p_0+1\le n \le N)$: 
\item[]
\hspace{5mm} Iterate the recursion (\ref{eqn:recursive_alg}) below for $n=p_0+1,p_0+2,\ldots , N$, with the matrix elements $b_m^n$ calculated by (\ref{eqn:expansion_matrix_element}) and the result of {\bf Step 0}, in order to obtain a basis system $\{F_n^{(1)}\}_{n=0}^{N} ,\ldots, \{F_n^{(D)}\}_{n=0}^{N}$ of $\Pi _N V$.
\item[Step 3] Removal of components from $\Pi _N V$ corresponding to non-$\ell ^2$-ones in $V$: 
\begin{description}
\item[Step 3.1.a] Integer-type quasi-orthogonalization of the basis system of $\Pi _N V$:
\item[]
\hspace{5mm} Find a system of linear combinations $\{E_n^{(1)}\}_{n=0}^{N} ,\ldots, \{E_n^{(D)}\}_{n=0}^{N}$ of \par\noindent $\{F_n^{(1)}\}_{n=0}^{N} ,\ldots, \{F_n^{(D)}\}_{n=0}^{N}$  which is sufficiently close to an orthogonal system with respect to the inner product $\langle \cdot ,\,\cdot \rangle _{Q,N}$, by the procedures %given in Table 5.1
explained below which is based on an intermediate idea between the Gram-Schmidt process and the Euclidean algorithm. 
\item[Step 3.1.b] Selection of minimum-ratio vector:
\item[]
\hspace{5mm} Find $\displaystyle \{G_n^{(1)}\}_{n=0}^{N}:=\{E_n^{(d_{\rm opt.})}\}_{n=0}^{N}$ with $d_{\rm opt.}:=\mathop{\rm argmin}_{1\le d\le D}\,\,\sigma_{K,N}^{(\Omega)}(\vec{E}^{(d)})$.
%with minimum ratio $\displaystyle\frac{\|\cdot\|_{Q,N}}{\|\cdot\|_{\ell ^2,K}}$ in the linear combinations obtained by {\bf Step 3.a}. 
\item[Step 3.1.c] Truncation (projection) by $\Pi _K$:
\item[]
\hspace{5mm}
Project the result of {\bf Step 3.1.b} to $\Pi _K V$. 
\end{description}
%\hspace{5mm} Find %the vector $\vec{f}_{K,N}^{(Q)}$ in (\ref{eqn:def_ratio2}), or find 
% a vector in the set $O_{K,N}^{(Q,c)}$ in {\bf C8} (which is identical to $O_{K,N}^{(Q,c)}$ in Condition {\bf C8}) from the basis vectors $\{F_n^{(1)}\}_{n=0}^{N} ,\ldots, \{F_n^{(D)}\}_{n=0}^{N}$ of $\Pi _N V$.% instead. As will be shown later, the latter alternative choice spares the quantity of calculation remarkably. 
\end{description}
\vspace{3mm}
From Conditions {\bf C2}, {\bf C5} and the definition (\ref{eqn:def_sp_sol_infty}), the calculations in {\bf Step 2} can be performed by the recursion 
\begin{eqnarray}
F_n^{(d)} = -\,\frac{1}{b_{n-\ell _0 }^n} \sum_{r=n-2\ell _0 }^{n-1} b_{n-\ell _0}^r F_r^{(d)} 
% F_n^{(d)} = \frac{1}{b_{n-\ell _0 }^n} \sum_{r=n-2\ell _0 }^{n-1} b_{n-\ell _0}^r F_r^{(d)} \,\,\,\,\,\,  (d=1,2,\ldots D)   
\label{eqn:recursive_alg}\end{eqnarray} 
with unchanged $F_s^{(d)}$ $(0\le s\le n-1)$, for $n=p_0+1,\, p_0+2,\, \ldots , \, N$.

The results of {\bf Step 0}-{\bf Step 2} belong to $\Pi _N V$. However, they contain components in $\Pi _N V\backslash (\Pi _N (V \cap \ell ^2(\mathbb{Z}^+)))$ $\bigl(=\{\Pi _N \vec{f} \, |\, \vec{f}\in V\backslash (V\cap \ell ^2(\mathbb{Z}^+ )\}\bigr)$ which have nothing to do with the  true solutions in $C^M(\mathbb{R} )\cap L_{(k)}^2(\mathbb{R} )$ of the differential equation. Hence, we should remove the components in $\Pi _N V\backslash (\Pi _N (V \cap \ell ^2(\mathbb{Z}^+)))$. {\bf Step 3} can almost remove them in the following sense, though we will prove it in detail later in this section. 
The orthogonalization with respect to  $\langle \cdot, \cdot\rangle _{Q,N}$ of the basis system of $\Pi _N V$ provides us with vectors sufficiently close to $\Pi _N (V \cap \ell ^2(\mathbb{Z}^+))$ with respect to  $\|\cdot\|_{\ell ^2,K}$ such that they belong to 
$((\sigma_{K,N}^{(\Omega )})^{-1}[0,c \underline{\sigma_{K,N}^{(\Omega )}} ])$ in Section 2.
%$O_{K,N}^{(Q,c)}$ 
%in {\bf C8}. %, as will be shown later.  
For this, a `quasi-orthogonalization' is sufficient, where the angles between any pair of vectors are sufficiently close to $\pi /2$, as will be shown later. Since exact orthogonalization (without round-off errors) by the Gram-Schmidt process requires many calculations for large $N$ and $D$, we will use the quasi-orthogonalization, without round-off errors but with fewer calculations. 

In the following, we explain in detail how the procedures in {\bf Step 3.1.a}-{\bf Step 3.1.c} can be realized by integer-type programs. 

Since the vector $\vec{F}^{(d)}$ is rational-(complex-)valued, there is an integer $C_d$ such that $\vec{F}_{\rm int.}^{(d)}:=C_d\vec{F}^{(d)}$ is (complex-)integer-valued. Then, {\bf Step 3.1.a} can be performed by the replacement procedures in Table \ref{tbl:2} for $\vec{v}_d^{\rm \, (initial)}=\vec{F}_{\rm int.}^{(d)}$ $(d=1,\ldots ,D)$, with  a  positive integer $h$. This is because of the inequality 
$\displaystyle\frac{|\langle \vec{v}_j, \vec{v}_\ell \rangle _{Q,N} |}{\|\vec{v}_j\|_{Q,N} \cdot\|\vec{v}_\ell \|_{Q,N}} \le \frac{1}{h}$ for $\, 1 \le j<\ell \le D$, which is guaranteed after the procedures in Table \ref{tbl:2}, whose 
%termination 
halting  
will be proved in Section \ref{sec:halt}. 

The iteration of {\bf Q1} in Table \ref{tbl:2} looks somewhat like the `lattice reduction problem'~\cite{lat1}~\cite{lat2} (known to be an  NP-hard problem), but this iteration  is an `imperfect' lattice reduction with few complex calculations which guarantees only the inequalities 
\begin{eqnarray*}
\Bigl|{\rm Re}\,  \Bigl(\vec{v}_j,\,  \vec{v}_\ell \Bigr)\Bigr|\le \frac{1}{2} \min (\|\vec{v}_j\|^2, \|\vec{v}_\ell \|^2) \,\,\,\,\,\, {\rm and} \,\,\,\,\,\, 
\Bigl|{\rm Im}\, \Bigl(\vec{v}_j,\,  \vec{v}_\ell \Bigr)\Bigr|\le \frac{1}{2} \min (\|\vec{v}_j\|^2, \|\vec{v}_\ell \|^2)
\end{eqnarray*}
for $1\le j<\ell \le D$, which are derived from the inequalities 
\begin{eqnarray*}
\Bigl|{\rm Re}\,  \Bigl(\vec{a},\,  \vec{b}-\Bigl[\frac{(\vec{b}, \vec{a})}{(\vec{a}, \vec{a})}\Bigr]_\mathbb{C} \vec{a}\Bigr)\Bigr|\le \frac{1}{2} \|\vec{a}\|^2 \,\,\,\,\,\, {\rm and} \,\,\,\,\,\, 
\Bigl|{\rm Im}\, \Bigl(\vec{a},\,  \vec{b}-\Bigl[\frac{(\vec{b}, \vec{a})}{(\vec{a}, \vec{a})}\Bigr]_\mathbb{C} \vec{a}\Bigr)\Bigr|\le \frac{1}{2} \|\vec{a}\|^2 \,\, .
\end{eqnarray*}
This iteration is used only as a preliminary procedure for the iteration of {\bf Q2}, where we are not aiming for the exactly minimal basis system of the lattice but instead, good enough orthogonality. Moreover, the iteration of {\bf Q1} can be regarded as a combination of the Gram-Schmidt process and a multidimensional complex version of the Euclidean algorithm. 
The final results of these procedures give the vectors $\vec{E}^{(d)}:=\vec{v}_d^{\rm \, (final)}$ $(d=1,2,\ldots ,D)$. 

\begin{table}[thb]
\caption{Basic idea of integer-valued quasi-orthogonalization}
\begin{center}
\begin{tabular}{|l@{}|}
\hline
$\left(\begin{array}{@{\,}ll}
\mbox{In the following, } \\ 
\mbox{ $ [z]_\mathbb{C} := -  {\rm sgn} ({\rm Re} \, z) \lfloor\frac{1}{2} - |{\rm Re} \, z|\rfloor  - i\,\, {\rm sgn} ({\rm Im} \, z)  \lfloor\frac{1}{2} - |{\rm Im} \, z| \rfloor $ }
\end{array}\right)$ \\ 
\hline
$\begin{array}{@{\,}ll}    
\mbox{{\bf Iterate Q1} below until nothing is changed.} \\ 
\,\,\, \left\lfloor\begin{array}{@{\,}ll} 
\mbox{{\bf Q1}} \left\{\begin{array}{@{\,}ll} 
\bullet \mbox{\, sorting and renumbering  of $\vec{v}_{1}, \, \vec{v}_{2}, \, \ldots \, , \, \vec{v}_D$ in order that } \\ \mbox{$\,\,\,\,\,\,\,\,\, \| \vec{v}_{1}\|_{Q,N} \le \| \vec{v}_{2}\|_{Q,N} \le \ldots \le \| \vec{v}_D\|_{Q,N} $} \\  
\bullet \mbox{\, {\bf for} $j=2$ to $D$} \\ 
\,\,\,\,\, 
\,\,\, \left\lfloor\begin{array}{@{\,}ll}
\mbox{{\bf for} $\ell =1$ to $j-1$} \\  
\,\,\, \left\lfloor\begin{array}{@{\,}ll}
\mbox{$\displaystyle \vec{v}_j\leftarrow\vec{v}_j-\Bigl[\frac{\langle\vec{v}_j, \vec{v}_\ell \rangle_{Q,N}}{\langle\vec{v}_\ell , \vec{v}_\ell \rangle_{Q,N}} \Bigr]_\mathbb{C} \vec{v}_\ell $} \\ 
\end{array}\right. \\ 
\end{array}\right. \\ 
\end{array}\right. \\  
\end{array}\right. \\  
\left(\begin{array}{@{\,}ll}
\mbox{This is a preparatory `imperfect'  lattice  } \\  \mbox{ reduction for the iteration of {\bf Q2} below.}
\end{array}\right)
\end{array}$ \\  
\hline 
$\begin{array}{@{\,}ll} \
\mbox{{\bf Iterate Q2} below until nothing is changed.} \\ 
\,\,\, \left\lfloor\begin{array}{@{\,}ll} 
\mbox{{\bf Q2}} 
\left\{\begin{array}{@{\,}ll} \mbox{\, {\bf for} $j=2$ to $D$} \\ 
\,\,\, \left\lfloor\begin{array}{@{\,}ll} 
\mbox{{\bf for} $\ell =1$ to $j-1$} \\  
\,\,\, \left\lfloor\begin{array}{@{\,}ll}  
\mbox{{\bf if} $\, h^2\cdot |\langle \vec{v}_j, \vec{v}_\ell \rangle_{Q,N}|^2 \ge  \|\vec{v}_j\|_{Q,N}^2 \cdot \|\vec{v}_\ell \|_{Q,N}^2 $,} \\  
\, \left\lfloor\begin{array}{@{\,}ll} 
\,\,\,\mbox{$\vec{v}_j \leftarrow 2 \vec{v}_j$} \\
\end{array}\right. \\ 
\mbox{$\displaystyle \vec{v}_j\leftarrow\vec{v}_j-\Bigl[\frac{\langle\vec{v}_j, \vec{v}_\ell \rangle_{Q,N}}{\langle\vec{v}_\ell , \vec{v}_\ell \rangle_{Q,N}} \Bigr]_\mathbb{C} \vec{v}_\ell $} \\ 
\end{array}\right. \\ 
\end{array}\right. \\ 
\end{array}\right. \\ 
\end{array}\right. \\  
\left(\begin{array}{@{\,}ll}
\mbox{This process provides us with the vectors } \\  \mbox{$\vec{v}_j$ $(1\le j \le D)$ \, s.t. $\displaystyle\frac{|\langle \vec{v}_j, \vec{v}_\ell \rangle _{Q,N} |}{\|\vec{v}_j\|_{Q,N} \cdot\|\vec{v}_\ell \|_{Q,N}} \le \frac{1}{h}$}\\ 
 \,\,\,\,\, \mbox{for $\, 1 \le j<\ell \le D$\,. }
 \\ 
\end{array}\right)   
\end{array}$ \\ 
\hline
\end{tabular}
\end{center}
\label{tbl:2}
\end{table}

For {\bf Step 3.1.b}, with $\displaystyle d_{\rm opt.}:=\mathop{\rm argmin}_{d\in\{1,\ldots ,D\}} \displaystyle \frac{\|\vec{v}_d^{\rm \, (final)}\|_{{Q,N}}}{\|\vec{v}_d^{\rm \, (final)}\|_{\ell ^2,K}}$ for the vectors $\vec{v}_d^{\rm \, (final)}$ $(d=1,\ldots D)$ after the replacement procedures  in Table \ref{tbl:2},  we then define $\vec{G}^{(1)}:=\vec{v}_{d_{\rm opt.}}^{\rm \, (final)}=\vec{E}^{d_{\rm opt.}}$.
%, and then we obtain $\Pi _K \vec{G}^{(1)}$ 
%define  $\Pi _K \vec{G}:=\Pi_K\vec{v}_{d_{\rm opt.}}=\Pi_K\vec{G}$ 
for {\bf Step 3.1.c}. 
% \begin{eqnarray} 
% \begin{eqnarray} 
% A_{K,N}^{(Q,c)}&:=&\{ a\vec{f}\, |\, \vec{f}\in O_{K,N}^{(Q,c)},\, a>0  \} \\ &=& 
% \biggl\{ \vec{f}\in \Pi _N V\backslash \{ 0\} \,\, \biggl|\, \,  \frac{\|\vec{f}\|_{Q,N}}{\|\vec{f}\|_{\ell ^2,K}}\, \le  \sqrt{c}  \inf_{ \vec{g} \in V\backslash\{ 0\}} \frac{\|\vec{g}\|_{Q,N}}{\|\vec{g}\|_{\ell ^2,K}}
% \biggr\}  \nonumber 
% \label{eqn:def_AKNQc}\end{eqnarray}
% with $O_{K,N}^{(Q,c)}$ 
% as 
% defined  
% in {\bf C8}. 
Then, we have the following theorem: 
\begin{theorem}
\label{thm:subopt}
When $h\ge D$, the vector $\vec{G}^{(1)}$
%$\vec{G}$  
belongs to $((\sigma_{K,N}^{(\Omega )})^{-1}[0,c \underline{\sigma_{K,N}^{(\Omega )}} ])$
%$A_{K,N}^{(Q,c)}$ with \par\noindent  $\displaystyle c=\frac{D}{1- \frac{D-1}{h}} \,  $,  
\end{theorem}
whose proof  will be given in Section \ref{sec:suboptimality}.

% Since there exists a nonzero positive number $a$ such that $\displaystyle \frac{1}{a}\Pi _K \vec{G}\in O_{K,N}^{(Q,c)}$, the 
Hence, the following  Theorem \ref{thm:opt_ratio_gen} shows the approach  of $\Pi_K\vec{G}^{(1)}$%\Pi _K \vec{G}$ 
(obtained in {\bf Step 3.1.c}) to the  vector space corresponding to the function space of true solutions in $C^M(\mathbb{R} \setminus S)\cap L_{(k)}^2(\mathbb{R} )$ (of the differential equation) projected to  the $K$-dimensional subspace  $<e_0,e_1,\ldots ,e_K>$, in the sense that 
% $\displaystyle \sup_{
% \vec{x} \in 
% (\sigma_{K,N}^{(\Omega )})^{-1} [0,c \underline{\sigma_{K,N}^{(\Omega )}} ] }
% \frac{ \|P_{V,K} \vec{x} -\vec{x}\|_{\ell ^2,K} }{\|\vec{x}\|_{\ell ^2,K} }$ 
$\displaystyle 
\frac{\|
P_{V,K}\Pi_K\vec{G}^{(1)}-\Pi _K\vec{G}^{(1)}
\|_{\ell ^2}}{\|\Pi_K\vec{G}^{(1)}\|_{\ell ^2}}$ 
converges to $0$ as $N$ tends to infinity: 

\begin{theorem}
\label{thm:opt_ratio_gen}
For fixed $K$, when $N$ goes to infinity,
the convergence
\begin{eqnarray}
\sup_{
\vec{x} \in 
(\sigma_{K,N}^{(\Omega )})^{-1} [0,c \underline{\sigma_{K,N}^{(\Omega )}} ] }
\frac{ \|P_{V,K} \vec{x} -\vec{x}\|_{\ell ^2,K} }{\|\vec{x}\|_{\ell ^2,K} }
\to 0
% With a constant parameter $c>1$, under {\bf C6} and with $O_{K,N}^{(Q,c)}$ as in ${\bf C8}$, 
% \begin{eqnarray*}
% \lim_{N\to\infty}\,\, \max_{\vec{g}\in O_{K,N}^{(Q,c)}} \,\, \inf_{\vec{f}\in V\cap \ell ^2(\mathbb{Z}^+ )}
% \|
% \Pi_K \vec{g}
% -\Pi_K \vec{f}
% \|_{\ell ^2} = 0.
\end{eqnarray}
\end{theorem}  \par\noindent 
Theorem 2.5 of~\cite{paper1} is equivalent to this, and its proof is given in Appendix of D of~\cite{paper1}.

For the practical algorithm, 
the procedures {\bf Q1} and {\bf Q2} can be executed with fewer calculations by replacements of the coefficients, instead of vectors, in the expansion of the vectors in terms of the initial basis vectors, as shown in Table \ref{tbl:3}. 
Moreover, we can propose a recursive 
%innovation 
renewal  
over $n$ for the product $\Phi_{n,r}^{(d)}:=F_n^{(d)} \overline{F_{n-r}^{(d)} }$ $(r=0,1,\ldots , 2\ell _0+1; \, d=1,2,\ldots ,D)$ instead of direct calculation of the inner products, which reduces the order of the required 
%number 
amount  
of calculations of inner products between the vectors. 

The order the number of calculations required for obtaining the inner products, which is empirically the narrowest bottle neck of our method, is $O\bigl(N^3(\log N)^2\bigr)$, because the number of digts of numerators and denominators of $F_n^{(d)}$ $(n\in \mathbb{Z}^+)$ can be bounded by $O(n \log n)$ (and the denominators of larger $n$ are always multiples of the denominators of smaller $n$). However, by these modifications, the order of the number of calculations required for obtaining the inner products, which is empirically the narrowest bottle neck of our method, is $O\bigl(N^2(\log N)^2\bigr)$. 
This is not too large because the number of significant digits of the solution obtained 
by our method 
is not fixed, but is an  increasing function of $N$ 
there 
%in our method 
because of the  accuracy increasing as $N$ increases, as is expected from the discussion before Theorem \ref{thm:opt_ratio_gen}. 
In practical use, we augment the dimension $K$ as well as $N$ so that they may be proportional to each other, for high accuracy. Then, the number of significant digits of the solution obtained increases without limit, which is observed empirically in numerical experiments.  
(Empirically, in most cases, the number of significant digits increases almost in proportion to a power of $N$.)  
A similar fact  can be shown mathematically with the order of limits `$\displaystyle\lim_{K\to\infty } \lim_{N\to\infty }$' because $\displaystyle\frac{\|(1-\Pi _K)\vec{f}\|_{\ell ^2}}{\|\vec{f}\|_{\ell ^2}}$ converges to $0$ for $\vec{f}\in V\cap \ell ^2(\mathbb{Z}^+ )$ as $K$ tends to infinity. Some of these numerical  results will be presented in Section \ref{sec:nm}. 

\begin{table}[p]
\caption{Practical operations for integer-type quasi-orthogonalization }
\begin{center}
\scalebox{.86}{
\begin{tabular}{|l@{}|}
\hline 
$\begin{array}{@{\,}ll}  
\mbox{{\bf for} $j=1$ to $D$} \\  
\,\,\, \left\lfloor\begin{array}{@{\,}ll} 
\mbox{{\bf for} $m=1$ to $D$} \\  
\,\,\, \left\lfloor\begin{array}{@{\,}ll} 
\mbox{$c_{jm} \leftarrow \delta _{jm}$} \\ 
\mbox{$p_{jm} \leftarrow \langle \vec{v}_j^{\, ({\rm initial})} , \vec{v}_{m}^{\, ({\rm initial})} \rangle _{Q,N}$} \\
\mbox{$q_{jm} \leftarrow \langle \vec{v}_j^{\, ({\rm initial})} , \vec{v}_{m}^{\, ({\rm initial})} \rangle _{\ell^2,K}$} \\
\end{array}\right. \\ 
\end{array}\right. \\ 
\end{array} $ \\ 
\,\,\,\,\,\,\,\,\, \mbox{$\left( \, \begin{array}{@{\,}ll} c_{jm}: \mbox{complex-integer-velued coefficients} \\ \displaystyle 
\,\,\,\,\,\,\,\,\, \,\,\,\,\,\, \mbox{ for } \,\,\,\,\, \vec{v}_j^{\, ({\rm new})}=\sum_{m=1}^D c_{jm} \vec{v}_m^{\, ({\rm initial})} \\ 
\end{array} \, \right) $} \\ 
\hline   
$\begin{array}{@{\,}ll} 
\mbox{{\bf Iterate P1} below until nothing is changed.} \\  
\,\,\, \left\lfloor\begin{array}{@{\,}ll}
\mbox{{\bf P1}} \left\{\begin{array}{@{\,}ll} \bullet \mbox{\, with a permutation $\,  n_{1}, n_{2},\ldots,n_D$\,  of } \\ 
\,\,\, \mbox{$\, 1,2,\ldots,D $\, s.t. $\displaystyle \, p_{1\, 1} \le p_{2\, 2} \le \ldots \le p_{D D}$}, \,\,\,\,\,\,\,\,\,\,\,\,\,\,\,\, \\ 
\,\,\,\,\,\,\,\,
\mbox{{\bf for} $j=1$ to $D$} \\ 
\,\,\,\,\,\,\,\,
\,\,\left\lfloor\begin{array}{@{\,}ll} 
\mbox{$c_{jm} \leftarrow c_{n_j m} ,\,\, p_{jm} \leftarrow p_{n_j n_m} ,\,\,\, q_{jm} \leftarrow p_{n_j n_m}$} \\ 
\end{array}\right.  \\ 
\bullet \mbox{\, {\bf for} $j=2$ to $D$} \\ 
\,\,\,\,\, 
\,\,\, \left\lfloor\begin{array}{@{\,}ll}
\mbox{{\bf for} $\ell =1$ to $j-1$} \\  
\,\,\, \left\lfloor\begin{array}{@{\,}ll} 
\mbox{$\displaystyle r\leftarrow \Bigl[ \frac{p_{j\ell }}{p_{\ell \ell}} \Bigr]_\mathbb{C} $} \\ 
\mbox{{\bf for} $m=1$ to $D$} \\  
\,\,\, \left\lfloor\begin{array}{@{\,}ll}
\mbox{$ c_{jm}\leftarrow c_{jm}-r\, c_{\ell m} $,\,\, $\displaystyle p_{jm}\leftarrow p_{jm}-r\, p_{\ell m}$,\,\, $q_{jm}\leftarrow q_{jm}-r\,  q_{\ell m}\!\!\!\!\!\!\!\!\!\!\!\!\!\!\!\!\!\!\!\! $}\\ 
\end{array}\right.\\ 
\mbox{{\bf for} $m=1$ to $D$} \\  
\,\,\,\left\lfloor\begin{array}{@{\,}ll}
\mbox{$\displaystyle p_{mj}\leftarrow p_{mj}-\overline{r}\, p_{m\ell } ,\,\,\, q_{mj}\leftarrow q_{mj}-\overline{r}\, q_{m\ell }$}\\ 
\end{array}\right. \\
\end{array}\right. \\ 
\end{array}\right. \\ 
\end{array}\right. \\  
\end{array}\right. \\  
\end{array}$
\\ 
\hline 
$\begin{array}{@{\,}ll} 
\mbox{{\bf Iterate P2} below until nothing is changed.} \\  
\,\,\, \left\lfloor\begin{array}{@{\,}ll}  
\mbox{{\bf P2}}  \left\{ \begin{array}{@{\,}ll} 
\mbox{\, {\bf for} $j=2$ to $D$} \\ 
%\,\,\, \left\lfloor\begin{array}{@{\,}ll}
\mbox{{\bf for} $\ell =1$ to $j-1$} \\  
\,\,\, \left\lfloor\begin{array}{@{\,}ll} 
\mbox{{\bf if} $\displaystyle\,\,\, h^2 \,\left[\frac{p_{j\ell }}{\widehat{N}^K}\right]_\mathbb{C}  \left[\frac{\overline{p_{j\ell }}}{\widehat{N}^K}\right]_\mathbb{C} \ge \left(\Bigl\lfloor\frac{p_{jj}}{\widehat{N}^K}\Bigr\rfloor +1\right)\left(\Bigl\lfloor\frac{p_{\ell\ell }}{\widehat{N}^K}\Bigr\rfloor +1\right)$,} \\ 
\, \left\lfloor\begin{array}{@{\,}ll} 
\mbox{{\bf for} $m=1$ to $D$} \\ 
\,\,\, \left\lfloor\begin{array}{@{\,}ll}
\mbox{$p_{mj} \leftarrow 2 p_{mj}$,\,\, $q_{mj} \leftarrow 2 q_{mj}$} \\
\mbox{$p_{jm} \leftarrow 2 p_{jm}$,\,\, $q_{jm} \leftarrow 2 q_{jm}$} \\
\end{array}\right. \\ 
\end{array}\right. \\ 
\mbox{$\displaystyle r\leftarrow \Bigl[ \frac{p_{j\ell }}{p_{\ell \ell}} \Bigr]_\mathbb{C} $} \\ 
\mbox{{\bf for} $m=1$ to $D$} \\  
\,\,\, \left\lfloor\begin{array}{@{\,}ll}
\mbox{$\displaystyle c_{jm}\leftarrow c_{jm}-r\, c_{\ell m}$, \,\, $\displaystyle p_{jm}\leftarrow p_{jm}-r\, p_{\ell m} ,\,\,\, q_{jm}\leftarrow q_{jm}-r\, q_{\ell m} $} \\ 
\end{array}\right.\\ 
\mbox{{\bf for} $m=1$ to $D$} \\  
\,\,\, \left\lfloor\begin{array}{@{\,}ll}
\mbox{$\displaystyle p_{mj}\leftarrow p_{mj}-\overline{r}\, p_{m\ell } ,\,\,\,  q_{mj}\leftarrow q_{mj}-\overline{r}\, q_{m\ell } $} \\ 
\end{array}\right. \\
\end{array}\right. \\ 
\end{array}\right. \\ 
\end{array}\right. \\
\end{array}$ \\ 
\hline 
$\begin{array}{@{\,}ll}  
\mbox{{\bf for} $j=1$ to $D$} \\ 
\,\,\, \left\lfloor\begin{array}{@{\,}ll} 
\mbox{$\displaystyle \vec{v}_j^{\rm \, (final)} \leftarrow \sum_{m=1}^D c_{jm} \vec{v}_m^{\rm \, (initial)}  $} \\ 
\mbox{{\bf for} $\ell =1$ to $D$} \\  
\,\,\, \left\lfloor\begin{array}{@{\,}ll}
\mbox{{\bf if} $j\ge 1$ or $m\ge 1$, } \\ 
\, \left\lfloor\begin{array}{@{\,}ll}
\mbox{$\displaystyle \langle \vec{v}_j^{\rm \, (final)} , \vec{v}_\ell^{\rm \, (final)} \rangle _{Q,N} \leftarrow \sum_{m=1}^D \sum_{n=1}^D c_{jm} \overline{c_{\ell n}} \, p_{mn}$} \\ 
\mbox{$\displaystyle \langle \vec{v}_j^{\rm \, (final)} , \vec{v}_\ell^{\rm \, (final)} \rangle _{\ell^2,K} \leftarrow \sum_{m=1}^D \sum_{n=1}^D c_{jm} \overline{c_{\ell n}} \, q_{mn}$}  \\
\end{array}\right. \\
\end{array}\right. \\
\end{array}\right. \\ 
\end{array} $ \\  
\hline
\end{tabular}
}
\end{center}
\label{tbl:3}
\end{table}

Huge integers can be treated by integer arrays for the base-$10^9$ positional notation $\displaystyle \sum_{\ell } c_\ell (10^{9})^\ell $ with integers $c_\ell $ which are smaller than $10^9$. We constructed practical program modules for the four arithmetical operations on these integer array expressions of huge integers.

Even when $D_{\ell ^2}:= \dim V\cap \ell ^2(\mathbb{Z}^+)\ge 2$, some modifications  enable us to obtain a quasi-orthogonal vector system   $\vec{G}^{(1)}, \vec{G}^{(2)}, \ldots ,  \vec{G}^{(D_{\ell ^2})}$ in $((\sigma_{K,N}^{(\Omega )})^{-1}[0,c \underline{\sigma_{K,N}^{(\Omega )}} ])$\, with respect to the inner product   $\langle \cdot , \cdot \rangle_{\ell ^2,K}$.\ and their projections $\Pi _K \vec{G}^{(1)}, \Pi _K \vec{G}^{(2)}, \ldots ,  \Pi _K \vec{G}^{(D_{\ell ^2})}$ to $\Pi _K V$. The details will be given in Section \ref{sec:multidim}.

\section{Suboptimality of the vectors obtained by {\bf Step 3.1.a} and {\bf Step 3.1.b}}
\label{sec:suboptimality}
In this section, we will prove Theorem \ref{thm:subopt}, which guarantees that the vectors $\vec{G}^{(1)}$ obtained by {\bf Step 3.1.a} and {\bf Step3.1.b} belong to $((\sigma_{K,N}^{(\Omega )})^{-1}[0,c \underline{\sigma_{K,N}^{(\Omega )}} ])$ with finite fixed $c$. For this, we begin with %two lemmata.
a lemma.  
%One of these lemmata will be essential later also for the upper bound on errors.

\begin{lemma}
\label{lemma:max_ratio_pl}
Let $U$ be a finite dimensional space and $(\cdot , \, \cdot)_\Lambda $ and $(\cdot , \, \cdot)_\Xi $ be inner products there. If nonzero vectors $\vec{f}_1,\, \vec{f}_2,\,  \ldots \vec{f}_{n}$\, $(n\le \dim U)$ in $U$ satisfy $\displaystyle\frac{|(\vec{f}_m,\, \vec{f}_\ell )_\Xi |}{\|\vec{f}_m\|_\Xi \cdot \|\vec{f}_m\|_\Xi } \le \zeta $\,   $(1\le m < \ell \le n)$ with fixed $\zeta$ such that $\displaystyle 0<\zeta<\frac{1}{n-1}$, then the inequality $\displaystyle \frac{\|\vec{f}\|_\Lambda}{\|\vec{f}\|_\Xi}\le \sqrt{\frac{\displaystyle\sum_{m=1}^{n}\frac{\|\vec{f}_m\|_\Lambda ^2}{\|\vec{f}_m\|_\Xi ^2}}{\,  1-(n-1)\zeta \, }}$ holds for any nonzero linear combination $\vec{f}$ of   $\vec{f}_1,\, \vec{f}_2,\, \ldots \vec{f}_{n}$.      
\end{lemma}
\par\noindent\noindent{\em Proof of Lemma }\ref{lemma:max_ratio_pl}: \quad
\rm 
For $\displaystyle\vec{f}:=\sum_{m=1}^{n}a_m\vec{f}_m$, define $\displaystyle A:=\sum_{m=1}^n|a_m|^2\|\vec{f}_m\|_\Xi ^2$. With the definitions $\displaystyle \vec{p}_m:=\frac{\vec{f}_m}{\|\vec{f}_m\|_\Xi }$,\, $\displaystyle C_m:=\frac{\|\vec{f}_m\|_\Lambda }{\|\vec{f}_m\|_\Xi}$ and $b_m:=a_m\, \|\vec{f}_m\|_\Xi $,\,\, we have \par\noindent  $\displaystyle\frac{|(\vec{p}_m,\, \vec{p}_\ell )_\Xi |}{\|\vec{p}_m\|_\Xi \cdot \|\vec{p}_m\|_\Xi } \le \zeta $,\, $\|\vec{p}_m\|_\Lambda = C_m $, \,  $\displaystyle\vec{f}=\sum_{m=1}^{n}b_m\vec{p}_m$ and $\displaystyle A=\sum_{m=1}^n|b_m|^2$. Then, 
\begin{eqnarray*}
\|\vec{f}\|_\Lambda ^2 = \bigl\|\sum_{m=1}^{n}b_m\vec{p}_m\bigr\|_\Lambda ^2 \le \Bigl( \sum_{m=1}^{n}\|b_m\vec{p}_m\|_\Lambda \Bigr)^2 
\le \Bigl( \sum_{m=1}^{n} |b_m|^2\Bigr) \Bigl( \sum_{m=1}^{n} \|\vec{p}_m\|_\Lambda ^2 \Bigr)
= A\sum_{m=1}^{n} C_m^2 .
\end{eqnarray*}
On the other hand, 
\begin{eqnarray*}
\|\vec{f}\|_\Xi ^2= \bigl\|\sum_{m=1}^{n}b_m\vec{p}_m\bigr\|_\Xi ^2 = \sum_{m=1}^{n}\sum_{\ell =1}^{n}b_m\overline{b}_\ell (\vec{p}_m , \vec{p}_{\ell } )_\Xi \ge  -\bigl( \sum_{m=1}^{n} |b_m| \bigr)^2 \zeta + \sum_{m=1}^n |b_m|^2 (1+\zeta) 
\\  \ge -\bigl(\sum_{m=1}^{n} |b_m|^2 \bigr) \, \bigl(\sum_{m=1}^{n} 1 \bigr) \,\zeta + A(1+\zeta) 
= -A\, n\,\zeta + A(1+\zeta) = A \bigl(1 - (n-1) \, \zeta \bigr)  . 
\end{eqnarray*}
Hence $\displaystyle \frac{\|\vec{f}\|_\Lambda ^2}{\|\vec{f}\|_\Xi ^2}\le \frac{\displaystyle A\sum_{m=1}^{n}C_m^2}{\,\,\,\,\, A\bigl(1-(n-1)\zeta \bigr)\,\,\,\,\, }=\frac{\displaystyle\sum_{m=1}^{n}\frac{\|\vec{f}_m\|_\Lambda ^2}{\|\vec{f}_m\|_\Xi ^2}}{\,\,\,\,\, 1-(n-1)\zeta \,\,\,\,\, }$.
\hfill\endproof 
\vspace{5mm}
\par\noindent This lemma enables 
us to prove the following theorem:
\begin{theorem}
\label{thm:min_ratio_n}
Let $U$ be a finite dimensional space and $(\cdot , \, \cdot)_\Lambda $ and $(\cdot , \, \cdot)_\Xi $ be inner products there. If nonzero vectors $\vec{f}_1,\, \vec{f}_2,\,  \ldots \vec{f}_{n}$\, $(n\le \dim U)$ in $U$ satisfy $\displaystyle\frac{|(\vec{f}_m,\, \vec{f}_\ell )_\Xi |}{\|\vec{f}_m\|_\Xi \cdot \|\vec{f}_m\|_\Xi } \le \zeta $\,\,   $(1\le m < \ell \le n)$ with  fixed $\zeta$ such that $\displaystyle 0<\zeta<\frac{1}{n-1}$,  then 
\begin{eqnarray*}
\min_j \frac{\|\vec{v}_j\|_\Xi }{\|\vec{v}_j\|_\Lambda } 
\,\, \le \,\, %\sqrt{\frac{1-(n-1)\zeta {n}}}
\sqrt{\frac{ n }{1-(n-1)\zeta }}
\, \, \inf_{\vec{u} \in U \backslash \{ 0 \}} \,\frac{\|\vec{u}\|_\Xi }{\|\vec{u}\|_\Lambda } \, . 
\end{eqnarray*}
\end{theorem}
\par\noindent\noindent{\em Proof of Theorem }\ref{thm:min_ratio_n}: \quad
\rm From %Lemmata \ref{lemma:max_ratio_pl} and \ref{lemma:min_ratio} 
Lemma \ref{lemma:max_ratio_pl} 
\begin{eqnarray*}
 \sup_{\vec{u} \in U \backslash \{ 0 \}} \,\frac{\|\vec{u}\|_\Lambda }{\|\vec{u}\|_\Xi }\le \sqrt{\frac{\displaystyle\sum_{m=1}^{n}\frac{\|\vec{f}_m\|_\Lambda ^2}{\|\vec{f}_m\|_\Xi ^2}}{\,  1-(n-1)\zeta \, }} \le \sqrt{\frac{\displaystyle n\, \max_j \frac{\|\vec{f}_j\|_\Lambda ^2}{\|\vec{f}_j\|_\Xi ^2}}{\,  1-(n-1)\zeta \, }} \le 
\sqrt{\frac{ n }{\,  1-(n-1)\zeta \, }}
\, \max_j \frac{\|\vec{f}_j\|_\Lambda }{\|\vec{f}_j\|_\Xi } .
\end{eqnarray*} 
This implies that 
\begin{eqnarray*}
 \inf_{\vec{u} \in U \backslash \{ 0 \}} \,\frac{\|\vec{u}\|_\Xi }{\|\vec{u}\|_\Lambda }\ge  \sqrt{\frac{\,  1-(n-1)\zeta \, }{ n }}\, \min_j \frac{\|\vec{f}_j\|_\Xi }{\|\vec{f}_j\|_\Lambda } \, .
\end{eqnarray*}  
\hfill\endproof
\par \par\noindent 

By means of Theorem  \ref{thm:min_ratio_n}, we can now prove Theorem  \ref{thm:subopt}. 
\par\noindent\noindent{\em Proof of Theorem }\ref{thm:subopt}: \quad
After the process {\bf Step 3.1.a}, the vectors $\vec{v}_1^{\rm \, (final)}, \ldots \vec{v}_D^{\rm \, (final)}$ satisfy the conditions in Theorem  \ref{thm:min_ratio_n}, with $\|\cdot\|_\Xi = \|\cdot\|_{Q,N} $ and $\displaystyle \zeta =\frac{1}{h}$. Hence, with $U=\Pi _N V$,  $\|\cdot\|_\Lambda = \|\cdot\|_{\ell ^2,K}$, $\|\cdot\|_\Lambda = \|\cdot\|_{Q,N}$, $n=D$ and $\displaystyle\zeta = \frac{1}{h}$, when $h\ge D$
\begin{eqnarray*}
 \sup_{\vec{g} \in \Pi _N V \backslash \{ 0 \}} \,\frac{\|\vec{g}\|_{\ell ^2,K} }{\|\vec{g}\|_{Q,N} }\le  \sqrt{\frac{ D }{\,  1-\frac{D-1}{h} \, }}\, \max_{d=1, \ldots ,D} \frac{\|\vec{v}_d^{\rm \, (final)}\|_{\ell ^2,K}}{\|\vec{v}_d^{\rm \, (final)}\|_{Q,N} }.
\end{eqnarray*} 
From the definition of $\vec{G}^{(1)}$ (made just before the theorem), this implies that 
% \begin{eqnarray*}
% \min_{d=1, \ldots ,D} \frac{\|\vec{v}_d\|_{Q,N} }{\|\vec{v}_d\|_{\ell ^2,K} }\le  \sqrt{\frac{ D }{\,  1-\frac{D-1}{h} \, }}\inf_{\vec{g} \in U \backslash \{ 0 \}} \,\frac{\|\vec{g}\|_{Q,N} }{\|\vec{g}\|_{\ell ^2,K} } \, . 
% \end{eqnarray*} 
\begin{eqnarray*}
\sqrt{\sigma_{K,N}^{(\Omega )}(\vec{G}^{(1)})}=
\frac{\|\vec{G}^{(1)}\|_{Q,N} }{\|\vec{G}^{(1)}\|_{\ell ^2,K}}\le 
\min_{d=1, \ldots ,D} \frac{\|\vec{v}_d^{\rm \, (final)}\|_{Q,N} }{\|\vec{v}_d^{\rm \, (final)}\|_{\ell ^2,K} }\le  \sqrt{\frac{ D }{\,  1-\frac{D-1}{h} \, }} 
\inf_{\vec{g} \in \Pi _N V \backslash \{ 0 \}} \,\frac{\|\vec{g}\|_{Q,N} }{\|\vec{g}\|_{\ell ^2,K} } \, . 
\end{eqnarray*} 
This inequality  
shows that the vector $\vec{G}^{(1)}$ belongs to $((\sigma_{K,N}^{(\Omega )})^{-1}[0,c \underline{\sigma_{K,N}^{(\Omega )}} ])$  with \par\noindent $\displaystyle c={\frac{D}{1- \frac{D-1}{h}}} \,  $, from the definition of $\vec{G}$ in the procedures for {\bf Step 3.2}.
\hfill\endproof\par 
Thus we have proved the validity of our method, by means of Theorem  \ref{thm:subopt} and Theorem  \ref{thm:opt_ratio_gen}.
\begin{rem}\begin{indention}{.8cm}\rm
For a practical algorithm, we propose an improvement where the iteration of {\bf P2} in Table \ref{tbl:3} is made also in the halfway steps in the recursion {\bf Step 2}. This improvement can reduce the size of the integers, and does not influence the structure of the algorithm because the iteration of {\bf P2} results only in a change of basis system.
\end{indention}\end{rem}

\section{Extension to the cases with $\dim V\cap \ell ^2(\mathbb{Z}^+ )\ge 2$}
\label{sec:multidim}
In this section, we explain how to extend the proposed method to the cases with $\dim V\cap \ell ^2(\mathbb{Z}^+ )\ge 2$, i.e. how to obtain approximately a quasi-orthogonal basis system of the subspace $V\cap \ell ^2(\mathbb{Z}^+ )$ with respect to $\langle \cdot , \cdot \rangle _{\ell ^2,K}$. In {\bf Algorithm} of our preceding paper~\cite{paper1}, for simplicity, we explained  {\bf Step 3.$n$} $(n\ge 2)$ using exact orthogonalization with respect to $\langle \cdot , \cdot \rangle _{\ell ^2,K}$. However, the exact orthogonalization requires a large amount of calculations and hence it is not practical. Moreover, the proposed method does not necessarily require the exact orthogonalization but quasi-orthogonalization is sufficient as is proved later. Therefore, in this paper, for {\bf Step 3.$n$}, we propose a practical procedures based on an quasi-orthogonalization which requires a relatively small amount of calculations.  This is a slight difference between the algorithms presented in~\cite{paper1} and this paper. For this difference, here we use step numbers {\bf Step 3$^\prime$.A} and {\bf Step 3$^\prime$.A.$\ldots$} 

When $\dim V\cap \ell ^2(\mathbb{Z}^+ )= 2$, the extension is possible, based on the idea of the quasi-minimization of the ratio $\displaystyle \frac{\|\vec{f}\|_{Q,N}}{\|\vec{f}\|_{\ell ^2,K}}$ in the subspace almost orthogonal to $\vec{G}^{\, (1)}$. Let $\vec{G}^{\, (2)}$ be the result of this quasi-minimization. Similar extension is possible even when $\dim V\cap \ell ^2(\mathbb{Z}^+ )\ge 3$, based on the idea of recursive iteration of the quasi-minimization of this ratio  in the subspace almost orthogonal to the span of the already obtained vectors $\vec{G}^{\, (1)},\vec{G}^{\, (2)},\ldots \,\,$. Thus we can obtain as many linearly independent vectors as $\dim V\cap \ell ^2(\mathbb{Z}^+ )$, each of which belongs to $((\sigma_{K,N}^{(\Omega )})^{-1}[0,c \underline{\sigma_{K,N}^{(\Omega )}} ])$ and hence is close to $V\cap \ell ^2(\mathbb{Z}^+ )$. This provides us with an approximate quasi-orthogonal  basis system of $V\cap \ell ^2(\mathbb{Z}^+ )$ with respect to $\langle \cdot , \cdot \rangle _{\ell ^2,K}$, i.e., with an approximation of the `general solution' of the differential equation.

Here, we should be careful of the fact that a linear combination of a set of vectors is not always close to $V\cap \ell ^2(\mathbb{Z}^+ )$ even if all the vectors in this set are close to $V\cap \ell ^2(\mathbb{Z}^+ )$ (in the sense of the angles between the vectors and the space). However, as is shown later in Section \ref{sec:conv_multidim}, when the set of vectors form a quasi-orthogonal system with respect to  $\langle \cdot , \cdot \rangle _{\ell ^2,K}$, any nonzero linear combination of them is close to $V\cap \ell ^2(\mathbb{Z}^+ )$. Hence, the idea mentioned above for the extension does not suffer from this problem.

In this section and the next one, we will explain the details of the extension. The procedures {\bf Step 0}-{\bf Step 2} do not change at all with this extension, because a basis system of $\Pi _N V$ is required there in the same sense as the one-dimensional case. Hence, we have only to explain how to modify {\bf Step 3} for this extension. In this section, we explain only how to modify the procedures; their validity will be proved in Section \ref{sec:conv_multidim}.

\begin{table}[p]
\caption{Basic idea of integer-valued quasi-orthogonalization (general)}
\begin{center}
\scalebox{.86}{
\begin{tabular}{|l@{}|}
\hline 
$ \begin{array}{@{\,}ll}  \\ 
\mbox{{\bf if} $\tilde{d}>0$} \\ 
\,\,\, \left\lfloor\begin{array}{@{\,}ll} 
\bullet \, 
\mbox{{\bf Iterate S1} below until nothing is changed.} \\  
\,\,\,\,\,\, \left\lfloor\begin{array}{@{\,}ll}
\mbox{{\bf S1}} \left\{\begin{array}{@{\,}ll} 
\bullet \mbox{\, sorting and renumbering of $\vec{v}_{\tilde{d}+1}, \, \vec{v}_{\tilde{d}+2}, \, \ldots \, , \, \vec{v}_D$} \\ 
\mbox{ in order that $ \| \vec{v}_{\tilde{d}+1}\|_{\ell ^2,K} \le \| \vec{v}_{\tilde{d}+2}\|_{\ell ^2,K} \le \ldots \le \| \vec{v}_D\|_{\ell ^2,K}$} \\ \\ 
\bullet \mbox{\, {\bf for} $j=\tilde{d}+2$ to $D$} \\ 
\,\,\,\,\,
\,\,\, \left\lfloor\begin{array}{@{\,}ll}  \mbox{{\bf for} $\ell =1$ to $j-1$} \\  
\,\,\, \left\lfloor\begin{array}{@{\,}ll} 
\mbox{$\displaystyle \vec{v}_j\leftarrow\vec{v}_j-\Bigl[\frac{\langle\vec{v}_j, \vec{v}_\ell \rangle_{\ell ^2,K}}{\langle\vec{v}_\ell , \vec{v}_\ell \rangle_{\ell ^2,K}} \Bigr]_\mathbb{C} \vec{v}_\ell $} \\ 
\end{array}\right. \\ 
\end{array}\right. \\ 
\end{array}\right. \\ 
\end{array}\right. \\    
\,\,\,\,\,\,\,\,\,\,\,\, \left(\begin{array}{@{\,}ll}
\mbox{This is a preparatory partial lattice reduction for {\bf Q2} below.}
\end{array}\right) 
\\ \hline 
\bullet 
\, \mbox{{\bf Iterate S2} below until nothing is changed.} \\  
\,\,\,\,\,\,  \left\lfloor\begin{array}{@{\,}ll}
\mbox{{\bf S2}} 
\left\{\begin{array}{@{\,}ll} \mbox{{\bf for} $j=\tilde{d}+2$ to $D$} \\ 
\,\,\, \left\lfloor\begin{array}{@{\,}ll}
\mbox{{\bf for} $\ell =1$ to $j-1$} \\  
\,\,\, \left\lfloor\begin{array}{@{\,}ll} 
\mbox{{\bf if} $\,\,\, g^2\,|\langle \vec{v}_j, \vec{v}_\ell \rangle_{\ell ^2,K}|^2 \ge %\bigl(\epsilon(\zeta,\, D-\tilde{d}+1)\bigr)^2 \,
 \|\vec{v}_j\|_{\ell ^2,K}^2 \|\vec{v}_\ell \|_{\ell ^2,K}^2 $,} \\  
\,\,\, \left\lfloor\begin{array}{@{\,}ll} 
\mbox{$\vec{v}_j \leftarrow 2 \vec{v}_j$} \\
\end{array}\right. \\ 
\mbox{$\displaystyle \vec{v}_j\leftarrow\vec{v}_j-\Bigl[\frac{\langle\vec{v}_j, \vec{v}_\ell \rangle_{\ell ^2,K}}{\langle\vec{v}_\ell , \vec{v}_\ell \rangle_{\ell ^2,K}} \Bigr]_\mathbb{C} \vec{v}_\ell $}\\ 
\end{array}\right. \\ 
\end{array}\right. \\  
\end{array}\right. \\  
\end{array}\right.  \\  
\,\,\,\,\,\,\,\,\,\,\,\, \left(\begin{array}{@{\,}ll}
\mbox{This process provides us with the vectors $\vec{v}_j$ $(1\le j\le D)$} \\ 
\mbox{ s.t. $\displaystyle\frac{|\langle \vec{v}_j, \vec{v}_\ell \rangle _{\ell ^2,K} |}{\|\vec{v}_j\|_{\ell ^2,K} \cdot\|\vec{v}_\ell \|_{\ell ^2,K}} \le \frac{1}{g} $}\\ %\epsilon _{\zeta,D-\tilde{d}+1}\,$} \\ 
\mbox{ for $\, \tilde{d}+1 \le j<\ell \le D$\, or $1\le j < \tilde{d}+1 \le \ell \le D$.} \\ 
\end{array}\right) \\  
\end{array}\right.  \\  
\end{array} $ \\
\hline 
$\begin{array}{@{\,}ll} \\   
\mbox{{\bf Iterate Q$^\prime$1} below until nothing is changed.} \\ 
\,\,\, \left\lfloor\begin{array}{@{\,}ll} 
\mbox{{\bf Q$^\prime$1}} \left\{\begin{array}{@{\,}ll} 
\bullet \mbox{\, sorting and renumbering  of $\vec{v}_{\tilde{d}+1}, \, \vec{v}_{\tilde{d}+2}, \, \ldots \, , \, \vec{v}_D$} \\ 
\,\,\,\,\, \mbox{ in order that $ \| \vec{v}_{\tilde{d}+1}\|_{Q,N} \le \| \vec{v}_{\tilde{d}+2}\|_{Q,N} \le \ldots \le \| \vec{v}_D\|_{Q,N} $} \\  
\bullet \mbox{\, {\bf for} $j=\tilde{d}+2$ to $D$} \\ 
\,\,\,\,\, 
\,\,\, \left\lfloor\begin{array}{@{\,}ll}
\mbox{{\bf for} $\ell =\tilde{d}+1$ to $j-1$} \\  
\,\,\, \left\lfloor\begin{array}{@{\,}ll}
\mbox{$\displaystyle \vec{v}_j\leftarrow\vec{v}_j-\Bigl[\frac{\langle\vec{v}_j, \vec{v}_\ell \rangle_{Q,N}}{\langle\vec{v}_\ell , \vec{v}_\ell \rangle_{Q,N}} \Bigr]_\mathbb{C} \vec{v}_\ell $} \\ 
\end{array}\right. \\ 
\end{array}\right. \\ 
\end{array}\right. \\  
\end{array}\right. \\  
\,\,\,\,\,\,\,\,\,\,\,\, \left(\begin{array}{@{\,}ll}
\mbox{This is a preparatory partial lattice reduction for {\bf P2} below.}
\end{array}\right)
\end{array}$ \\  
\hline 
$\begin{array}{@{\,}ll} \\ 
\mbox{{\bf Iterate Q$^\prime$2} below until nothing is changed.} \\ 
\,\,\, \left\lfloor\begin{array}{@{\,}ll} 
\mbox{{\bf Q$^\prime$2}} 
\left\{\begin{array}{@{\,}ll} \mbox{\, {\bf for} $j=\tilde{d}+2$ to $D$} \\ 
\,\,\, \left\lfloor\begin{array}{@{\,}ll} 
\mbox{{\bf for} $\ell =\tilde{d}+1$ to $j-1$} \\  
\,\,\, \left\lfloor\begin{array}{@{\,}ll}  
\mbox{{\bf if} $\,\,\, h^2\, |\langle \vec{v}_j, \vec{v}_\ell \rangle_{Q,N}|^2 \ge %\bigl(\epsilon(\delta,D-\tilde{d})\bigr)^2\, 
\|\vec{v}_j\|_{Q,N}^2 \|\vec{v}_\ell \|_{Q,N}^2 $,} \\  
\, \left\lfloor\begin{array}{@{\,}ll} 
\,\,\,\mbox{$\vec{v}_j \leftarrow 2 \vec{v}_j$} \\
\end{array}\right. \\ 
\mbox{$\displaystyle \vec{v}_j\leftarrow\vec{v}_j-\Bigl[\frac{\langle\vec{v}_j, \vec{v}_\ell \rangle_{Q,N}}{\langle\vec{v}_\ell , \vec{v}_\ell \rangle_{Q,N}} \Bigr]_\mathbb{C} \vec{v}_\ell $} \\ 
\end{array}\right. \\ 
\end{array}\right. \\ 
\end{array}\right. \\ 
\end{array}\right. \\  
\,\,\,\,\,\,\,\,\,\,\,\, \left(\begin{array}{@{\,}ll}
\mbox{This process provides us with the vectors $\vec{v}_j$ $(1\le j \le D)$} \\ 
\mbox{ s.t. $\displaystyle\frac{|\langle \vec{v}_j, \vec{v}_\ell \rangle _{Q,N} |}{\|\vec{v}_j\|_{Q,N} \cdot\|\vec{v}_\ell \|_{Q,N}} \le \frac{1}{h}%\epsilon (\delta,D-\tilde{d})\,
 $ for $\, \tilde{d}+1 \le j<\ell \le D$\, and } \\ 
\mbox{$\displaystyle\frac{|\langle \vec{v}_j, \vec{v}_\ell \rangle _{\ell ^2,K} |}{\|\vec{v}_j\|_{\ell ^2,K} \cdot\|\vec{v}_\ell \|_{\ell ^2,K}} \le \frac{1}{g} \, \sqrt{\frac{D-\tilde{d}}{1-\frac{\displaystyle D-\tilde{d}-1}{\displaystyle g}}} $ for $1\le j < \tilde{d}+1 \le \ell \le D$.} \\ 
\end{array}\right)   
\end{array}$ \\ 
\hline
\end{tabular}
}
\end{center}
\label{tbl:4}
\end{table}

\begin{table}[p]
\caption{Practical operations for integer-type quasi-orthogonalization (general)}
\begin{center}
\scalebox{.86}{
\begin{tabular}{|l@{}|}
\hline 
$\begin{array}{@{\,}ll}  
\mbox{{\bf for} $j=1$ to $D$} \\  
\,\,\, \left\lfloor\begin{array}{@{\,}ll} 
\mbox{{\bf for} $m=1$ to $D$} \\  
\,\,\, \left\lfloor\begin{array}{@{\,}ll} 
\mbox{$c_{jm} \leftarrow \delta _{jm}$} \\ 
\mbox{$p_{jm} \leftarrow \langle \vec{v}_j^{\, ({\rm initial})} , \vec{v}_{m}^{\, ({\rm initial})} \rangle _{Q,N}$} \\
\mbox{$q_{jm} \leftarrow \langle \vec{v}_j^{\, ({\rm initial})} , \vec{v}_{m}^{\, ({\rm initial})} \rangle _{\ell^2,K}$} \\
\end{array}\right. \\ 
\end{array}\right. \\ 
\end{array} $ \\ 
\,\,\,\,\,\,\,\,\, \mbox{$\left( \, \begin{array}{@{\,}ll} c_{jm}: \mbox{complex-integer-valued coefficients} \\ \displaystyle 
\,\,\,\,\,\,\,\,\, \,\,\,\,\,\, \mbox{ for } \,\,\,\,\, \vec{v}_j^{\, ({\rm new})}=\sum_{m=1}^D c_{jm} \vec{v}_m^{\, ({\rm initial})} \\ 
p_{jm},q_{jm}: \mbox{ inner products } \\ 
\end{array} \, \right) $} \\ 
\hline 
$\begin{array}{@{\,}ll} \\ 
\mbox{{\bf if} $\tilde{d}>0$} \\ 
\,\,\, \left\lfloor\begin{array}{@{\,}ll} 
\bullet \, \mbox{{\bf Iterate R1} below until nothing is changed.} \\  
\,\,\,\,\,\,  \left\lfloor\begin{array}{@{\,}ll}
\mbox{{\bf R1}} \left\{\begin{array}{@{\,}ll} \bullet \mbox{\, with a permutation $\,  n_{\tilde{d}+1}, n_{\tilde{d}+2},\ldots,n_D \,$ of } \\ 
\,\,\, \mbox{$\, \tilde{d}+1,\tilde{d}+2,\ldots,D $\, s.t. $\displaystyle \, q_{\tilde{d}+1\, \tilde{d}+1} \le q_{\tilde{d}+2\, \tilde{d}+2} \le \ldots \le q_{D D}$}, \\ 
\,\,\,\,\,\,\,\,\,\,\,\,\,
\mbox{{\bf for} $j=1$ to $D$} \\ 
\,\,\,\,\,\,\,\,\,\,\,\,\,
\,\,\left\lfloor\begin{array}{@{\,}ll} 
\mbox{$c_{jm} \leftarrow c_{n_j m}, \,\,\, p_{jm} \leftarrow p_{n_j n_m}, \,\,\, q_{jm} \leftarrow p_{n_j n_m}$} \\ 
\end{array}\right.\\  \\ 
\bullet \mbox{\, {\bf for} $j=\tilde{d}+2$ to $D$} \\ 
\,\,\,\,\,
\,\,\, \left\lfloor\begin{array}{@{\,}ll}
\mbox{{\bf for} $\ell =1$ to $j-1$} \\  
\,\,\,\,\,\, \left\lfloor\begin{array}{@{\,}ll} 
\mbox{$\displaystyle r\leftarrow \Bigl[ \frac{q_{j\ell }}{q_{\ell \ell}} \Bigr]_\mathbb{C} $} \\ 
\mbox{{\bf for} $m=1$ to $D$} \\  
\,\,\, \left\lfloor\begin{array}{@{\,}ll}
\mbox{$\displaystyle c_{jm}\leftarrow c_{jm}-r\, c_{\ell m} $} \\ 
\mbox{$\displaystyle p_{jm}\leftarrow p_{jm}-r\, p_{\ell m} , \,\,\,  q_{jm}\leftarrow q_{jm}-r\, q_{\ell m} $} \\ 
\end{array}\right.\\ 
\mbox{{\bf for} $m=1$ to $D$} \\  
\,\,\, \left\lfloor\begin{array}{@{\,}ll}
\mbox{$\displaystyle p_{mj}\leftarrow p_{mj}-\overline{r}\, p_{m\ell } , \,\,\,  q_{mj}\leftarrow q_{mj}-\overline{r}\, q_{m\ell } $} \\ 
\end{array}\right. \\
\end{array}\right. \\ 
\end{array}\right. \\ 
\end{array}\right. \\  \end{array}\right.\\ 
\hline  \\ 
\bullet \, \mbox{{\bf Iterate R2} below until nothing is changed} \\  
\,\,\,\,\,\,  \left\lfloor\begin{array}{@{\,}ll} 
\mbox{{\bf R2}}  
\left\{\begin{array}{@{\,}ll}
\mbox{{\bf for} $j=\tilde{d}+2$ to $D$} \\ 
\,\,\, \left\lfloor\begin{array}{@{\,}ll}
\mbox{{\bf for} $\ell =1$ to $j-1$} \\  
\,\,\, \left\lfloor\begin{array}{@{\,}ll} 
\mbox{{\bf if} $\displaystyle \,\,\, g^2\left[\frac{q_{j\ell }}{\widehat{N}^K}\right]_\mathbb{C}  \left[\frac{\overline{q_{j\ell }}}{\widehat{N}^K}\right]_\mathbb{C}\! \ge \left(\Bigl\lfloor\frac{q_{jj}}{\widehat{N}^K}\Bigr\rfloor \!+\! 1\right) \left(\Bigl\lfloor\frac{q_{\ell\ell }}{\widehat{N}^K}\Bigr\rfloor \! +\! 1\right)$,} \\ 
\, \left\lfloor\begin{array}{@{\,}ll}
\mbox{{\bf for} $m=1$ to $D$} \\ 
\,\,\, \left\lfloor\begin{array}{@{\,}ll}
\mbox{$p_{mj} \leftarrow 2 p_{mj}$,\,\, $q_{mj} \leftarrow 2 q_{mj}$} \\
\mbox{$p_{jm} \leftarrow 2 p_{jm}$,\,\, $q_{jm} \leftarrow 2 q_{jm}$} \\
\end{array}\right. \\ 
\end{array}\right. \\  
\mbox{$\displaystyle r\leftarrow \Bigl[ \frac{q_{j\ell }}{q_{\ell \ell}} \Bigr]_\mathbb{C} $} \\ 
\mbox{{\bf for} $m=1$ to $D$} \\  
\,\,\, \left\lfloor\begin{array}{@{\,}ll}
\mbox{$\displaystyle c_{jm}\leftarrow c_{jm}-r\, c_{\ell m} $} \\ 
\mbox{$\displaystyle p_{jm}\leftarrow p_{jm}-r\, p_{\ell m} $} \\ 
\mbox{$\displaystyle q_{jm}\leftarrow q_{jm}-r\, q_{\ell m} $} \\ 
\end{array}\right.\\ 
\mbox{{\bf for} $m=1$ to $D$} \\  
\,\,\, \left\lfloor\begin{array}{@{\,}ll}
\mbox{$\displaystyle p_{mj}\leftarrow p_{mj}-\overline{r}\, p_{m\ell } $} \\ 
\mbox{$\displaystyle q_{mj}\leftarrow q_{mj}-\overline{r}\, q_{m\ell } $} \\ 
\end{array}\right. \\
\end{array}\right. \\ 
\end{array}\right. \\ 
\end{array}\right. \\ 
\end{array}\right. \\ 
\end{array}\right. \\  \end{array} $
\\ 
\hline 
%\mbox{} \\ 
\multicolumn{1}{|c|}{(to be continued to the next page)}
\end{tabular}
}
\end{center}
\label{tbl:5}
\end{table}
\begin{table}[p]
\begin{center}
\scalebox{.86}{
\begin{tabular}{|l@{}|}
\multicolumn{1}{|c|}{(continued from the previous page)} \\ 
%\mbox{} \\ 
\hline   
$\begin{array}{@{\,}ll} \\ 
\mbox{{\bf Iterate P$^\prime$1} below until nothing is changed.} \\  
\,\,\, \left\lfloor\begin{array}{@{\,}ll}
\mbox{{\bf P$^\prime$1}} \left\{\begin{array}{@{\,}ll} \bullet \mbox{\, with a permutation $\,  n_{\tilde{d}+1}, n_{\tilde{d}+2},\ldots,n_D$\,  of } \\ 
\,\,\, \mbox{$\, \tilde{d}+1,\tilde{d}+2,\ldots,D $\, s.t. $\displaystyle \, p_{\tilde{d}+1\, \tilde{d}+1} \le p_{\tilde{d}+2\, \tilde{d}+2} \le \ldots \le p_{D D}$}, \,\,\,\,\,\,\,\,\,\,\,\,\,\,\,\, \\ 
\,\,\,\,\,\,\,\,
\mbox{{\bf for} $j=1$ to $D$} \\ 
\,\,\,\,\,\,\,\,
\,\,\left\lfloor\begin{array}{@{\,}ll} 
\mbox{$c_{jm} \leftarrow c_{n_j m} ,\,\, p_{jm} \leftarrow p_{n_j n_m} ,\,\,\, q_{jm} \leftarrow p_{n_j n_m}$} \\ 
\end{array}\right.  \\ 
\bullet \mbox{\, {\bf for} $j=\tilde{d}+2$ to $D$} \\ 
\,\,\,\,\, 
\,\,\, \left\lfloor\begin{array}{@{\,}ll}
\mbox{{\bf for} $\ell =\tilde{d}+1$ to $j-1$} \\  
\,\,\, \left\lfloor\begin{array}{@{\,}ll} 
\mbox{$\displaystyle r\leftarrow \Bigl[ \frac{p_{j\ell }}{p_{\ell \ell}} \Bigr]_\mathbb{C} $} \\ 
\mbox{{\bf for} $m=1$ to $D$} \\  
\,\,\, \left\lfloor\begin{array}{@{\,}ll}
\mbox{$\displaystyle c_{jm}\leftarrow c_{jm}-r\, c_{\ell m} $} \\ 
\mbox{$\displaystyle p_{jm}\leftarrow p_{jm}-r\, p_{\ell m} ,\,\,\,  q_{jm}\leftarrow q_{jm}-r\, q_{\ell m} $} \\ 
\end{array}\right.\\ 
\mbox{{\bf for} $m=1$ to $D$} \\  
\,\,\,\left\lfloor\begin{array}{@{\,}ll}
\mbox{$\displaystyle p_{mj}\leftarrow p_{mj}-\overline{r}\, p_{m\ell } ,\,\,\, q_{mj}\leftarrow q_{mj}-\overline{r}\, q_{m\ell } $} \\ 
\end{array}\right. \\
\end{array}\right. \\ 
\end{array}\right. \\ 
\end{array}\right. \\  
\end{array}\right. \\  
\end{array}$
\\ 
\hline 
$\begin{array}{@{\,}ll} \\ 
\mbox{{\bf Iterate P$^\prime$2} below until nothing is changed.} \\  
\,\,\, \left\lfloor\begin{array}{@{\,}ll}  
\mbox{{\bf P$^\prime$2}} \left\{ \begin{array}{@{\,}ll} 
\mbox{\, {\bf for} $j=\tilde{d}+2$ to $D$} \\ 
%\,\,\, \left\lfloor\begin{array}{@{\,}ll}
\mbox{{\bf for} $\ell =\tilde{d}+1$ to $j-1$} \\  
\,\,\, \left\lfloor\begin{array}{@{\,}ll} 
\mbox{{\bf if} $\displaystyle \,\,\, h^2\,\left[\frac{p_{j\ell }}{\widehat{N}^K}\right]_\mathbb{C}  \left[\frac{\overline{p_{j\ell }}}{\widehat{N}^K}\right]_\mathbb{C} \ge \left(\Bigl\lfloor\frac{p_{jj}}{\widehat{N}^K}\Bigr\rfloor +1\right)\left(\Bigl\lfloor\frac{p_{\ell\ell }}{\widehat{N}^K}\Bigr\rfloor +1\right)$,} \\ 
\, \left\lfloor\begin{array}{@{\,}ll} 
\mbox{{\bf for} $m=1$ to $D$} \\ 
\,\,\, \left\lfloor\begin{array}{@{\,}ll}
\mbox{$p_{mj} \leftarrow 2 p_{mj}$,\,\, $q_{mj} \leftarrow 2 q_{mj}$} \\
\mbox{$p_{jm} \leftarrow 2 p_{jm}$,\,\, $q_{jm} \leftarrow 2 q_{jm}$} \\
\end{array}\right. \\ 
\end{array}\right. \\ 
\mbox{$\displaystyle r\leftarrow \Bigl[ \frac{p_{j\ell }}{p_{\ell \ell}} \Bigr]_\mathbb{C} $} \\ 
\mbox{{\bf for} $m=1$ to $D$} \\  
\,\,\, \left\lfloor\begin{array}{@{\,}ll}
\mbox{$\displaystyle c_{jm}\leftarrow c_{jm}-r\, c_{\ell m} $} \\ 
\mbox{$\displaystyle p_{jm}\leftarrow p_{jm}-r\, p_{\ell m} ,\,\,\, q_{jm}\leftarrow q_{jm}-r\, q_{\ell m} $} \\ 
\end{array}\right.\\ 
\mbox{{\bf for} $m=1$ to $D$} \\  
\,\,\, \left\lfloor\begin{array}{@{\,}ll}
\mbox{$\displaystyle p_{mj}\leftarrow p_{mj}-\overline{r}\, p_{m\ell } ,\,\,\,  q_{mj}\leftarrow q_{mj}-\overline{r}\, q_{m\ell } $} \\ 
\end{array}\right. \\
\end{array}\right. \\ 
\end{array}\right. \\ 
\end{array}\right. \\
\end{array}$ \\ 
\hline 
$\begin{array}{@{\,}ll}  
\mbox{{\bf for} $j=\tilde{d}+1$ to $D$} \\ 
\,\,\, \left\lfloor\begin{array}{@{\,}ll} 
\mbox{$\displaystyle \vec{v}_j^{\rm \, (final)} \leftarrow \sum_{m=1}^D c_{jm} \vec{v}_m^{\rm \, (initial)}  $} \\ 
\mbox{{\bf for} $\ell =1$ to $D$} \\  
\,\,\, \left\lfloor\begin{array}{@{\,}ll}
\mbox{{\bf if} $j\ge \tilde{d}+1$ or $m\ge \tilde{d}+1$, } \\ 
\, \left\lfloor\begin{array}{@{\,}ll}
\mbox{$\displaystyle \langle \vec{v}_j^{\rm \, (final)} , \vec{v}_\ell^{\rm \, (final)} \rangle _{Q,N} \leftarrow \sum_{m=1}^D \sum_{n=1}^D c_{jm} \overline{c_{\ell n}} \, p_{mn}$} \\ 
\mbox{$\displaystyle \langle \vec{v}_j^{\rm \, (final)} , \vec{v}_\ell^{\rm \, (final)} \rangle _{\ell^2,K} \leftarrow \sum_{m=1}^D \sum_{n=1}^D c_{jm} \overline{c_{\ell n}} \, q_{mn}$}  \\
\end{array}\right. \\
\end{array}\right. \\
\end{array}\right. \\ 

\mbox{{\bf for} $j=1$ to $D$} \\ 
\,\,\, \left\lfloor\begin{array}{@{\,}ll} 
\mbox{{\bf for} $d=1$ to $D$} \\  
\,\,\, \left\lfloor\begin{array}{@{\,}ll}
\mbox{$\displaystyle a_{j,d}^{\rm \, (final)} \leftarrow \sum_{m=1}^D  c_{jm} a_{m,d}^{\rm \, (initial)}
 $} 
\\ 
\end{array}\right.  \\
\end{array}\right. \,\,\, \left(\begin{array}{@{\,}ll}\displaystyle a_{j,d}:\mbox{complex-integer-valued} \\ \displaystyle\,\,\,\,\,\,\,\,\,\,\,\,  \mbox{ coefficients s.t. } \vec{v}_j = \sum_{d=1}^D a_{j,d} \vec{F}_{\rm int.}^{\, (d)} \!\!\! \end{array}\right) \\ 

\end{array} $ \\  
\hline
\end{tabular}
}
\end{center}
\end{table}

For a concrete description of this, we provide some preliminary notation. 
Let $\tilde{d}$ be the number of `already obtained vectors' $\vec{G}^{\, (d)}$ $(d=1,2, \ldots ,\tilde{d})$ by the method mentioned above, and let $T_{\tilde{d}}:=< \vec{G}^{\, (1)}, \vec{G}^{\, (2)},\ldots , \vec{G}^{\, (\tilde{d})}>$ where we set $T_0:=\{ 0\}$. 
Let $R_{\tilde{d}}$ denote the subspace of $\Pi _N V$ satisfying $\Pi _N V=T_{\tilde{d}}\oplus R_{\tilde{d}}$, in which we find another quasi-minimum-ratio vector at the next step.  By the procedures in {\bf Step 3$^\prime $.A.1}  below, 
this subspace   $R_{\tilde{d}}$ is chosen 
to be very close to the orthogonal complement $\Pi _N V^{\bot T_{\tilde{d}}}$ of $T_{\tilde{d}}$
%such that $\Pi _N V=T_{\tilde{d}}\oplus \Pi _N V^{\bot T_{\tilde{d}}}$ and $T_{\tilde{d}}\bot \Pi _N V^{\bot T_{\tilde{d}}}$ 
with respect to $\langle \cdot , \cdot \rangle _{\ell ^2,K}$, but it is not always exactly equal to the latter, where `closeness' is used  in the sense that the angles between nonzero vectors in $T_{\tilde{d}}$ and nonzero vectors in $R_{\tilde{d}}$ are close to $\displaystyle\frac{\pi}{2}$ with respect to $\langle \cdot , \cdot \rangle _{\ell ^2,K}$. Obviously, $\dim T_{\tilde{d}}=\tilde{d}$ and $\dim R_{\tilde{d}}=N-\tilde{d}$. 
Moreover, let $D_{\ell ^2}=\dim V\cap \ell ^2(\mathbb{Z}^+ )$.

With these notations, 
for the extension to the cases where $D_\ell ^2\ge 2$, the procedures in {\bf Step3} are replaced by 

\vspace{2mm}

\begin{description}
\item[Step 3$^\prime $] Removal of components from $\Pi _N V$ corresponding to non-$\ell ^2$-ones in $\Pi _N V$: 
\begin{description}
\item[Step 3$^\prime $.A] Integer-type extraction of a quasi-orthogonal basis system for $\Pi _N (V \cap {\ell ^2}(\mathbb{Z}^+))$:
\item[]
\hspace{5mm} Iterate the series of steps {\bf Step 3$^\prime $.A.a1}$-${\bf Step  3$^\prime $.A.b2} below (once in this order) for $\tilde{d}=0,1, ,\ldots ,D_{\ell ^2}-1$ with %the initial basis vectors $\vec{V}_d^{<0>}:=\vec{F}_{\rm int.}^{\, (d)}$ $(d=1,2,\ldots ,D)$ and 
the initial subspaces $T_0=\{ 0\}$. (If we require an upper bound for errors, iterate this for $\tilde{d}=0,1,,\ldots ,D_{\ell ^2}$.)\, The result of each step of this iteration gives a quasi-orthogonal basis system  $\vec{G}^{\, (1)}, \vec{G}^{\, (2)},\dots \vec{G}^{\, (\tilde{d})},$
%$\vec{V}_1^{<\tilde{d}>},\vec{V}_2^{<\tilde{d}>},\ldots ,\vec{V}_{\tilde{d}}^{<\tilde{d}>}$ 
of $T_{\tilde{d}}$ with respect to $\langle \cdot , \cdot \rangle _{\ell ^2,K}$.
%, as well as  $\vec{V}_1^{<\tilde{d}>},\vec{V}_2^{<\tilde{d}>},\ldots ,\vec{V}_{D}^{<\tilde{d}>}$ is a basis system of $\Pi _N V$. 
(This substep is an extension of {\bf Step 3.1.a}+{\bf Step 3.1.b}. Each iteration with $\tilde{d}=n$ of this substep is corresponding to {\bf Step 3.$n$} of~\cite{paper1}, with a slight difference between exact orthogonalization and quasi-orthogonalization with respect to $\langle \cdot ,\,\cdot \rangle _{\ell ^2,K}$.)
\begin{description}
\item[Step 3$^\prime $.A.a1] Integer-type quasi-orthogonalization of the basis system of $R_{\tilde{d}} $:
\item[]
\hspace{5mm} If $\tilde{d}\ge 1$, find a system of $N-\tilde{d}$ linear combinations $\vec{v}_{1}^{<\tilde{d}>},\vec{v}_{2}^{<\tilde{d}>},\ldots ,\vec{v}_{N-\tilde{d}}^{<\tilde{d}>}$ of $\vec{F}_{\rm int.}^{\, (1)}, \vec{F}_{\rm int.}^{\, (2)},\ldots , \vec{F}_{\rm int.}^{\, (D)}$ %$\vec{V}_{1}^{<\tilde{d}>},\vec{V}_{2}^{<\tilde{d}>},\ldots ,\vec{V}_{D}^{<\tilde{d}>}$ 
which is sufficiently close to an orthogonal system with respect to the inner product $\langle \cdot ,\,\cdot \rangle _{\ell ^2,K}$ and also almost orthogonal to $T_{\tilde{d}}$ 
%the `already obtained quasi-minimum-ratio vectors' $\vec{G}^{\, (d)}$ $(d=1,2, \ldots , \tilde{d})$  
with respect to  $\langle \cdot ,\,\cdot \rangle _{\ell ^2,K}$, %by the iterations of the procedures {\bf S$^\prime $1} and {\bf S$^\prime $2}  %given in Table 5.1 
by the procedures explained below which are based on an idea intermediate between the Gram-Schmidt process and the Euclidean algorithm. 
% There, the vectors $\vec{V}_1^{<\tilde{d}>},\vec{V}_2^{<\tilde{d}>},\ldots ,\vec{V}_{\tilde{d}}^{<\tilde{d}>}$ are unchanged. 
Then, if $\tilde{d}\ge 1$, choose the subspace $R_{\tilde{d}}$ to be the span of the $N-\tilde{d}$ linear combinations $\vec{v}_{1}^{<\tilde{d}>},\vec{v}_{2}^{<\tilde{d}>},\ldots ,\vec{v}_{N-\tilde{d}}^{<\tilde{d}>}$  obtained by this substep. This substep is omitted by the simple substitution $\vec{v}_{d}^{<0>}:=\vec{F}_{\rm int.}^{\, (d)}$ $(d=1,2,\ldots ,D)$ and   $R_0:=\Pi _N V$ when $\tilde{d}=0$.
\item[Step 3$^\prime $.A.a2] Integer-type quasi-orthogonalization of basis system of $R_{\tilde{d}}$:
\item[]
\hspace{5mm} Find a system of $N-\tilde{d}$ linear combinations $\vec{u}_{1}^{<\tilde{d}>},\vec{u}_{2}^{<\tilde{d}>},\ldots ,\vec{u}_{N-\tilde{d}}^{<\tilde{d}>}$  of the $N-\tilde{d}$ basis vectors $\vec{v}_{1}^{<\tilde{d}>},\vec{v}_{2}^{<\tilde{d}>},\ldots ,\vec{v}_{N-\tilde{d}}^{<\tilde{d}>}$ of $R_{\tilde{d}}$ obtained in {\bf Step 3$^\prime $.A.a1} 
% of $\vec{V}_{\tilde{d}+1}^{<\tilde{d}>},\vec{V}_{\tilde{d}+2}^{<\tilde{d}>},\ldots ,\vec{V}_{D}^{<\tilde{d}>}$ in $R_{\tilde{d}}$ 
which is sufficiently close to an orthogonal system with respect to the inner product $\langle \cdot ,\,\cdot \rangle _{Q,N}$, %by the iterations of the procedures {\bf Q$^\prime $1} and {\bf Q$^\prime $2} %given in Table 5.1 
by the procedures explained below which is based on an intermediate idea between the Gram-Schmidt process and the Euclidean algorithm. 
% There, the vectors $\vec{V}_1^{<\tilde{d}>},\vec{V}_2^{<\tilde{d}>},\ldots ,\vec{V}_{\tilde{d}}^{<\tilde{d}>}$ are unchanged.
\item[Step 3$^\prime $.A.b1] Selection of minimum-ratio vector:
\item[]
\hspace{5mm} Find the linear combination $\vec{G}^{\, (\tilde{d}+1)}$ with minimum ratio $\displaystyle\frac{\|\cdot\|_{Q,N}}{\|\cdot\|_{\ell ^2,K}}$ in the $N-\tilde{d}$ linear combinations $\vec{u}_{1}^{<\tilde{d}>},\vec{u}_{2}^{<\tilde{d}>},\ldots ,\vec{u}_{N-\tilde{d}}^{<\tilde{d}>}$ obtained by {\bf Step 3$^\prime $.A.a2}. 
\item[Step 3$^\prime $.A.b2] Innovation of space $T_{\tilde{d}}\, $:
\item[]
\hspace{5mm} 
Let $T_{\tilde{d}+1}:= T_{\tilde{d}}\oplus \{ a\vec{G}^{\, (\tilde{d}+1)}\, |\, a\in \mathbb{C} \} $ with $\vec{G}^{\, (\tilde{d}+1)}$ obtained by {\bf Step 3$^\prime $.A.b1}.
\item[]
\hspace{5mm} 
\end{description}
\item[Step 3$^\prime $.B.c] Truncation (projection) to $\Pi _K V$:
\item[]
\hspace{5mm}
Project the vectors $\vec{G}^{\, (d)}$ $(d=1,2, \ldots , D_{\ell ^2})$ obtained by {\bf Step 3$^\prime $.A} by $\Pi _K$. 
\end{description}
%\hspace{5mm} Find %the vector $\vec{f}_{K,N}^{(Q)}$ in (\ref{eqn:def_ratio2}), or find 
% a vector in the set $O_{K,N}^{(Q,c)}$ in {\bf C8} (which is identical to $O_{K,N}^{(Q,c)}$ in Condition {\bf C8}) from the basis vectors $\{F_n^{(1)}\}_{n=0}^{N} ,\ldots, \{F_n^{(D)}\}_{n=0}^{N}$ of $\Pi _N V$.% instead. As will be shown later, the latter alternative choice spares the quantity of calculation remarkably. 
\end{description}
\vspace{3mm}

In the following, we explain how to realize these procedures, in detail. These procedures can be described in a unified framework of an iterative change of the basis systems of $\Pi _N V$. In this framework, 
all the basis vectors of $\Pi _N V$ 
at the intermediate steps with $\tilde{d}=0,1,\ldots , D_{\ell ^2}$ in {\bf Step 3$^\prime $.A} above are denoted by $\vec{V}_d^{<\tilde{d}>}$  $(d=1,2,\ldots ,D)$. With these notations, the basis vectors of $\Pi _N V$ just after {\bf Step 3$^\prime $.A.a1} with $\tilde{d}$ are 
\begin{eqnarray*}
\vec{V}_d^{<\tilde{d}>}=\left\{ \begin{array}{@{\,}ll} \displaystyle \vec{G}^{\, (d)} & ( {\rm if }\, d\le \tilde{d}) \\ \\ \displaystyle \vec{v}_{d-\tilde{d}}^{<\tilde{d}>} & ( {\rm if }\, d \ge  \tilde{d}+1) \,\, , \end{array}\right. 
\end{eqnarray*}
while they are changed to 
\begin{eqnarray*}
\vec{V}_d^{<\tilde{d}>}=\left\{ \begin{array}{@{\,}ll} \displaystyle \vec{G}^{\, (d)} & ( {\rm if }\, d\le \tilde{d}) \\ \\ \displaystyle \vec{u}_{d-\tilde{d}}^{<\tilde{d}>} & ( {\rm if }\, d \ge  \tilde{d}+1) \end{array}\right. 
\end{eqnarray*}
by {\bf Step 3$^\prime $.A.a2} with $\tilde{d}$. 
% $\vec{V}_j^{<\tilde{d}>}$  $(d=1,2,\ldots ,D)$ of $\Pi _N V$ at the intermediate steps with $\tilde{d}=0,1,\ldots , D_{\ell ^2}$ in {\bf Step 3$^\prime $.A} above 
Hence, $T_{\tilde{d}}=<\vec{V}_1^{<\tilde{d}>}, \vec{V}_2^{<\tilde{d}>},\ldots , \vec{V}_{\tilde{d}}^{<\tilde{d}>}>$ while \par\noindent $R_{\tilde{d}}=<\vec{V}_{\tilde{d}+1}^{<\tilde{d}>}, \vec{V}_{\tilde{d}+2}^{<\tilde{d}>},\ldots , \vec{V}_D^{<\tilde{d}>}>$.
In each substep, they are expressed as linear combinations  $\displaystyle \vec{V}_j^{<\tilde{d}>}=\sum_{d=0}^D a_{j,d}^{<\tilde{d}>} \, \vec{F}_{\rm int.}^{\, (d)}$ of the initial basis vectors  $\vec{F}_{\rm int.}^{\, (1)},\vec{F}_{\rm int.}^{\, (2)},\ldots ,\vec{F}_{\rm int.}^{\, (D)}$ (obtained by {\bf Step 2}) with the coefficients $a_{j,d}^{<\tilde{d}>}$ $(j,d=1,2,\ldots ,D;\, \tilde{d}=0,1,\ldots ,D_{\ell ^2})$.%, and we have only to explain how the coefficients $a_{j,d}^{<\tilde{d}>}$ $(j,d=1,2,\ldots ,D;\, \tilde{d}=0,1,\ldots ,D_{\ell ^2})$ are changed at each substep. 

At the initial step of {\bf Step 3$^\prime $}, with $\tilde{d}=0$, let $\vec{V}_d^{<0>}:=\vec{F}_{\rm int.}^{\, (d)}$ $(d=1,2,\ldots ,D)$. From this initial basis system of $\Pi _N V$, perform the following concrete procedures:

{\bf (*)} \, For {\bf Step 3$^\prime $.A.a1}, with $\vec{v}_d = \vec{V}_d^{<\tilde{d}>}$ $(d=1,2,\ldots ,D)$, do the iteration of {\bf R1} and then the iteration of {\bf R2} in Table \ref{tbl:4} with a sufficiently large integer $g$. (How to choose $g$ will be explained later.)\, These procedures are omitted for the exceptional case of $\tilde{d}=0$. In the result of these iterations, 
the vectors  $\vec{v}_j$ $(\tilde{d}+1\le j\le D)$ satisfy  
$\displaystyle\frac{|\langle \vec{v}_j, \vec{v}_\ell \rangle _{\ell ^2,K} |}{\|\vec{v}_j\|_{\ell ^2,K} \cdot\|\vec{v}_\ell \|_{\ell ^2,K}} \le \frac{1}{g} $ for $\, \tilde{d}+1 \le j<\ell \le D$ or $1\le j < \tilde{d}+1 \le \ell \le D$ when $\tilde{d}\ge 1$. 

Next, for {\bf Step 3$^\prime $.A.2}, with $\vec{v}_d$ $(d=1,2,\ldots ,D)$ 
%innovated 
renewed   
above, do the iteration of {\bf P$^\prime$1} and then the iteration of {\bf P$^\prime$2} in Table \ref{tbl:4}. 

Next,  for {\bf Step 3$^\prime $.A.3}, with $\displaystyle d_{\rm opt.}^{<\tilde{d}>}:=\mathop{\rm argmin}_{d\in\{\tilde{d}+1,\ldots ,D\}} \displaystyle \frac{\|\vec{v_d}\|_{{Q,N}}}{\|\vec{v_d}\|_{\ell ^2,K}}$ for the vectors $\vec{v}_d$ \par\noindent $(d=\tilde{d}+1,\ldots D)$ after these procedures, define $\vec{G}^{\, (\tilde{d}+1)}:=\vec{v}_{d_{\rm opt.}^{<\tilde{d}>}}$ $\left(=\vec{u}_{d_{\rm opt.}^{<\tilde{d}>}-\tilde{d}}^{<\tilde{d}>}\right)$. In the unified framework of change of the basis system mentioned above, this is equivalent to  $\vec{V}_{\tilde{d}+1}^{<\tilde{d}+1>}:=\vec{v}_{d_{\rm opt.}^{<\tilde{d}>}}$. This is an implicit process for {\bf Step 3$^\prime $.A.4}. Here note that the basis vector $\vec{V}_{\tilde{d}+1}^{<\tilde{d}+1>}$ is fixed exactly then and it remains fixed. At any steps in the procedures in Table \ref{tbl:4}, the vectors $\vec{V}_{d}^{<\tilde{d}>}$ with $d\le \tilde{d}$ are not changed. Moreover, define $\vec{V}_{d}^{<\tilde{d}+1>}:=\vec{V}_{d}^{<\tilde{d}>}$ for $d=1,2,\ldots \tilde{d}$.

With the above innovation of the basis system,   
the basis vectors $\vec{V}_j^{<\tilde{d}+1>}$ \par\noindent $(1\le j \le D)$ satisfy  $\displaystyle\frac{|\langle \vec{V}_j^{<\tilde{d}+1>}, \vec{V}_\ell ^{<\tilde{d}+1>}\rangle _{Q,N} |}{\|\vec{V}_j^{<\tilde{d}+1>}\|_{Q,N} \cdot\|\vec{V}_\ell ^{<\tilde{d}+1>}\|_{Q,N}} \le \frac{1}{h}
 $ for $\, \tilde{d}+1 \le j<\ell \le D$\, and they satisfy 
$\displaystyle\frac{|\langle \vec{V}_j^{<\tilde{d}+1>}, \vec{V}_\ell ^{<\tilde{d}+1>}\rangle _{\ell ^2,K} |}{\|\vec{V}_j^{<\tilde{d}+1>}\|_{\ell ^2,K} \cdot\|\vec{V}_\ell ^{<\tilde{d}+1>}\|_{\ell ^2,K}} \le \frac{1}{g} \,\sqrt{\frac{D-\tilde{d}}{1-\frac{\displaystyle D-\tilde{d}-1}{\displaystyle g}}} $ for  \par\noindent $1\le j < \tilde{d}+1 \le \ell \le D$ when $g\ge D-\tilde{d}$, where the latter inequality will be proved in Section \ref{sec:conv_multidim}. 

Since the vectors $\vec{V}_{d}^{<\tilde{d}+1>}$ with $d=1,2,\ldots \tilde{d}$ have been fixed to be $\vec{G}^{\, (d)}$, from the discussions above,  the following quasi-orthogonalities (a)-(e) are guaranteed  simultaneously when $\tilde{d}\ge 1$ and $g\ge D-\tilde{d}$. 
\begin{lemma}
\label{lemma:quasi_ort}
For the vector systems defined above with $\tilde{d}\ge 1$, the inequalities below are satisfied:
\begin{description}
\item[]$(a)$: The basis system of $T_{\tilde{d}}$ satisfies \par\noindent $\displaystyle\frac{|\langle \vec{G}^{\, (j)}, \vec{G}^{\, (\ell )}\rangle _{\ell ^2,K} |}{\|\vec{G}^{\, (j)}\|_{\ell ^2,K} \cdot\|\vec{G}^{\, (\ell )}\|_{\ell ^2,K}} \le \frac{1}{g} \,\sqrt{\frac{D-\tilde{d}}{1-\frac{\displaystyle D-\tilde{d}-1}{\displaystyle g}}} $
for $1\le j < \ell \le \tilde{d}$. 
\item[]$(b)$:
The old basis system of $R_{\tilde{d}}$ satisfies 
 $\displaystyle\frac{|\langle \vec{v}_j^{<\tilde{d}>}, \vec{v}_\ell ^{<\tilde{d}>}\rangle _{\ell ^2,K} |}{\|\vec{v}_j^{<\tilde{d}>}\|_{\ell ^2,K} \cdot\|\vec{v}_\ell ^{<\tilde{d}>}\|_{\ell ^2,K}} \le \frac{1}{g} $ for  \par\noindent $1\le j < \ell \le D-\tilde{d}$.
\item[]$(c)$: The new basis system of $R_{\tilde{d}}$ satisfies 
$\displaystyle\frac{|\langle \vec{u}_j^{<\tilde{d}>}, \vec{u}_\ell ^{<\tilde{d}>}\rangle _{Q,N} |}{\|\vec{u}_j^{<\tilde{d}>}\|_{Q,N} \cdot\|\vec{u}_\ell ^{<\tilde{d}>}\|_{Q,N}} \le \frac{1}{h}$ for \par\noindent $1\le j < \ell \le D-\tilde{d}$. 
\item[]$(d)$: Between the basis system of $T_{\tilde{d}}$  and the old basis system of $R_{\tilde{d}}$,  the inequality \par\noindent
$\displaystyle\frac{|\langle \vec{G}^{\, (j)}, \vec{v}_\ell ^{<\tilde{d}>}\rangle _{\ell ^2,K} |}{\|\vec{G}^{\, (j)}\|_{\ell ^2,K} \cdot\|\vec{v}_\ell ^{<\tilde{d}>}\|_{\ell ^2,K}} \le \frac{1}{g} $ holds  for $1\le j \le  \tilde{d}$ and  $1\le \ell \le N-\tilde{d}$. 
\item[]$(e)$: Between the basis system of $T_{\tilde{d}}$  and the new basis system of $R_{\tilde{d}}$,  the inequality \par\noindent
$\displaystyle\frac{|\langle \vec{G}^{\, (j)}, \vec{u}_\ell ^{<\tilde{d}>}\rangle _{\ell ^2,K} |}{\|\vec{G}^{\, (j)}\|_{\ell ^2,K} \cdot\|\vec{u}_\ell ^{<\tilde{d}>}\|_{\ell ^2,K}} \le \frac{1}{g} \,\sqrt{\frac{D-\tilde{d}}{1-\frac{\displaystyle D-\tilde{d}-1}{\displaystyle g}}} $ holds  for $1\le j \le  \tilde{d}$ and \par\noindent $1\le \ell \le N-\tilde{d}$. 
\end{description}
\end{lemma} \par
\noindent 
The proof of (a) and (e) is given in the first part of Section \ref{sec:conv_multidim}, whereas (b)-(d) is obvious by Table \ref{tbl:4}.  
Here, note that the quasi-orthogonality (c) is with respect to $\langle \cdot , \cdot \rangle _{Q,N}$ while the others are with respect to $\langle \cdot , \cdot \rangle _{\ell ^2,K}$. 
With $h\ge D-\tilde{d}$, the quasi-orthogonality (c) guarantees that $\vec{G}^{\, (\tilde{d}+1)}=\vec{u}_{d_{\rm opt.}^{<\tilde{d}>}-\tilde{d}}^{<\tilde{d}>}$ belongs to $((\sigma_{K,N}^{(\Omega )})^{-1}[0,c \underline{\sigma_{K,N}^{(\Omega )}} ])$  with $\displaystyle c={\frac{D-\tilde{d}}{1- \frac{D-\tilde{d}-1}{h}}}$, from a similar discussion to Section \ref{sec:suboptimality}, with $D-\tilde{d}$ instead of $D$. Hence, the convergence to $\Pi _K(V\cap \ell ^2(\mathbb{Z}^+))$ is similarly guaranteed by Theorem \ref{thm:opt_ratio_gen}. 
The quasi-orthogonality (a) is essential later in order to prove that any nonzero vector in $T_{D_{\ell ^2}}$ (i.e., any nonzero linear combination of $\vec{G}^{\, (1)}, \vec{G}^{\, (2)}, \ldots \vec{G}^{\, (D_{\ell ^2})}$) belongs to $((\sigma_{K,N}^{(\Omega )})^{-1}[0,c \underline{\sigma_{K,N}^{(\Omega )}} ])$ with a fixed finite $c$, as was mentioned briefly at the beginning of this section.

Next, as a preparation for the next step, we choose the basis vectors  $\vec{V}_{\tilde{d}+2}^{<\tilde{d}+1>},\vec{V}_{\tilde{d}+3}^{<\tilde{d}+1>},\ldots , \vec{V}_{D}^{<\tilde{d}+1>}$ so that $\vec{V}_{1}^{<\tilde{d}+1>},\vec{V}_{2}^{<\tilde{d}+1>},\ldots , \vec{V}_{D}^{<\tilde{d}+1>}$ may be a basis system for $\Pi _N V$.  For this, linear independence suffices, because $\dim \Pi _N V=D$. 
An easy method is the choice of the vectors other than $\vec{u}_{d_{\rm opt.}^{<\tilde{d}>}-\tilde{d}}^{<\tilde{d}>}$ in $\vec{u}_1^{<\tilde{d}>}, \vec{u}_2^{<\tilde{d}>}, \ldots \vec{u}_{D-\tilde{d}}^{<\tilde{d}>}$ (i.e., the vectors other than $\vec{v}_{d_{\rm opt.}^{<\tilde{d}>}}$ int $\vec{v}_{\tilde{d}+1}, \vec{v}_{\tilde{d}+2}, \ldots \vec{v}_{D}$ after {\bf Step 3$^\prime$.A.3} )  for the vectors  $\vec{V}_{\tilde{d}+2}^{<\tilde{d}+1>},\vec{V}_{\tilde{d}+3}^{<\tilde{d}+1>},\ldots , \vec{V}_{D}^{<\tilde{d}+1>}$. 

However, this method requires more calculations than the choice from the initial basis vectors $\vec{F}_{\rm int.}^{\, (1)},\vec{F}_{\rm int.}^{\, (2)},\ldots ,\vec{F}_{\rm int.}^{\, (D)}$, because the quasi-orthogonalization procedures for {\bf Step 3$^\prime $.2} %have augmented the size of the integers in the numerators and the denominators of the vectors. Moreover, the quasi-orthogonalization
 with respect to $\langle \cdot , \cdot \rangle _{Q,N}$ have made the basis vectors almost parallel to one another with respect to $\langle \cdot , \cdot \rangle _{\ell ^2,K}$. Hence, a choice from $\vec{F}_{\rm int.}^{\, (1)},\vec{F}_{\rm int.}^{\, (2)},\ldots ,\vec{F}_{\rm int.}^{\, (D)}$ is desirable.  The check of linear independence can be made then with only the coefficients  $a_{j,d}^{<\tilde{d}+1>}$ $(j,d=1,2,\ldots ,D)$, as follows; with the definitions of vectors $\vec{a}_d$ by  $(\vec{a}_d)_j:=\tilde{a}_{j,d}$ and the  set $\tilde{S}$ of the $D-\tilde{d}-1$ numbers of the chosen vectors from $\vec{F}_{\rm int.}^{\, (1)},\vec{F}_{\rm int.}^{\, (2)},\ldots ,\vec{F}_{\rm int.}^{\, (D)}$  such that the statement $n\in\tilde{S}$ is equivalent to the statement that the vector $\vec{F}_{\rm int.}^{\, (n)}$ is chosen as one of  the basis vectors,  where the linear independence is guaranteed if the vector system  $\bigl\{ \vec{a}_d\, \bigl|\, d\in \{1,2,\ldots ,D\}\backslash \tilde{S} \bigr\}$ is linearly independent.
Since there exists at least one choice of $\tilde{S}$ such that this system is linearly independent, 
 which is easily shown from the linear independence of $\vec{V}_{1}^{<\tilde{d}+1>},\vec{V}_{2}^{<\tilde{d}+1>},\ldots , \vec{V}_{\tilde{d}+1}^{<\tilde{d}+1>}$, this type of choice always exists.

In other words, there exist integers $n_{\tilde{d}+2},n_{\tilde{d}+3},\ldots ,n_D$ such that 
the vector system $\vec{V}_{1}^{<\tilde{d}+1>},\vec{V}_{2}^{<\tilde{d}+1>},\ldots , \vec{V}_{\tilde{d}+1}^{<\tilde{d}+1>}, \vec{F}_{\rm int.}^{\, (n_{\tilde{d}+2})}, \vec{F}_{\rm int.}^{\, (n_{\tilde{d}+3})}, \ldots ,\vec{F}_{\rm int.}^{\, (n_D)}$ is a basis system of $\Pi _N V$. The set of these  integers should satisfy only the following conditions (i) and (ii):
\begin{description}
\item[]
(i): The vector system  $\bigl\{ \vec{a}_d\, \bigl|\, d\in \{1,2,\ldots ,D\}\backslash \{n_{\tilde{d}+2},n_{\tilde{d}+3},\ldots ,n_D\} \bigr\}$ is linearly independent. 
\item[]
(ii): $n_m\ne n_\ell $ if $m\ne \ell $.
\end{description} 
By the check of coefficients $a_{j,d}^{<\tilde{d}+1>}$, these integers can be easily chosen.  Empirically, except for very special cases with simple symmetry, most choices of $D-\tilde{d}-1$ distinct numbers from $\{1,2,\ldots ,D\}$ give the linear independence of (i). We have only to check the linear independence of (i) with an arbitrary choice of $\{n_{\tilde{d}+2},n_{\tilde{d}+3},\ldots ,n_D\}$ satisfying (ii), by the coefficients $a_{j,d}^{<\tilde{d}+1>}$. If, exceptionally, linear independence is not satisfied, replace one of these numbers by another, and try it again. How to determine the coefficients $a_{j,d}^{<\tilde{d}+1>}$ easily will be explained later, in the explanation of practical operations with a reduced amount of calculations. 

With these integers, define $\vec{V}_{d}^{<\tilde{d}+1>}:=\vec{F}_{\rm int.}^{\, (n_{d})}$ for $d=\tilde{d}+2, \, \tilde{d}+3, \ldots , D$.  Then, with the increment $\tilde{d}\to \tilde{d}+1$, return to the top {\bf (*)} of the procedures for {\bf Step 3$^\prime $.A.a1} if $\tilde{d}\le D_{\ell ^2}$. 

By means of the iterations explained above, we can obtain a quasi-orthogonal vector system  (with respect to $\langle \cdot ,\,\cdot \rangle _{\ell ^2,K}$) which satisfies the following theorem:
\begin{theorem}
\label{thm:app_multidim}
When $h\ge D_{\ell ^2}$ and \par\noindent $\displaystyle g>\frac{1}{2}\, \Bigl (\sqrt{(D_{\ell ^2}-1)^2+4D_{\ell ^2}(D_{\ell ^2}-1)}+(D_{\ell ^2}-1)\Bigr)$, 
any nonzero linear combination of $\,\, \vec{G}^{\, (1)},\vec{G}^{\, (2)},\ldots \vec{G}^{\, (D_{\ell ^2})}$ obtained by the above iterations $(${\bf Step 3$^\prime$}$)$ belongs to $\displaystyle ((\sigma_{K,N}^{(\Omega )})^{-1}[0,c \underline{\sigma_{K,N}^{(\Omega )}} ])$ with 
\begin{eqnarray*}
c(g,h,D_{\ell ^2}):={\frac{\displaystyle D_{\ell ^2}\, \sqrt{\frac{D_{\ell ^2}}{1- \frac{D_{\ell ^2}-1}{h}}}}{\displaystyle \,  1-\frac{D_{\ell ^2}-1}{g}\sqrt{\frac{D_{\ell ^2}}{1- \frac{D_{\ell ^2}-1}{g}}} \, }}\, . 
\end{eqnarray*}
\end{theorem}

\par\noindent The proof is given in the last part of Section \ref{sec:conv_multidim}. This theorem implies that the vector system $\,\, \Pi _K \vec{G}^{\, (1)},\Pi _K \vec{G}^{\, (2)},\ldots \Pi _K \vec{G}^{\, (D_{\ell ^2})}$ is approximately a quasi-orthogonal basis system for $\Pi _K (V \cap \ell ^2(\mathbb{Z}^+))$, because of Theorem \ref{thm:opt_ratio_gen}.

As was mentioned in Section \ref{sec:cp} for one-dimensional cases, the number of calculations can be remarkably reduced by some modifications. With these modifications, we propose a practical realization in Table \ref{tbl:5} instead of the procedures in Table \ref{tbl:4}. In this practical realization, the innovation of the coefficients $a_{j,d}^{<\tilde{d}+1>}$ is made at the last step.

\section{Suboptimality  under the extension to multidimensional cases}
\label{sec:conv_multidim} 
In this section, proving Lemma \ref{lemma:quasi_ort} and Theorem \ref{thm:app_multidim}, 
we show the validity of the procedures in extension {\bf Step 3$^\prime $} which has been proposed in Section \ref{sec:multidim} as an extension of {\bf Step 3} to the cases where 
$\dim V\cap \ell ^2(\mathbb{Z}^+ )\ge 2$.

First, as a preliminary process, we will show the quasi-orthogonality 
\begin{eqnarray*}
\displaystyle\frac{|\langle \vec{V}_j^{<\tilde{d}+1>}, \vec{V}_\ell ^{<\tilde{d}+1>}\rangle _{\ell ^2,K} |}{\|\vec{V}_j^{<\tilde{d}+1>}\|_{\ell ^2,K} \cdot\|\vec{V}_\ell ^{<\tilde{d}+1>}\|_{\ell ^2,K}} \le \frac{1}{g} \,\sqrt{\frac{D-\tilde{d}}{1-\frac{\displaystyle D-\tilde{d}-1}{\displaystyle g}}} \,\,\,\,\,\, (1\le j < \tilde{d}+1 \le \ell \le D)
\end{eqnarray*}
with respect to $\langle \cdot , \cdot \rangle _{\ell ^2,K}$ 
%mentioned in Section \ref{sec:multidim}. 
in (a) and (e) of Lemma \ref{lemma:quasi_ort}.  
Since  the vectors $\vec{V}_\ell ^{<\tilde{d}+1>}$ with \par\noindent $\tilde{d}+1\le \ell \le D$ are  linear combinations of the vectors $\vec{v}_{1}^{<\tilde{d}>},\vec{v}_{2}^{<\tilde{d}>},\ldots ,\vec{v}_{N-\tilde{d}}^{<\tilde{d}>}$ (i.e., of the vectors $\vec{v}_{\tilde{d}+1}, \vec{v}_{\tilde{d}+2}, \ldots \vec{v}_D$ after the iterations of {\bf S1} and {\bf S2} in Table \ref{tbl:4}) where the quasi-orthogonality 
$\displaystyle\frac{|\langle \vec{v}_j^{<\tilde{d}>}, \vec{v}_\ell ^{<\tilde{d}>}\rangle _{\ell ^2,K} |}{\|\vec{v}_j^{<\tilde{d}>}\|_{\ell ^2,K} \cdot\|\vec{v}_\ell ^{<\tilde{d}>}\|_{\ell ^2,K}} \le \frac{1}{g}$  $(1\le j<\ell \le D-\tilde{d})$ %with respect to $\langle \cdot , \cdot \rangle _{\ell ^2,K}$ 
is guaranteed, the above quasi-orthogonality for $\vec{V}_j^{<\tilde{d}+1>}$
, i.e., (a) and (e) of Lemma \ref{lemma:quasi_ort}  
can be shown directly from the following lemma, with $n=D-\tilde{d}$,  $\vec{v}_m= \vec{u}_m^{<\tilde{d}>}$ $(m=1,2,\ldots ,D-\tilde{d}\, )$ and $\displaystyle \zeta=\epsilon = \frac{1}{g}$\, :
\begin{lemma}
\label{lemma:min_dir_cos_2} 
When any pair of two distinct vectors in the set $\{ \vec{v}_1, \vec{v}_2, \ldots \vec{v}_n \}$  satisfies
$\displaystyle \frac{|( \vec{v}_m,\, \vec{v}_\ell ) |}{\| \vec{v}_m \| \cdot \| \vec{v}_\ell \|} \le  \zeta \, $ $\bigl($with fixed $\zeta$ such that $\displaystyle 0<\zeta<\frac{1}{n-1}$ $\bigr)$ and there is a vector $\vec{u}$ such that  the inequality  $\displaystyle \frac{|( \vec{v}_m,\, \vec{u} ) |}{\| \vec{v}_m \| \cdot \| \vec{u}\|} \le  \epsilon \, $ holds  for $m=1,2,\ldots ,n$, 
then the inequality  
$\, \displaystyle  \frac{|(\vec{w},\, \vec{u})|}{\|\vec{w}\|\cdot \|\vec{u}\|} \le \epsilon \, \, \sqrt{\frac{n\,  }{1- (n-1)\, \zeta  }} $ holds for any vector $\, \vec{w} \in <\vec{v}_1, \vec{v}_2, \ldots \vec{v}_n>\backslash \{ 0\}\,$.
\end{lemma} \par\noindent 
\par\noindent\noindent{\em Proof of Lemma }\ref{lemma:min_dir_cos_2}: \quad
 \rm
Define $\displaystyle \, \vec{p}_j := \frac{\vec{v}_j}{\|\vec{v}_j\|} \,\,\, (j=1,2,...,n) \, $ and $\displaystyle \vec{q}:=\frac{\vec{u}}{\|\vec{u}\|}.$ Then, from the condition of the lemma, for any pair of distinct vectors in the set  $\, \{ \vec{p}_j \,|\, j=1,2,...,n\}\, $,  the inequality $\, |( \vec{p}_m,\, \vec{p}_\ell ) | \le  \zeta \, $  holds. Moreover, the inequality $\, |( \vec{p}_m,\, \vec{q}) | \le  \epsilon  \, $ holds  for $m=1,2,\ldots ,n$.  Since $\vec{w} \in <\vec{v}_1, \vec{v}_2, \ldots \vec{v}_n>\backslash \{ 0\}$,
there is a set of complex coefficients $\, b_j \,\,\, (j=1,2,..,n) \, $ such that $\displaystyle\, \vec{w}\, = \sum_{\ell =1}^n b_\ell \vec{p}_\ell \, $  and 
$\displaystyle\, A:=\sum_{\ell =1}^n|b_\ell |^2  >0 $.  From the Schwarz inequality and the condition of this lemma, the following two inequalities can be derived:
\begin{eqnarray*}
|(\vec{w}, \vec{u})|^2 \, =\,  \bigl| \sum_{\ell =1}^n \overline{b_\ell} \, (\vec{p}_\ell , \vec{q} ) \bigr|^2 
\, \le \, \bigl( \sum_{\ell =1}^n |b_\ell |^2 \bigr)\, \cdot \, \bigr( \sum_{\ell =1}^n \epsilon ^2 \bigr) \, = \, n\, \epsilon ^2 A
\end{eqnarray*}
\begin{eqnarray*} \|\vec{w}\|^2 \, &=& \, \sum_{\ell =1}^n \sum_{\ell '=1}^n b_\ell \overline{b_{\ell '}} (\vec{p}_\ell , \vec{p}_{\ell '} )  
\, \ge \, \, -\,  \bigl( \sum_{\ell =1}^n |b_\ell | \bigr)^2 \, \zeta + \sum_{\ell =1}^n |b_\ell |^2 \, (1+\zeta) \\ 
\, &\ge &\, -\, \bigl(\sum_{\ell =1}^n |b_\ell |^2 \bigr) \, \bigl(\sum_{\ell =1}^n 1 \bigr) \, \zeta + (1+\zeta)A 
\, = \, \left(1 \, - \, (n-1) \, \zeta \right) A
\end{eqnarray*}
where we have utilized the relation $\, (\vec{p}_\ell , \vec{p}_\ell)=1=(1+\zeta)-\zeta\, $. Since $A>0$, 
\begin{eqnarray*}
\frac{|(\vec{v}_m, \vec{u})|^2}{||\vec{v}_m||^2 \, ||\vec{u}||^2} 
\, = \, \frac{|(\vec{p}_m, \vec{w})|^2}{||\vec{w}||^2} \le \frac{n\, \epsilon ^2}{1- (n-1) \, \zeta}  . 
\end{eqnarray*}
\hfill\endproof
\par\par\noindent 

In the following, by means of this quasi-orthogonality and Lemma \ref{lemma:max_ratio_pl}, %we show 
we prove Theorem \ref{thm:app_multidim},  
%to show 
i.e.,  
the convergence of {\it any nonzero linear combination} of the vectors \par\noindent $\Pi _K \vec{G}^{\, (1)},\Pi _K \vec{G}^{\, (2)},\ldots \Pi _K \vec{G}^{\, (D_{\ell ^2})}$ to 
%$\Pi _N (V\cap\ell ^2(\mathbb{Z}^+))$, 
$\Pi _K (V \cap \ell ^2(\mathbb{Z}^+))$,   
under the choice of the integers $h$ and $g$ such that  $h\ge D_{\ell ^2}$ and $\displaystyle g>\frac{1}{2}\, \Bigl (\sqrt{(D_{\ell ^2}-1)^2+4D_{\ell ^2}(D_{\ell ^2}-1)}+(D_{\ell ^2}-1)\Bigr)$
%, with which the inequalities $\displaystyle1-\frac{D_{\ell ^2}-1}{h}>0$, $\displaystyle 1-\frac{D_{\ell ^2}-1}{g}>0$ and $\displaystyle 1-\frac{D_{\ell ^2}-1}{g}\sqrt{\frac{D_{\ell ^2}}{1- \frac{D_{\ell ^2}-1}{g}}}>0$ are guaranteed.
. 
\par\noindent\noindent{\em Proof of Theorem }\ref{thm:app_multidim}: \quad
\rm 
From the conditions for $h$ and $g$, the inequalities \par\noindent $\displaystyle1-\frac{D_{\ell ^2}-1}{h}>0$, $\displaystyle 1-\frac{D_{\ell ^2}-1}{g}>0$ and $\displaystyle 1-\frac{D_{\ell ^2}-1}{g}\sqrt{\frac{D_{\ell ^2}}{1- \frac{D_{\ell ^2}-1}{g}}}>0$ hold.  
Let $W$ be the span of the vectors $\vec{G}^{\, (1)},\vec{G}^{\, (2)},\ldots \vec{G}^{\, (D_{\ell ^2})}$. This is a subspace of $\Pi _N V$. As has been explained in Section \ref{sec:multidim}, the vector $\vec{G}^{\, (\tilde{d})}\ $ belongs to $((\sigma_{K,N}^{(\Omega )})^{-1}[0,c \underline{\sigma_{K,N}^{(\Omega )}} ])$  with $\displaystyle c={\frac{D_{\ell ^2}-\tilde{d}}{1- \frac{D_{\ell ^2}-\tilde{d}-1}{h}}}$, 
and hence it belongs to $\displaystyle ((\sigma_{K,N}^{(\Omega )})^{-1}[0, \frac{D_{\ell ^2}}{1- \frac{D_{\ell ^2}-1}{h}}\cdot \underline{\sigma_{K,N}^{(\Omega )}} ])$ 
because \par\noindent $\displaystyle {\frac{D_{\ell ^2}-\tilde{d}}{1- \frac{D_{\ell ^2}-\tilde{d}-1}{h}}}\le {\frac{D_{\ell ^2}}{1- \frac{D_{\ell ^2}-1}{h}}}$. Moreover, as has been shown in Lemma \ref{lemma:quasi_ort} with the above proof,  these vectors satisfy the quasi-orthogonality 
\begin{eqnarray*}
\displaystyle\frac{|\langle \vec{G}^{\, (j)}, \vec{G}^{\, (\ell )}\rangle _{\ell ^2,K} |}{\|\vec{G}^{\, (j)}\|_{\ell ^2,K} \cdot\|\vec{G}^{\, (\ell )}\|_{\ell ^2,K}} \le \!\frac{1}{g} \sqrt{\frac{D_{\ell ^2}-\tilde{d}}{1-\frac{\displaystyle D_{\ell ^2}-\tilde{d}-1}{\displaystyle g}}}\le \!\frac{1}{g} \sqrt{\frac{D_{\ell ^2}}{1- \frac{D_{\ell ^2}-1}{g}}} \,\, (1\le j < \ell \le D_{\ell ^2}).
\end{eqnarray*}
Hence, from Lemma \ref{lemma:max_ratio_pl} with $\langle \cdot ,\cdot \rangle _\Lambda = \langle \cdot ,\cdot \rangle _{\ell ^2, K}$, $\langle \cdot ,\cdot \rangle _\Xi = \langle \cdot ,\cdot \rangle _{Q,N}$ and \par\noindent $\displaystyle \xi = \frac{1}{g} \,\sqrt{\frac{D_{\ell ^2}}{1- \frac{D_{\ell ^2}-1}{g}}}$, we can show the inequalities 
\begin{eqnarray*}
\frac{\|\vec{v}\|_{Q,N}^2}{\|\vec{v}\|_{\ell ^2,K}^2}\le {\frac{\displaystyle\sum_{m=1}^{D_{\ell ^2}}
%\frac{\|\vec{v}_m\|_{Q,N} ^2}{\|\vec{v}_m\|_{\ell ^2,K} ^2}}
\frac{\|\vec{G}^{\, (m)}\|_{Q,N} ^2}{\|\vec{G}^{\, (m)}\|_{\ell ^2,K} ^2}} 
{\displaystyle \,  1-\frac{D_{\ell ^2}-1}{g}\sqrt{\frac{D_{\ell ^2}}{1- \frac{D_{\ell ^2}-1}{g}}} \, }}\le {\frac{\displaystyle D_{\ell ^2}\, \sqrt{\frac{D_{\ell ^2}}{1- \frac{D_{\ell ^2}-1}{h}}}}{\displaystyle \,  1-\frac{D_{\ell ^2}-1}{g}\sqrt{\frac{D_{\ell ^2}}{1- \frac{D_{\ell ^2}-1}{g}}} \, }}
\end{eqnarray*}
for any vector $\vec{v}$ in $< \vec{G}^{(1)}, \vec{G}^{(2)},\ldots , \vec{G}^{(D_{\ell ^2})}>  \backslash \{ 0\}$.
This implies that any nonzero linear combination of $\vec{G}^{\, (1)},\vec{G}^{\, (2)},\ldots \vec{G}^{\, (D_{\ell ^2})}$ belongs to $\displaystyle ((\sigma_{K,N}^{(\Omega )})^{-1}[0,c(g,h,D_{\ell ^2}) \underline{\sigma_{K,N}^{(\Omega )}} ])$  with 
\begin{eqnarray*}
c(g,h,D_{\ell ^2}):={\frac{\displaystyle D_{\ell ^2}\, \sqrt{\frac{D_{\ell ^2}}{1- \frac{D_{\ell ^2}-1}{h}}}}{\displaystyle \,  1-\frac{D_{\ell ^2}-1}{g}\sqrt{\frac{D_{\ell ^2}}{1- \frac{D_{\ell ^2}-1}{g}}} \, }}.
\end{eqnarray*}
\hfill\endproof 
\par\noindent 
This fact guarantees the convergence of {\it any nonzero linear combination} of the vectors $\Pi _K \vec{G}^{\, (1)},\Pi _K \vec{G}^{\, (2)},\ldots \Pi _K \vec{G}^{\, (D_{\ell ^2})}$ to %$\Pi _N (V\cap\ell ^2(\mathbb{Z}^+))$ 
$\Pi _K (V\cap\ell ^2(\mathbb{Z}^+))$  
in the same sense as Section \ref{sec:suboptimality}, under the choice of integers $h$ and $g$ such that $h\ge D_{\ell ^2}$ and  \par\noindent $\displaystyle g>\frac{1}{2}\, \Bigl (\sqrt{(D_{\ell ^2}-1)^2+4D_{\ell ^2}(D_{\ell ^2}-1)}+(D_{\ell ^2}-1)\Bigr)$.

\section{Proof of the halting of {\bf Step 3.1}}
\label{sec:halt}
Since the procedures in Table \ref{tbl:2} and Table \ref{tbl:3} for the process {Step 3.1} contain iterations which finish only under certain conditions, we should verify that they 
%terminate. 
halt.  
Otherwise, the algorithm could not be executed.  
In this section, we prove that they 
%terminate.
halt in a finite number of steps.  

The main idea used in the proof of 
%termination 
halting  
is based on the finiteness of the number of vectors in a lattice with bounded norms and the monotonic decrease of the norm, except for finitely many times, in the execution of $\vec{v}_j\leftarrow 2\vec{v}_j$ in {\bf Q2} in Table \ref{tbl:2}.

%In the following, we will use the extended definition of the angle $\theta $ between {\it complex} vectors $\vec{v}$ and $\vec{u}$ by  $\displaystyle\cos \theta = \frac{|(\vec{v},\vec{u})|}{\|\vec{v}\|\cdot\|\vec{u}\|}$ and $0\le \theta \le \frac{\pi}{2}$. 
It is easily shown that the vectors $\vec{v}_1, \vec{v}_1,\ldots ,\vec{v}_D$ satisfy 
\begin{eqnarray} 
|{\rm Re} \, ( \vec{v}_m,\, \vec{v}_\ell )_u | \le \frac{1}{2} \, \| \vec{v}_\ell \|_u^2 \,\, , 
\,\,\,\,\, 
|{\rm Im} \, ( \vec{v}_m,\, \vec{v}_\ell )_u | \le \frac{1}{2} \, \| \vec{v}_\ell \|_u^2 
\\    {\rm and \,\, hence } \,\,\,\, 
|( \vec{v}_m,\, \vec{v}_\ell )_u |^2 \le \frac{1}{2} \, \| \vec{v}_\ell \|_u^4 \,\, , 
\,\,\,\,\,\,\,\,\,\,\, {\rm if} \,\,\,\,\, m\ne\ell \,\, . \nonumber 
\label{eqn:qsot}\end{eqnarray}
after the iteration of {\bf Q1} of Table \ref{tbl:2}. 

This inequality yields the following theorem:
\begin{theorem}
\label{thm:halt}
The procedures in Table \ref{tbl:2}  
%terminate 
halt  
within a finite number of steps.
\end{theorem} \par\noindent 
For the proof of this, we begin with some preliminary definitions and a lemma:
\begin{definition}
\label{definition:lattice}
For a set of $n$ vectors $\vec{u}_1,\vec{u}_2,\ldots ,\vec{u}_n$, define %the lattice generated by them 
\begin{eqnarray*}
{\rm Latt}(\vec{u}_1,\vec{u}_2,\ldots ,\vec{u}_n):=\displaystyle\bigl\{\sum_{\ell =1}^n z_\ell  \vec{u}_\ell  \, \bigl| \, z_\ell \in \mathbb{Z} +\mathbb{Z} i \,\, (m=1,\dots n)\bigr\}\, . 
\end{eqnarray*}
Obviously, ${\rm Latt}(\vec{u}_1,\vec{u}_2,\ldots ,\vec{u}_n)\subset <\vec{u}_1,\vec{u}_2,\ldots ,\vec{u}_n>$
\end{definition}
\begin{definition}
\label{definition:lattice_finite}
For a set of $n$ vectors $\vec{u}_1,\vec{u}_2,\ldots ,\vec{u}_n$ in a linear space $U$ with norm $\|\cdot \|_\Xi $, define 
$\displaystyle s_\Xi(\vec{u}_1,\ldots ,\vec{u}_n):=\sum_{j=1}^D \|\vec{u}_j\|_\Xi$
and 
\begin{eqnarray*}
\tilde{T}_\Xi^{(s)}(\vec{u}_1,\ldots ,\vec{u}_n):=\bigl\{ \, (\vec{f}_1,\ldots ,\vec{f}_D) \in \bigl({\rm Latt}(\vec{u}_1,\ldots ,\vec{u}_n)\bigr)^n \bigl|\,s_\Xi(\vec{f}_1,\ldots ,\vec{f}_n) \le s \bigr\} \, .
\end{eqnarray*}
\end{definition}
\begin{lemma}
\label{lemma:max_distance_lattice}
Let $U$ be a linear space and $\langle \cdot , \cdot \rangle$ be an inner product there.
Then, for a set of $n$ vectors $\vec{u}_1,\vec{u}_2,\ldots ,\vec{u}_n$ in $U$, with norm $\|\vec{f}\|:=\sqrt{\langle\vec{f},\vec{f}\rangle}$,
\begin{eqnarray*}
\sup_{\vec{f}\in <\vec{u}_1,\ldots ,\vec{u}_n>}\inf_{\vec{g}\in {\rm Latt}(\vec{u}_1,\ldots ,\vec{u}_n)}\|\vec{f}-\vec{g}\|\le \sqrt{\frac{n}{2}}\max _{r\in\{1,\ldots n\}}\|\vec{u}_r\|\, .
\end{eqnarray*}
\end{lemma}
\par\noindent\noindent{\em Proof of Lemma }\ref{lemma:max_distance_lattice}: \quad
\rm
This lemma will be proved by mathematical induction. 
When $n=1$, the statement in the lemma holds because 
\begin{eqnarray*}
\sup_{\vec{f}\in <\vec{u}_1>}\inf_{\vec{g}\in {\rm Latt}(\vec{u}_1)}\|\vec{f}-\vec{g}\|=\left(\max_{-\frac{1}{2}\le x \le \frac{1}{2},\, -\frac{1}{2}\le y \le \frac{1}{2}} |x+iy| \right)\|\vec{u}_1||=\sqrt{\frac{1}{2}}\,\, \|\vec{u}_1|| \, .
\end{eqnarray*}
Let $P_{<\vec{u}_1,\ldots ,\vec{u}_r>}$ be the orthogonal projector to $<\vec{u}_1,\ldots ,\vec{u}_r>$ with respect to $\langle \cdot , \cdot \rangle$. %Then, $ ^\forall \vec{f}\in U,\,\, $ $P_{{\rm Span}(\vec{u}_1,\ldots ,\vec{u}_r)}\vec{f}\in {\rm Span}(\vec{u}_1,\ldots ,\vec{u}_r) $ $\|P_{{\rm Span}(\vec{u}_1,\ldots ,\vec{u}_r)}\vec{f}\|\le \|\vec{f}\|$.
If the statement in the lemma holds with $n=n'$, then 
\begin{eqnarray*}
&& ^\forall \vec{f}\in <\vec{u}_1,\ldots ,\vec{u}_{n'}>, \,\,\,  -\frac{1}{2}\le ^\forall x \le \frac{1}{2} \,\, \mbox{ and } \,\, -\frac{1}{2}\le ^\forall y \le \frac{1}{2}, \\ \\&& \displaystyle
\sup_{\vec{f}\in <\vec{u}_1,\ldots ,\vec{u}_{n'}>}\inf_{\vec{g}\in {\rm Latt}(\vec{u}_1,\ldots ,\vec{u}_{n'})}\|\vec{f}+(x+iy)\vec{u}_{n'+1}-\vec{g}\| \\ 
&&\le \left(\sqrt{\frac{n'}{2}}\max _{r\in\{1,\ldots ,n'\}}\|\vec{u}_r\|\right)^2+\|P_{<\vec{u}_1,\ldots ,\vec{u}_r>}(x+iy)\vec{u}_{n'+1}\|^2 \\  
&&\le \frac{n'}{2}\max _{r\in\{1,\ldots ,n'\}}\|\vec{u}_r\|^2+\frac{1}{2}\, \|\vec{u}_{n'+1}\|^2 \, 
\le \,\frac{n'+1}{2}\max _{r\in\{1,\ldots ,n'+1\}}\|\vec{u}_r\|^2 .
\end{eqnarray*} 
This implies that the statement in the lemma holds also for $n=n'+1$, 
because any vector $\vec{h}$ in $<\vec{u}_1,\ldots ,\vec{u}_{n'+1}>$ can be decomposed as $\vec{h}=\vec{f}+(x+iy)\vec{u}_{n'+1}+z\vec{u}_{n'+1}$ with $\vec{f}\in <\vec{u}_1,\ldots ,\vec{u}_{n'}>$, $-\frac{1}{2}\le  x \le \frac{1}{2}$, $-\frac{1}{2}\le y \le \frac{1}{2}$ and $z\in\mathbb{Z} + \mathbb{Z} i$.
\hfill\endproof
\par \par\noindent 
\par\noindent\noindent{\em Proof of Theorem }\ref{thm:halt}: \quad
\rm
Let $\widetilde{\vec{V}}_1,\widetilde{\vec{V}}_2,\ldots ,\widetilde{\vec{V}}_D$ be the vectors $\vec{v}_1,\vec{v}_2,\ldots ,\vec{v}_D$ after the iteration of {\bf Q1} and before the iteration of {\bf Q2} in Table \ref{tbl:2}, and let $\vec{V}_1,\vec{V}_2,\ldots ,\vec{V}_D$ be the vectors $\vec{v}_1,\vec{v}_2,\ldots ,\vec{v}_D$ after all the procedures in Table \ref{tbl:2}. 
Since $\vec{F}_{\rm int.}^{(1)},\vec{F}_{\rm int.}^{(2)},\ldots ,\vec{F}_{\rm int.}^{(D)}$ are linearly independent, so also are $\widetilde{\vec{V}}_1,\widetilde{\vec{V}}_2,\ldots ,\widetilde{\vec{V}}_D$, because {\bf Q1} does not change linear independence. 
Hence, for $2\le j\le D$,\, with the orthogonal projector \par\noindent $P_{<\widetilde{\vec{V}}_1,\ldots ,\widetilde{\vec{V}}_{j-1}>}$ to $<\widetilde{\vec{V}}_1,\ldots ,\widetilde{\vec{V}}_{j-1}>$ with respect to $\langle\cdot ,\cdot\rangle_{Q,N}$, \, \par\noindent $\displaystyle L_j:=\bigl\|\bigl(1-P_{<\widetilde{\vec{V}}_1,\ldots ,\widetilde{\vec{V}}_{j-1}>}\bigr)\widetilde{\vec{V}}_j\bigr\|_\Xi  >0$, and $\displaystyle\frac{\displaystyle\max _{r\in\{1,\ldots ,j-1\}}\|\vec{u}_r\|_{Q,N}}{\, L_j \, \epsilon}$ is finite, which guarantees the existence of an integer $K_j$ such that $\displaystyle\,\, 2^{K_j} L_j \, \epsilon \ge \sqrt{\frac{1}{2}}\max _{r\in\{1,\ldots ,j-1\}}\!\!\|\vec{u}_r\|_{Q,N}$ 

On the other hand, from Lemma \ref{lemma:max_distance_lattice} with $n=1$, for $\ell =1,2,\ldots j-1$, 
\begin{eqnarray*}
\inf_{\vec{g}\in {\rm Latt}(\widetilde{\vec{V}}_\ell )}\|P_{<\widetilde{\vec{V}}_\ell >}(2^{K_j}\widetilde{\vec{V}}_j)-\vec{g}\|\le \sqrt{\frac{1}{2}} \,\, \|\widetilde{\vec{V}}_\ell \|_{Q,N}\le \sqrt{\frac{1}{2}} \max _{r\in\{1,\ldots j-1\}}\|\widetilde{\vec{V}}_r\|\, .
\end{eqnarray*}
These inequalities result in 
\begin{eqnarray*}
\inf_{\vec{g}\in {\rm Latt}(\widetilde{\vec{V}}_\ell )}\|P_{<\widetilde{\vec{V}}_\ell >}(2^{K_j}\widetilde{\vec{V}}_j)-\vec{g}\|\le \, \epsilon \, \bigl\|P_{<\widetilde{\vec{V}}_\ell >}(2^{K_j}\widetilde{\vec{V}}_j)\bigr\|_{Q,N} \, 
\end{eqnarray*}
which shows the existence of a vector $\vec{g}_\ell $ in ${\rm Latt}(\widetilde{\vec{V}}_\ell )$ such that  \par\noindent $\displaystyle \frac{\langle (2^{K_j}\widetilde{\vec{V}}_j)-\vec{g}_\ell , \,\, \vec{g}_\ell \rangle _{Q,N}}{\|(2^{K_j}\widetilde{\vec{V}}_j)-\vec{g}_\ell \|_{Q,N}\cdot\,\|\vec{g}_\ell \|_{Q,N}}\le \epsilon $, for $\ell =1,2,\ldots j-1$. Because the relation \par\noindent $\vec{v}_j\in {\rm Latt}(\widetilde{\vec{V}}_1,\widetilde{\vec{V}}_2,\ldots ,\widetilde{\vec{V}}_{j-1},2^{\kappa_j}\widetilde{\vec{V}}_{j})$ with the frequency $\kappa _j$ of the substitution  $\vec{v}_j\leftarrow2\vec{v}_j$ made already for $\vec{v}_j$ is always guaranteed at any step of the iteration of {\bf Q2}, this fact implies that the substitution  $\vec{v}_j\leftarrow2\vec{v}_j$ in the iteration of {\bf Q2} in Table \ref{tbl:2} cannot be repeated more than $K_j+1$ times for $\vec{v}_j$. 
Since the procedures other than the substitution $\vec{v}_j\leftarrow2\vec{v}_j$ 
do not increase $\displaystyle s_{Q,N}(\vec{v}_1,\ldots ,\vec{v}_D):=\sum_{r=1}^D \|\vec{v}_r\|_{Q,N}$, 
at any step of the iteration of {\bf Q2}, $s_{\vec{v}_1,\ldots ,\vec{v}_D}$ is always bounded by  $\displaystyle \tilde{s}:=\sum_{r=1}^D 2^{K_r+1} \|\widetilde{\vec{V}}_r\|_{Q,N}$. 
Moreover, any process in {\bf Q2} gives a vector $\vec{v}_j$ in ${\rm Latt}(\widetilde{\vec{V}}_1,\ldots ,\widetilde{\vec{V}}_D)$. 
Hence, at any step of the iteration of {\bf Q2}, $ (\vec{v}_1, \ldots , \vec{v}_D)$ belongs to  $\tilde{T}_{Q,N}^{(\tilde{s}\, )}(\widetilde{\vec{V}}_1,\ldots ,\widetilde{\vec{V}}_D)$ 
 which is a finite set, where $\tilde{T}_\Xi^{(s)}(\vec{u}_1,\ldots ,\vec{u}_n)$ has been defined in Definition \ref{definition:lattice_finite}.
Similarly to this, because any procedure in {\bf Q1} does not increase $ s_{Q,N}(\vec{v}_1,\ldots ,\vec{v}_D)$ and gives a vector $v_j$ in ${\rm Latt}(\vec{F}_{\rm int.}^{(1)},\ldots ,\vec{F}_{\rm int.}^{(D)})$, at any step in the iteration of {\bf Q1}, $ (\vec{v}_1, \ldots , \vec{v}_D)$ belongs to the finite set $\tilde{T}_{Q,N}^{(\sigma)}(\vec{F}_{\rm int.}^{(1)},\ldots ,\vec{F}_{\rm int.}^{(D)})$ with $\displaystyle \sigma:=s_{Q,N}(\vec{F}_{\rm int.}^{(1)},\ldots ,\vec{F}_{\rm int.}^{(D)})=\sum_{r=1}^D \|\vec{F}_{\rm int.}^{(r)}\|_{Q,N}$.

Then, the procedures in Table \ref{tbl:2} except for finitely many (not greater than $K_j+1$ for each $j$) times of the execution of $\vec{v}_j\leftarrow 2\vec{v}_j $ $(j=2,3,\ldots ,D)$ do not increase  $s_{Q,N}(\vec{v}_1,\ldots ,\vec{v}_D)$ and the substitutions  $\displaystyle \vec{v}_j\leftarrow\vec{v}_j-\Bigl[\frac{\langle\vec{v}_j, \vec{v}_\ell \rangle_{Q,N}}{\langle\vec{v}_\ell , \vec{v}_\ell \rangle_{Q,N}} \Bigr]_\mathbb{C} \vec{v}_\ell $ in {\bf Q1} and {\bf Q2} always decrease  $s_{Q,N}(\vec{v}_1,\ldots ,\vec{v}_D)$,  unless $\vec{v}_j$ is changed. Since the sets $\tilde{T}_{Q,N}^{(\tilde{s}\, )}(\widetilde{\vec{V}}_1,\ldots ,\widetilde{\vec{V}}_D)$ and  $\tilde{T}_{Q,N}^{(\sigma)}(\vec{F}_{\rm int.}^{(1)},\ldots ,\vec{F}_{\rm int.}^{(D)})$ are finite, this process is also carried out only finitely many times. Hence, the procedures in Table \ref{tbl:2} 
%terminate 
halt  
in a finite number of steps. 
\hfill\endproof
\par
\begin{rem}\begin{indention}{.8cm}\rm
In spite of the complications of the above proof of 
halt,  
%termination, 
the 
%number 
amount  
of calculations required for the iterations of the substitution processes in Table \ref{tbl:2} is much smaller than the 
%number 
amount  
of calculations required for the calculations of the inner products themselves; this has been observed empirically. Hence, the iterations of the substitution processes in Table \ref{tbl:2} are not `bottlenecks' of our method at all, though proof of their 
%termination 
halting  
is 'logically' necessary.
\end{indention}\end{rem}

\par\noindent Even for the cases when $\dim V\cap \ell ^2(\mathbb{Z}^+ )\le 2$, we can prove %termination 
the halting  
of the processes in Table \ref{tbl:4} and Table \ref{tbl:5} for the process {\bf Step 3$^\prime $} in a similar manner to this, with some modifications, because the basic structure of the procedures is almost the same as for the one-dimensional case. Here, we omit it because of the complexity of the notations.

\section{Upper bound on errors}
\label{sec:ul}
In numerical methods, it is very important to know the precision of the results. 
Here we will give an error bound for our method.

For this, we begin with two lemmata, with  $R$ as in (\ref{eqn:w_n}) 
and \par\noindent $\displaystyle \Delta _K := \sup _{\vec{f}\in V\cap \ell ^2(\mathbb{Z}^+ )\backslash \{ 0\}} \frac{\|(1-\Pi _K)\vec{f}\|_{\ell ^2}}{\|\Pi _K\vec{f}\|_{\ell ^2}}$\,  which is the upper limit of the truncation error, normalized in the subspace  $<e_0,e_1, \ldots e_K>$.
\begin{lemma}
\label{lemma:error_bound_1}
Assume that $\dim V\cap\ell ^2(\mathbb{Z}^+)=1$. 
%Let $P_{U\cap \ell ^2}$ be the orthogonal projector to $U\cap\ell ^2(\mathbb{Z}^+)$ with respect to  $\langle \cdot ,\cdot\rangle _{\ell ^2,K}$. 
For a nonzero vector $\vec{w}$ in $\Pi _N V$, let  $W:=\{ a \vec{w}\, | \, a\in \mathbb{C} \}$ and $U$ be a subspace of $\Pi _N V$ such that $\Pi _N V=W\oplus U$ and  $\displaystyle \sup_{\vec{u}\in U\backslash \{ 0\}}\frac{|( \vec{w},\, \vec{u} )_{\ell ^2,K} |}{|| \vec{w} ||_{\ell ^2,K} \cdot || \vec{u} ||_{\ell ^2,K} } \, \le \,\, \xi $.  Moreover, let  $\displaystyle \, C_1:=\frac{\,\, \|\vec{w}\|_{Q,N}\, }{\|\vec{w}\|_{\ell ^2,K}}\, $ and $\displaystyle \Gamma_1:=\inf_{\vec{u}\in U\backslash \{ 0\}} \frac{\|\vec{u}\|_{Q,N} }{\|\vec{u}\|_{\ell ^2,K} }$.\, $\,\,\,\, $ If $\displaystyle\, \Gamma_1>\, \xi C_1$ and   $\langle\vec{f}, \vec{w}\rangle _{\ell, K}\ne 0$ for  $ \vec{f}\in V\cap\ell ^2(\mathbb{Z}^+)\backslash \{ 0\}$, then 
\begin{eqnarray*}
\frac{\displaystyle\inf_{\vec{g}\in V\cap {\ell ^2}(\mathbb{Z}^+ )} \|\Pi_K\vec{w}-\vec{g}\|_{\ell ^2}}{\|\Pi_K\vec{w}\|_{\ell ^2}} 
\le \frac{C_1+1}{\Gamma_1-\xi C_1}+\Bigl(1+\frac{R}{\Gamma_1-\xi C_1}\Bigr)\Delta_K\, .
\end{eqnarray*}
\end{lemma} \par\noindent 
This theorem can be generalized to the cases where $\dim V\cap\ell ^2(\mathbb{Z}^+)\ge 2$, as follows:
\begin{lemma}
\label{lemma:error_bound}Let  $P_{\Pi _N (V\cap\ell ^2(\mathbb{Z}^+))}$ be the orthogonal projector on $\Pi _N V$ to $\Pi _N (V\cap\ell ^2(\mathbb{Z}^+))$  with respect to  $\langle \cdot ,\cdot\rangle _{\ell ^2,K}$.  
%Let $P_{U\cap \ell ^2}$ be the orthogonal projector to $U\cap\ell ^2(\mathbb{Z}^+)$ with respect to  $\langle \cdot ,\cdot\rangle _{\ell ^2,K}$. 
With a positive integer $\widetilde{D}$ not greater than $\dim V\cap\ell ^2(\mathbb{Z}^+)$ and $\widetilde{D}$ linearly independent vectors $\vec{w}_1\! , \vec{w}_2\! , \dots , \vec{w}_{\widetilde{D}}$ in  $W_{\widetilde{D}}:=<\vec{w}_1, \vec{w}_2, \dots , \vec{w}_{\widetilde{D}}>$ and $U_{\widetilde{D}}$ be a subspace of $\, \Pi _N V$ such that  $\Pi _N V=W_{\widetilde{D}}\oplus U_{\widetilde{D}}$ and \par\noindent  $\displaystyle \sup_{\vec{w}\in W_{\widetilde{D}}\backslash \{ 0\}, \, \vec{u}\in U_{\widetilde{D}}\backslash \{ 0\}}\frac{|( \vec{w},\, \vec{u} )_{\ell ^2,K} |}{|| \vec{w} ||_{\ell ^2,K} \cdot || \vec{u} ||_{\ell ^2,K} } \, \le \,\, \xi $.  Moreover, let  $\displaystyle C_{\widetilde{D}}:=\sup_{\vec{w}\in W_{\widetilde{D}}\backslash \{ 0\}} \frac{\|\vec{w}\|_{Q,N} }{\|\vec{w}\|_{\ell ^2,K} }$ and  $\displaystyle \Gamma_{\widetilde{D}}:=\inf_{\vec{u}\in U_{\widetilde{D}}\backslash \{ 0\}} \frac{\|\vec{u}\|_{Q,N} }{\|\vec{u}\|_{\ell ^2,K} }$.\, $\,\,\,\, $ If $\, \Gamma_{\widetilde{D}}>\xi C_{\widetilde{D}}$ and $\bigl( (1-P_{\Pi _N (V\cap\ell ^2(\mathbb{Z}^+))}\Pi _N V\bigr) \cap W_{\widetilde{D}}=\{ 0\}$, then 
\begin{eqnarray*}
 \sup_{\vec{w}\in W_{\widetilde{D}}\backslash \{ 0\}} \!\!\frac{\displaystyle\inf_{\vec{g}\in V\cap {\ell ^2}(\mathbb{Z}^+ )}\|\Pi_K\vec{w}-\vec{g}\|_{\ell ^2}}{\|\Pi_K\vec{w}\|_{\ell ^2}} 
\le \frac{C_{\widetilde{D}}+1}{\Gamma_{\widetilde{D}}-\xi C_{\widetilde{D}}}+\Bigl(1+\frac{R}{\Gamma_{\widetilde{D}}-\xi C_{\widetilde{D}}}\Bigr)\Delta_K\, .
\end{eqnarray*}
\end{lemma}  
\par\noindent\noindent{\em Proof of Lemma }\ref{lemma:error_bound}: \quad
\rm 
Let  $\vec{w}$ be a nonzero vector in $W_{\widetilde{D}}$ $(\subset \Pi _N V)$. Since \par\noindent  $P_{\Pi _N (V\cap\ell ^2(\mathbb{Z}^+))}\vec{w}\in \Pi _N (V\cap\ell ^2(\mathbb{Z}^+))$, there exists a vector $\vec{g}_0$ in $V\cap \ell ^2(\mathbb{Z}^+ )$ such that $P_{\Pi _N (V\cap\ell ^2(\mathbb{Z}^+))}\vec{w}=\Pi_K\vec{g}_0 $. From the definition of $\Delta _K$, 
\begin{eqnarray*}
\|\Pi_K\vec{w}-\vec{g}_0 \|_{\ell ^2}\!\!&=&\! \|\Pi_K(\vec{w} - \vec{g}_0 )\|_{\ell ^2,K}\!\!+\!\|(1-\Pi_K)\vec{g}_0 \|_{\ell ^2} \\ 
\!\!&\le &\!\|\vec{w}-P_{\Pi _N (V\cap\ell ^2(\mathbb{Z}^+))}\vec{w} \|_{\ell ^2,K}\!\!+\!\Delta_K \|\vec{g}_0 \|_{\ell ^2,K}\!\! \\ 
&=&\|(1-P_{\Pi _N (V\cap\ell ^2(\mathbb{Z}^+))})\vec{w} \|_{\ell ^2,K}+\Delta_K \|P_{\Pi _N (V\cap\ell ^2(\mathbb{Z}^+))}\vec{w} \|_{\ell ^2,K}\\  
&\le &\|(1-P_{\Pi _N (V\cap\ell ^2(\mathbb{Z}^+))})\vec{w} \|_{\ell ^2,K}+\Delta_K \|\vec{w} \|_{\ell ^2,K}  .
\end{eqnarray*}
Since $\|\vec{w}\|_{\ell ^2,K}=\|\Pi_K\vec{w}\|_{\ell ^2}$, therefore, proof of the statement 
\begin{eqnarray*}
\sup_{\vec{w}\in W_{\widetilde{D}}\backslash \{ 0\}} \!\!\frac{\|(1-P_{\Pi _N (V\cap\ell ^2(\mathbb{Z}^+))})\vec{w}\|_{\ell ^2,K}}{\|\vec{w}\|_{\ell ^2,K}} 
\le \frac{C_{\widetilde{D}}+R\Delta _K+1}{\Gamma_{\widetilde{D}}-\xi C_{\widetilde{D}}}
\end{eqnarray*}
suffices. In the following parts of this proof, we show this statement.

Let $P_{W_{\widetilde{D}}}$ be the orthogonal projector to $W_{\widetilde{D}}$ with respect to  $\langle \cdot ,\cdot\rangle _{\ell ^2,K}$, and let \par\noindent $W_{\widetilde{D}}^\bot :=(1-P_{W_{\widetilde{D}}})\Pi _N V$. 
If $\bigl( (1-P_{\Pi _N (V\cap\ell ^2(\mathbb{Z}^+))})\, \Pi _N V\bigr) \cap W_{\widetilde{D}}=\{ 0\}$, then  $\Pi _N (V\cap\ell ^2(\mathbb{Z}^+))+W_{\widetilde{D}}^\bot =\Pi _N V$ because the orthogonal complements (in $\Pi _N V$) of $(1-P_{\Pi _N (V\cap\ell ^2(\mathbb{Z}^+))})\, \Pi _N V $ and $W_{\widetilde{D}}$ are $\Pi _N (V\cap\ell ^2(\mathbb{Z}^+))$ and $W_{\widetilde{D}}^\bot$, respectively, with respect to  $\langle \cdot ,\cdot\rangle _{\ell ^2,K}$. Hence, for any nonzero vector $\vec{w}$ in $W_{\widetilde{D}}$, there exist vectors $\vec{v}\in \Pi _N (V\cap\ell ^2(\mathbb{Z}^+))$ and $\vec{w}^\bot \in W_{\widetilde{D}}^\bot $ such that  $\vec{w}=\vec{v}+\vec{w}^\bot$. Then $\vec{w}=P_{W_{\widetilde{D}}}\vec{v}$. Since $\langle \vec{w}^\bot, \vec{w}\, \rangle_{\ell ^2,K}=0$, the inequalities $\|\vec{v}\|_{\ell ^2,K}\ge \|\vec{w}\|_{\ell ^2,K}>0$ 
 hold. 

If $\vec{v}\in W_{\widetilde{D}}$, then $\vec{w}^\bot =0$ and $\vec{w}=\vec{v}$. This  results in \par\noindent $\|(1-P_{\Pi _N (V\cap\ell ^2(\mathbb{Z}^+))})\vec{w}\,\|_{\ell ^2}\le \|(1-P_{\vec{v}})\vec{w}\,\|_{\ell ^2,K}=0$, which satisfies the statement of the lemma obviously. Therefore, in the following, we prove the lemma for the cases  $\vec{v}\notin W_{\widetilde{D}}$

From the condition $\Pi _N V=W_{\widetilde{D}}\oplus U_{\widetilde{D}}$, there exists vectors $\widehat{\vec{w}}\in W_{\widetilde{D}}$ and $\widehat{\vec{u}}\in U_{\widetilde{D}}$ such that $\vec{v}=\widehat{\vec{w}}+\widehat{\vec{u}}$ and $\displaystyle \frac{|( \vec{f},\, \widehat{\vec{u}} )_{\ell ^2,K} |}{|| \vec{f} ||_{\ell ^2,K} \cdot || \widehat{\vec{u}} ||_{\ell ^2,K} }  \le  \xi $ for any $\vec{f}\in W_{\widetilde{D}}$. 
Hence, \par\noindent $\displaystyle\|P_{W_{\widetilde{D}}}\widehat{\vec{u}}\|_{\ell ^2,K}=\frac{\langle P_{W_{\widetilde{D}}}\widehat{\vec{v}}, \widehat{\vec{u}}\rangle _{\ell ^2,K}}{\|P_{W_{\widetilde{D}}}\widehat{\vec{u}}\|_{\ell ^2,K}}\le\xi \, \|\widehat{\vec{u}}\|_{\ell ^2,K}.$   
From the assumption that $\vec{v}\notin W_{\widetilde{D}}$, the vector $\widehat{\vec{u}}$ is not $0$. 
Hence, the ratios  $\displaystyle z:=\frac{\|\widehat{\vec{w}}-\vec{w}\|_{\ell ^2,K}}{\|\widehat{\vec{u}}\|_{\ell ^2,K}}$ and $\displaystyle d:=\frac{\|\vec{w}^\bot\|_{\ell ^2,K}}{\|\widehat{\vec{u}}\|_{\ell ^2,K}}$ are well-defined. 
Then, the inequality $d\le 1$ holds because $-\vec{w}^\bot$ is the perpendicular from $\vec{v}$ to $W_{\widetilde{D}}$.
Since $\widehat{\vec{w}}-\vec{w}=P_{W_{\widetilde{D}}}(\widehat{\vec{w}}-\vec{v})=P_{W_{\widetilde{D}}}\widehat{\vec{u}}$, the above inequality results in \par\noindent $\displaystyle z=\frac{\|P_{W_{\widetilde{D}}}\widehat{\vec{u}}\|_{\ell ^2,K}}{\|\widehat{\vec{u}}\|_{\ell ^2,K}}\le \xi $.

The trigonometric inequality $\|\widehat{\vec{u}}\|_{Q,N}\le \|\widehat{\vec{w}}\|_{Q,N}+\|\vec{v}\|_{Q,N}$ and the definitions of $C_{\widetilde{D}}$ and $\Gamma_{\widetilde{D}}$ result in the inequalities 
\begin{eqnarray*}
 \Gamma_{\widetilde{D}}\|\widehat{\vec{u}}\|_{\ell ^2,K}&\le& C_{\widetilde{D}}\|\widehat{\vec{w}}\|_{\ell ^2,K}+(R\Delta_K+1) \|\vec{v}\|_{\ell ^2,K} \\ &\le& C_{\widetilde{D}}\|\widehat{\vec{w}}-\vec{w}\|_{\ell ^2,K}+C_{\widetilde{D}}\|\vec{w}\|_{\ell ^2,K}+(R\Delta_K+1) \|\vec{v}\|_{\ell ^2,K}\, \, ,
\end{eqnarray*}
because $\displaystyle \frac{\|\vec{v}\|_{Q,N}}{\|\vec{v}\|_{\ell ^2,K}}\le \frac{\displaystyle\sum_{n=0}^K |v_n|^2+ \displaystyle\sum_{n=K+1}^\infty R|v_n|^2}{\displaystyle\sum_{n=0}^K |v_n|^2}=(R\Delta_K+1)\,$. Hence, \par\noindent $\displaystyle \Gamma_{\widetilde{D}}\|\widehat{\vec{u}}\|_{\ell ^2,K}\le zC_{\widetilde{D}} \|\widehat{\vec{u}}\|_{\ell ^2,K}+ (C_{\widetilde{D}}+R\Delta_K+1) \|\vec{v}\|_{\ell ^2,K}$, because $\|\vec{w}\|_{\ell ^2,K}\le \|\vec{v}\|_{\ell ^2,K}$.  From the condition $\Gamma_{\widetilde{D}}-\xi C_{\widetilde{D}}>0$ and the inequality $z\le\xi$ above, $\Gamma_{\widetilde{D}}-z C_{\widetilde{D}}>0$. Hence, 
\begin{eqnarray*}
 \|\vec{w}^\bot\|_{\ell ^2,K}\! =d\,\|\widehat{\vec{u}}\|_{\ell ^2,K}\! \le \|\widehat{\vec{u}}\|_{\ell ^2,K}\!\le \frac{C_{\widetilde{D}}+R\Delta_K+1}{\Gamma_{\widetilde{D}}-zC_{\widetilde{D}}}\|\vec{v}\|_{\ell ^2,K}\!\le \frac{C_{\widetilde{D}}+R\Delta_K+1}{\Gamma_{\widetilde{D}}-\xi C_{\widetilde{D}}}\|\vec{v}\|_{\ell ^2,K} , 
\end{eqnarray*}
i.e. $\displaystyle \frac{\|\vec{w}^\bot\,\|_{\ell ^2,K}}{\|\vec{v}\|_{\ell ^2,K}}\le \frac{C_{\widetilde{D}}+R\Delta_K+1}{\Gamma_{\widetilde{D}}-\xi C_{\widetilde{D}}}$ where $\|\vec{v}\|_{\ell ^2,K}> 0$ is guaranteed as has been shown above. 

Let $P_{\vec{v}}$ be the orthogonal projector to the subspace $S_{\vec{v}}:=\{a\vec{v}\, |\, a\in\mathbb{C}\}$. Then from a geometrical comparison of length among the perpendiculars, the relation  $\displaystyle \frac{\|(1-P_{\Pi _N (V\cap\ell ^2(\mathbb{Z}^+))})\vec{w}\,\|_{\ell ^2,K}}{\|\vec{w}\,\|_{\ell ^2,K}}\le \frac{\|(1-P_{\vec{v}})\vec{w}\,\|_{\ell ^2,K}}{\|\vec{w}\,\|_{\ell ^2,K}}=\frac{\|\vec{w}^\bot\,\|_{\ell ^2,K}}{\|\vec{v}\|_{\ell ^2,K}}$ holds, because \par\noindent $S_{\vec{v}}\subset \Pi _N (V\cap\ell ^2(\mathbb{Z}^+))$ and $\vec{w}$ $(=\vec{v}+\vec{w}^\bot)$ is the orthogonal projection of $\vec{v}$ to the subspace $S_{\vec{w}}:=\{a\vec{w}\, |\, a\in\mathbb{C}\}$ with respect to  $\langle \cdot ,\cdot\rangle _{\ell ^2,K}$ as well as the projection $P_{W_{\widetilde{D}}}\vec{v}$ to $W_{\widetilde{D}}$.

Hence, the inequality $\displaystyle \frac{\|(1-P_{\Pi _N (V\cap\ell ^2(\mathbb{Z}^+))})\vec{w}\,\|_{\ell ^2,K}}{\|\vec{w}\,\|_{\ell ^2,K}}\le \frac{C_{\widetilde{D}}+R\Delta_K+1}{\Gamma_{\widetilde{D}}-\xi C_{\widetilde{D}}}$ holds for any \par\noindent $\vec{w}\in W_{\widetilde{D}}$. The lemma has been proved by the combination of this fact and the discussion in the first part of this proof.
\hfill\endproof 
\par \par\noindent 

The proof of Lemma \ref{lemma:error_bound_1} is just the proof above with $\widetilde{D}=1$.

These propositions lead to the following theorem, which gives the upper bound of the error:

\begin{theorem}
\label{thm:error_bound}
Let  $P_{\Pi _N (V\cap\ell ^2(\mathbb{Z}^+))}$ be the orthogonal projector on $\Pi _N V$ to $\Pi _N (V\cap\ell ^2(\mathbb{Z}^+))$ with respect to  $\langle \cdot ,\cdot\rangle _{\ell ^2,K}$. 
Let  $T:=<\vec{G}^{\, (1)},\vec{G}^{\, (2)}, \ldots ,\vec{G}^{\, (D_{\ell ^2})}>$ $(=T_{D_{\ell ^2}}$ in Section \ref{sec:multidim} $)$  
and $R:=<\vec{v}_1^{<D_{\ell ^2}>},\vec{v}_2^{<D_{\ell ^2}>},\ldots ,\vec{v}_{D-D_{\ell ^2}}^{<D_{\ell ^2}>}>$ $(=R_{D_{\ell ^2}}$ in Section \ref{sec:multidim} $)$   be  subspaces of $\Pi _N V$ such that  $\Pi _N T=W\oplus R$ and suppose that  the quasi-orthogonalities $(a), (b)$ and $(d)$ in Lemma \ref{lemma:quasi_ort} are satisfied for $\tilde{d}=D_{\ell ^2}$ .  Moreover, let \par\noindent $\displaystyle C:=\sup_{\vec{w}\in T\backslash \{ 0\}} \frac{\|\vec{w}\|_{Q,N} }{\|\vec{w}\|_{\ell ^2,K} }$ and $\displaystyle \Gamma:=\inf_{\vec{u}\in R\backslash \{ 0\}} \frac{\|\vec{u}\|_{Q,N} }{\|\vec{u}\|_{\ell ^2,K} }$.\, $\,\,\,\, $ If $\, \Gamma>\xi (g) \, C$ with  
\begin{eqnarray*}
\xi (g):= \frac{1}{g} \, \, \sqrt{\frac{D\, D_{\ell ^2}}{\left(1- \frac{\displaystyle D-1}{\displaystyle g}\right)\Biggl(1- \frac{\textstyle D_{\ell ^2}-1}{\displaystyle g}\sqrt{\frac{\displaystyle D  }{\displaystyle  1- \frac{\displaystyle D-1}{\displaystyle g}}}\Biggr)}} \,\,\,    
\end{eqnarray*}
and $\bigl( (1-P_{\Pi _N (V\cap\ell ^2(\mathbb{Z}^+))})\, \Pi _N V\bigr) \cap W_T=\{ 0\}$, then the inequality 
\begin{eqnarray*} 
 \sup_{\vec{w}\in T\backslash \{ 0\}}\frac{\displaystyle\inf_{\vec{g}\in V\cap \ell ^2(\mathbb{Z}^+ )}\|\Pi_K\vec{w}-\vec{g}\|_{\ell ^2}}{\|\Pi_K\vec{w}\|_{\ell ^2}} 
\le \frac{C+1}{\Gamma-\xi (g)\, C}+\Bigl(1+\frac{R}{\Gamma-\xi (g) \, C}\Bigr)\Delta_K \, 
\end{eqnarray*}
holds.
\end{theorem} \par\noindent 
\par\noindent\noindent{\em Proof of Theorem }\ref{thm:error_bound}: \quad
\rm From Lemma \ref{lemma:min_dir_cos_2} and the quasi-orthogonalities (b) and (d), with $n=D-D_{\ell ^2}$, $\vec{u}_m=\vec{v}_m^{<D_{\ell ^2}>}$ $(m=1,2,\ldots ,D-D_{\ell ^2})$ and $\vec{v}=\vec{G}^{\, (\ell )}$, 
the inequality $\, \displaystyle  \frac{|\langle \vec{w},\vec{G}^{\, (\ell )}\rangle _{\ell ^2,K}|}{\|\vec{w}\|_{\ell ^2,K} \|\vec{G}^{\, (\ell )}\|_{\ell ^2,K}} \le \frac{1}{g} \, \, \sqrt{\frac{D\,  }{1- \frac{D-1}{g}}} $ holds for any vector in $\vec{w}\in R$ and for $\ell =1,2,\ldots D_{\ell ^2}$. Using Lemma \ref{lemma:min_dir_cos_2} again with $n=D_{\ell ^2}$ and  $\vec{u}_m=\vec{G}^{\, (m )}$ $(m=1,2,\ldots ,D_{\ell ^2})$, we have the inequality $\displaystyle \sup_{\vec{w}\in R\backslash \{ 0\},\, \vec{u}\in T\backslash \{ 0\}}\frac{|( \vec{w},\, \vec{u} )_{\ell ^2,K} |}{|| \vec{w} ||_{\ell ^2,K}  || \vec{u} ||_{\ell ^2,K} } \le \xi (g)$.

Hence, from Lemma \ref{lemma:error_bound} with $W_{\widetilde{D}}=T$, $V_{\widetilde{D}}=R$ and $\xi =\xi (g)$, we have the statement of the theorem 
\hfill \endproof 
\par \par 
The condition $\bigl( (1-P_{\Pi _N (V\cap\ell ^2(\mathbb{Z}^+))})\, \Pi _N V\bigr) \cap T=\{ 0\}$
 is satisfied except for very special cases with `artificially  bad' approximation where a linear combination of obtained vectors $\Pi _K \vec{G}^{\, (1)},\Pi _K \vec{G}^{\, (2)}, \ldots ,\Pi _K \vec{G}^{\, (D_{\ell ^2})}$ for solutions is orthogonal to the space $V\cap \ell ^2(\mathbb{Z}^+ )$ of true solutions (!) with respect to $\langle \cdot ,\cdot\rangle _{\ell ^2,K}$. Since  the algorithm is used under guaranteed convergence of any nonzero linear combination of $\Pi _K \vec{G}^{\, (1)},\Pi _K \vec{G}^{\, (2)}, \ldots ,\Pi _K \vec{G}^{\, (D_{\ell ^2})}$, which has been shown in the last part of Section \ref{sec:conv_multidim}, we may neglect this condition in a practical sense. 

For sufficiently large $g$, we can choose very small $\xi (g)$. Moreover, there is a method to obtain an upper limit of $C$ and a lower limit of $\Gamma$ using only the numerical results, without any knowledge about the true solutions of the differential equation, as follows: 

The use of  Lemma \ref{lemma:max_ratio_pl} with $U=\Pi _N V$, $(\cdot ,\cdot )_\Lambda=\langle \cdot ,\cdot \rangle _{Q,N}$,  $(\cdot ,\cdot )_\Xi=\langle \cdot ,\cdot \rangle _{\ell ^2,K}$ and $\vec{f}_m=\vec{G}^{\, (m)}$ $(m=1,2,\ldots D_{\ell ^2})$ leads us to the inequality \par\noindent $\displaystyle C\le \sqrt{\frac{\displaystyle \sum_{m=1}^{D_{\ell ^2}}\frac{\|\vec{G}^{\, (m)}\|_{Q,N}}{\|\vec{G}^{\, (m)}\|_{\ell ^2,K}}}{1-\frac{\displaystyle D_{\ell ^2}-1}{\displaystyle g}\sqrt{\displaystyle \frac{\displaystyle D\,  }{\displaystyle 1- \frac{\displaystyle D-1}{\displaystyle g}}}}}$\, . 
On the other hand, the use of Lemma \ref{lemma:max_ratio_pl} with $U=\Pi _N V$, $(\cdot ,\cdot )_\Lambda=\langle \cdot ,\cdot \rangle _{\ell ^2,K}$,  $(\cdot ,\cdot )_\Xi=\langle \cdot ,\cdot \rangle _{Q,N}$ and $\vec{f}_m=\vec{u}_m^{<D_{\ell ^2}>}$ $(=\vec{V}_{D_{\ell ^2}+m}^{<D_{\ell ^2}+1>})$  $(m=1,2,\ldots D-D_{\ell ^2})$ leads us the inequality 
$\displaystyle \frac{1}{\Gamma} \le \sqrt{\frac{\displaystyle \sum_{m=1}^{D-D_{\ell ^2}}\frac{\|\vec{u}_m^{<D_{\ell ^2}>}\|_{\ell ^2,K}}{\|\vec{u}_m^{<D_{\ell ^2}>}\|_{Q,N}}}{1-\frac{\displaystyle D-D_{\ell ^2}-1}{\displaystyle h}}}$,  and hence $\displaystyle \Gamma \ge \sqrt{\frac{1-\frac{\displaystyle D-D_{\ell ^2}-1}{\displaystyle h}}{\displaystyle \sum_{m=1}^{D-D_{\ell ^2}}\frac{\|\vec{u}_m^{<D_{\ell ^2}>}\|_{\ell ^2,K}}{\|\vec{u}_m^{<D_{\ell ^2}>}\|_{Q,N}}}}$\, .\,  
In these bounds, the factors  $1-\frac{\displaystyle D_{\ell ^2}-1}{\displaystyle g}\sqrt{\displaystyle \frac{\displaystyle D\,  }{\displaystyle 1- \frac{\displaystyle D-1}{\displaystyle g}}}$ and $1-\frac{\displaystyle D-D_{\ell ^2}-1}{\displaystyle h}$ are nearly equal to $1$, for  sufficiently large $g$ and $h$. In the usual situations with $\Gamma>>C$ and $\xi (g)<<1$ mentioned below,  the condition $\Gamma-\xi (g)\, C>0$ is satisfied.

These upper and lower bounds of $C$ and $\Gamma$, respectively, can be calculated numerically in the algorithm if we iterate {\bf Step 3.A} up to $\tilde{d}=D_{\ell ^2}$.  
Hence, this theorem gives an upper 
%limit 
bound  
for the $L_{(k)}^2$-norm of the error for any non-zero linear combination of the numerical solutions 
$\displaystyle f_{\rm app.}^{(d)}(x):=\sum_{n=0}^K(\Pi _K \vec{G}^{\, (d)})_n e_n(x) $ $(d=0,\ldots ,D_{\ell ^2})$, as the function  $\displaystyle\frac{C+1}{\Gamma-\xi (g)\, C}+\Bigl(1+\frac{R}{\Gamma-\xi (g)\, C}\Bigr)\Delta_K$ of the `worst magnitude of truncation error' $\Delta_K$. 

In usual circumstances  where the algorithm gives good convergence and with sufficiently large $g$, the parameters satisfy   
%$R\approx \Gamma>>C>>1$
$\Gamma>>R>>C>>1$  
and $\xi (g)<<1$. %Since $K$ is finite, we can also obtain a point-wise error bound, by the Schwarz inequality. 
Hence, 
%it 
the above upper bound  
is approximated roughly by $\displaystyle \frac{C}{\Gamma}
%+2\Delta _K
+\Delta _K 
$, which is approximately bounded by \par\noindent 
$\displaystyle \sqrt{\left({\displaystyle \sum_{m=1}^{D_{\ell ^2}}\frac{\|\vec{G}^{\, (m)}\|_{Q,N}}{\|\vec{G}^{\, (m)}\|_{\ell ^2,K}}}\right)\left({\displaystyle \sum_{m=1}^{D-D_{\ell ^2}}\frac{\|\vec{u}_m^{<D_{\ell ^2}>}\|_{\ell ^2,K}}{\|\vec{u}_m^{<D_{\ell ^2}>}\|_{Q,N}}}\right)}
%+2\Delta _K
+\Delta _K 
$. 

\begin{figure}[thb]
\begin{center}
\includegraphics{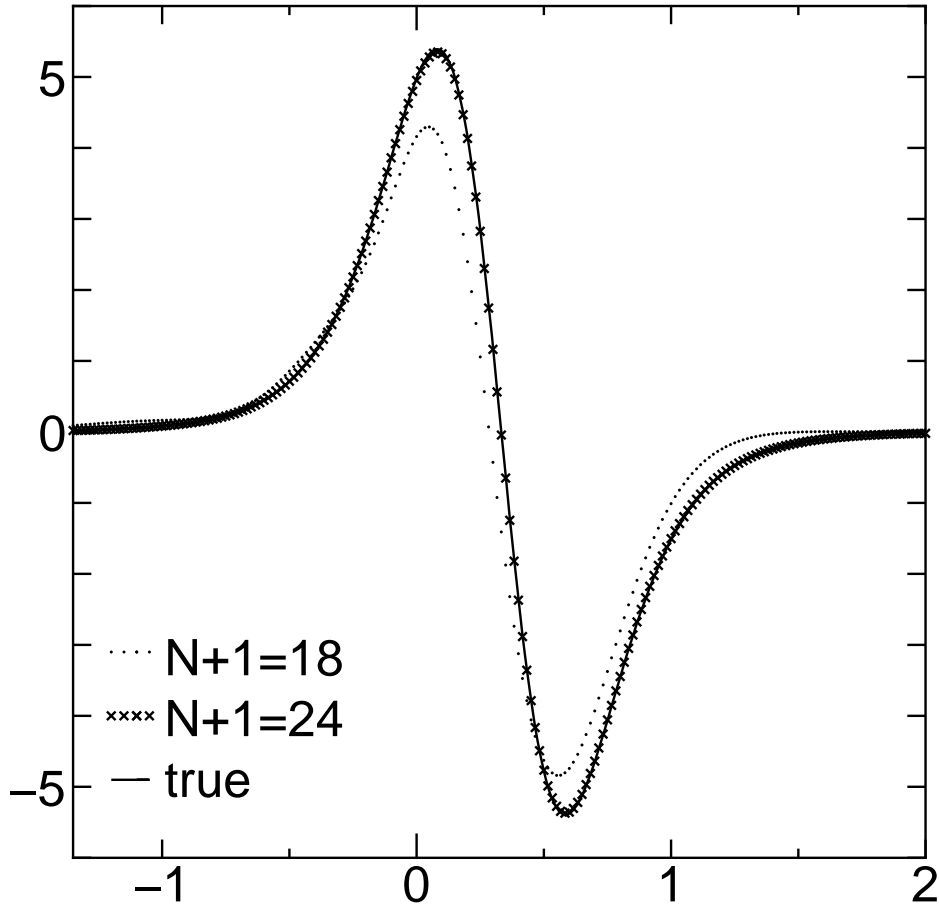}
\caption{Numerical results for the ODE $(9x^2-6x+5)f^{\prime\prime}+(90x-30)f^\prime+126f=0$}
\end{center}\end{figure}

\section{Numerical results}
\label{sec:nm}
In this section, we will give some numerical results of the proposed method. 
In Subsection \ref{subsec:np}, 
in order to show how accurate the results are, we will treat by intent some examples of ODEs which can be solved analytically because there we can compare the results with exact solutions up to arbitrary precision, though the proposed method can be widely applied to ODEs which can not be solved analytically at all. 
These results contain some examples where we are successful in attaining the accuracy with several hundreds or several thousand significant digits by an ordinary personal computer.

In Subsection \ref{subsec:nc}, we will compare theoretically the accuracy and the amount of calculation between some typical existing methods and the proposed method. However, direct numerical comparisons between the proposed method and existing methods are very difficult because usual existing methods with arbitrary-precision arithmetic (GNU Multi-Precision Library, for example) often require an astronomically large amount of calculations in order to attain such an extraordinarily high accuracy. Therefore, we will compare only the order of amount of calculations necessary for attaining a very large number of required significant digits.

\subsection{Numerical results by the proposed method}
\label{subsec:np}
In another paper\cite{paper1}, we have already shown how extraordinarily accurate results the proposed method gives for Weber's differential equation (Schr\"{o}dinger equation for harmonic oscillators) whose basic solutions in $L^2(\mathbb{R})$ are Hermite functions, and there we have given another example of a third order ODE and an example of the associate Legendre differential equation. 
In the former example of that paper, we were successful in obtaining the results where the ratio $f(x_0)/f(x_1)$ between values of solution function $f(x)$ coincides up to  $2599$ digits with the true ratio and the raio $f_n/f_{n^\prime}$ between coefficients $f_n$ ($n=0,1,\ldots $) in the expansion $f(x)=\sum_n f_n e_n(x)$ coincides up to $8783$ digits with the true ratio. There we observed that the number of significant digits are almost proportional a power of the dimension $N+1$ when the dimension ls very large. This implies that the amount of required calculations increases almost in a polynomial order of the number of required significant digits, from the reasons shown in Subsection \ref{subsec:nc} below. 

In this section, we will give two other examples than those. One is an example where we are successful in obtaining perfectly exact ratios between coefficients $f_n$ ($n=0,1,\ldots ,K$)  of the solution $f(x)$, and the other example is for an ODE whose true solutions are weighted associate Laguerre functions.
 
First, 
%Here
we will give numerical results for the second-order differential equation 
\begin{eqnarray*}
(9x^2-6x+5)f^{\prime\prime}+(90x-30)f^\prime +126f=0. 
\end{eqnarray*}
The space of its true solutions in $L_{(k_0)}^2(\mathbb{R} )$ is 
$\displaystyle \biggl\{ \frac{C (3x-1)}{\bigl((3x-1)^2+4\bigr)^4} \, \biggl|\, C\in\mathbb{C}\biggl\}$ for $k_0\le 3$. %and $\{ 0\}$ for $k\ge 4$. 
In Figure 3, the results with $k_0=2$, $N+1=18, 24$ 
$K=2\lfloor \frac{3N}{8} \rfloor + k_0$ and $J=2\lfloor \frac{7N}{16} \rfloor + k_0$ 
%and $J=2\lfloor \frac{7N}{16} \rfloor + k_0$   
are shown, under the normalization  $ \langle f, 
\frac{1}{2\pi}
(\psi_{k_0,0}+\psi_{k_0,-k_0-1})\rangle_{{\cal H}}=1$. The error of the result with $N+1=24$ is almost invisible there.

\begin{table}[bht]
\begin{center}
\begin{tabular}{|c|c|@{}c@{}|@{}c@{}|c|c|}
\hline 
$N\!\! +\!\! 1$ & ratio $\displaystyle \frac{f_2}{f_0}$ & \,  value of $\displaystyle {\rm Re} \, \frac{f_2}{f_0}$ & \,  value of  $\displaystyle  {\rm Im} \, \frac{f_2}{f_0}$ & \multicolumn{2}{|@{}c@{}|}{$\begin{array}{@{\,}ll}\mbox{number of } \\ \mbox{significant digits}\end{array}$} \\
\hline 
$\begin{array}{@{\,}ll}  18 \\ \\ \end{array}$  & $\displaystyle\frac{-59+31\,i}{33}$ &       $-1.78787878\ldots$ & $+0.93939393\ldots$ &  1  & 0 \\ 
\hline 
$\begin{array}{@{\,}ll}  24 \\ \\ \end{array}$  & $\displaystyle\frac{-2051+1976\, i}{1381}$ & $-1.48515568\ldots$ & $+1.43084721\ldots$& 2 & 2 \\
\hline 
$\begin{array}{@{\,}ll}  30 \\ \\ \end{array}$  & $\displaystyle\frac{-2249+2192\, i}{1520}$ & $-1.47960526\ldots$ & $+1.44210526\ldots$& 4 & 3 \\
\hline 
$\begin{array}{@{\,}ll}  36 \\ \\ \end{array}$  & $\displaystyle\frac{-2182+2126\, i}{1475}$ & $-1.47932203\ldots$ & $+1.44135593\ldots$ & 6  & 5 \\ 
\hline 
$\begin{array}{@{\,}ll}  48 \\ \\ \end{array}$  & $\displaystyle\frac{-42251+41166\, i}{28561}$ & $\begin{array}{cc} \mbox{(perfectly exact)} \\ -1.47932495\ldots\end{array}$ & $\begin{array}{cc} \mbox{(perfectly exact)} \\ +1.44133608\ldots\end{array}$ & \hspace{2.1mm} $\infty $ \hspace{2.1mm} & $\infty $\\
\hline 
\hline 
$\begin{array}{@{\,}ll}  {\rm true} \\ \\ \end{array}$  & $\displaystyle\frac{-42251+41166\, i}{28561}$ & $-1.47932495\ldots$  & $+1.44133608\ldots$ & $-$ & $-$\\
\hline 
\end{tabular}
\caption{Numerical results of the ratio $\displaystyle \frac{f_2}{f_0}$}
\end{center}
\label{tbl:6}
\end{table}

Moreover, we have investigated within how many digits the ratio between two coefficients $f_n$ and $f_{n^\prime}$ in the expansion $f(x)=\sum_n f_n e_n(x)$ coincides with the true ratio. For this differential equation, all the true ratios %$\displaystyle\frac{f_n}{f_{n'}}$ 
are rational-complex-valued. With $k=2$, true ratios are $\displaystyle \frac{f_1}{f_0} = 1$ and 
$\displaystyle \frac{f_2}{f_0} = \frac{-42251+41166 \, i}{28561}$, for example, which were obtained analytically using the computer algebra software system ``Mathematica''. As is shown in Table 6, for ${\rm Re}\, \displaystyle \frac{f_2}{f_0}$, the number of significant digits increases monotonically as $N$ increases. Moreover, the results with $N+1=48$  have attained the exact true ratio $\displaystyle\frac{-42251+41166\, i}{28561}$, where the ratio exactly coincides with the true value. In this example, with $N+1=48$, the other ratios among $f_n$ $(n=0,\ldots ,47)$ are all exact. The results with $N+1>48$ have also yielded this perfectly exact true ratio.
\begin{table}[ht]
\begin{center}
\begin{tabular}{|r|l|l|}
\hline
$n$ & real part of $f_n$ & imaginary part of $f_n$ \\
\hline
\hline
$0$ & $-$6.75984000378$\ldots\, $e$-$1& $\pm$0\\
\hline 
$1000$ & $-$2.84538929863$\ldots \, $e$-$271 & $-$9.83514249870$\ldots \, $e$-$272 \\ 
\hline 
$2000$ & +5.40708023241$\ldots \, $e$-$550 & $-$7.83428979204$\ldots \, $e$-$549 \\ 
\hline 
$3000$ & +8.49092503337$\ldots \, $e$-$827 & $-$1.66804004343$\ldots \, $e$-$827 \\ 
\hline 
$4000$ & +2.95369741917$\ldots \, $e$-$1105& +6.01269219021$\ldots \, $e$-$1105 \\ 
\hline 
\end{tabular}
\end{center}
\caption{Almost exponential decay of the coefficients}
\end{table}
Since the coefficients $f_n$ ($n\in \mathbb{Z}^+$) decay almost exponentially in this case as is shown in Table 7, we can easily attain the accuracy of the solution function $f(x)$ with more than one thousand significant digits. This exponential decay results from the relationship between the Fourier series and the basis functions used in our method (mentioned in Section \ref{sec:fb}). For example, the numerical results of the ratio $f(2)/f(1)$ with $N+1=10000$ coincides with the true ratio $\frac{5}{2}\left(\frac{8}{29}\right)^4=1.44779797562779150\ldots \times 10^{-2}$ up to 1942 significant digits. 

For ODEs where the ratios between coefficients are irrational, we can not attain such perfectly exact  ratios, because the results of our method are always rational-(complex-)valued. However, as has been shown in \cite{paper1}, the results have been successful in that very good rational approximations of the true irrational ratios were obtained, with 8783 significant digits 
as was mentioned above, for example.
%Other examples are given in~\cite{paper1}, which have irrational ratios between the coefficients. There the number of significant digits of the ratios increases monotonically as $N$ increases, one of which gives $340$ significant digits with $N+1=7000$. Though they cannot attain the exact true 
%values 
%ratios  
%because the results of our method are always rational-(complex-)valued, the results have been successful in that very good rational approximations of the true irrational ratios were obtained. 

\begin{figure}[tb]
\begin{center}
\begin{minipage}{0.4\linewidth}
\scalebox{.5}{\includegraphics{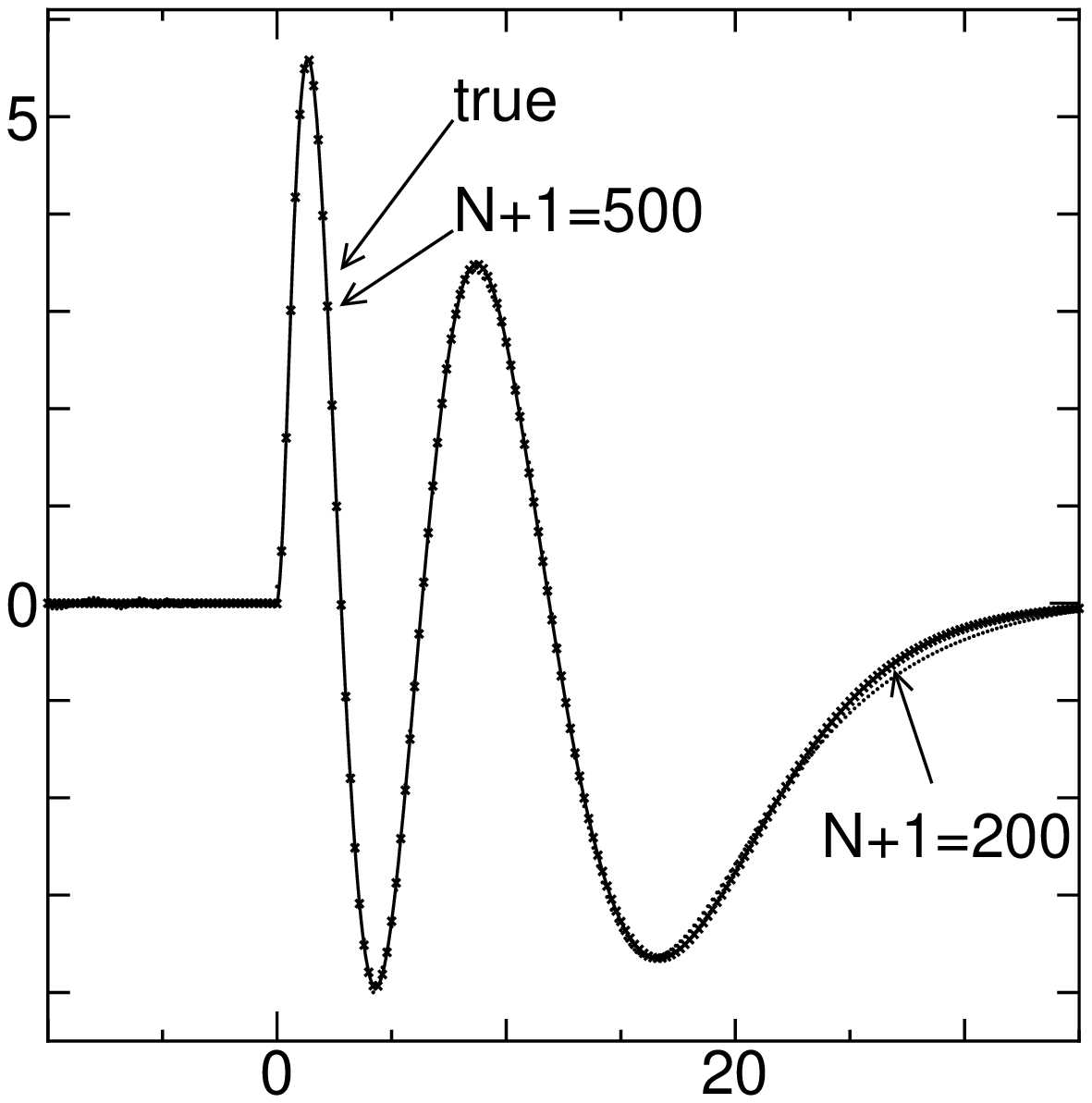}}
\caption{Numerical results for ODE (\ref{eqn:wtdasslag2}), i.e. for ODE (\ref{eqn:wtdasslag})}
\label{fig:wtdasslag}
\end{minipage}
\begin{minipage}{0.4\linewidth}
\includegraphics{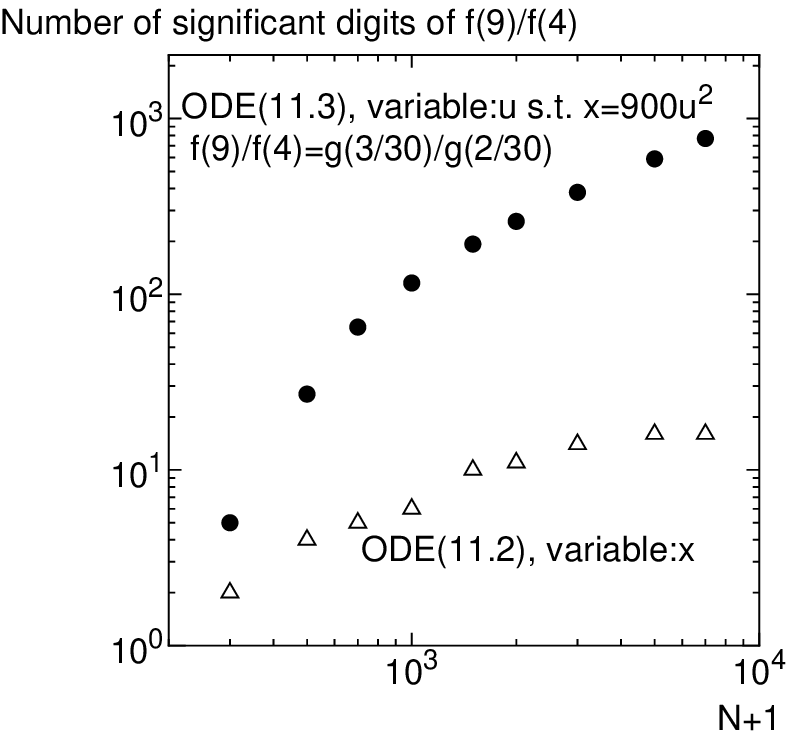}
\caption{Comparison of the number of significant digits of the ratio $f(9)/f(4)$ between ODE (\ref{eqn:wtdasslag2}) and ODE (\ref{eqn:wtdasslag3}) }
\label{fig:digits_v_wtdasslag}
\end{minipage}
\end{center}
\end{figure}

In the following, we will give the second example, for the Fuchsian-type ODE
\begin{eqnarray}
xf^{\prime\prime}(x)+f^{\prime}(x)+\left(-\frac{x}{4}+(\nu+\frac{\mu+1}{2})-\frac{\mu^2}{x}\right)f(x)=0
\label{eqn:wtdasslag}
\end{eqnarray}
with nonnegative integers $\mu$ and $\nu$. 
Since we have the Fuchsian-type 
ODE for $\phi(x):=x^{-\mu /2} \, e^{x/2} \, f(x)$ 
\begin{eqnarray*}
x\phi^{\prime\prime}(x)+(\mu+1-x)\phi^{\prime}(x)+\nu f(x)=0
\end{eqnarray*}
which is the associate Laguerre differential equation, it is easily shown that the solutions in $L_{(k_0)}^2(\mathbb{R})$ are proportional to the weighted associate Laguerre function 
\begin{eqnarray*}
f(x)=\left\{ \begin{array}{@{\,}ll} C\, 
x^{\mu /2} \, e^{-x/2} L_\nu^\mu (x)  &   (x\ge 0) \\ 
0 & (x<0)  \end{array} \right. \hspace{5mm}(C\in\mathbb{R})
\end{eqnarray*}
with the associate Laguerre polynomial $L_\nu^\mu (x)$. From the last discussion in Section \ref{sec:sv}, the ODE (\ref{eqn:wtdasslag}) can be treated by our algorithm as the 
Fuchsian-type 
ODE 
\begin{eqnarray}
x^2f^{\prime\prime}(x)+xf^{\prime}(x)+\left(-\frac{x^2}{4}+(\nu+\frac{\mu+1}{2})x-\mu^2\right)f(x)=0
\label{eqn:wtdasslag2}
\end{eqnarray}
whose coefficient functions are polynomials. In this ODE, the polynomial $p_2(x)=x^2$ of the highest order term has zero point at $x=0$. 
However, as has been found in many other 
Fuchsian-type ODEs, 
%ODEs with zero points of  $p_M(x)$, 
numerical results always converge to true solutions empirically. 
So is this case, and here we show how the results converge to true solutions in this case.

In Fig \ref{fig:wtdasslag}, the results with $\mu =4$, $\nu =3$, $k_0=0$, $N+1=200, 500$, \par\noindent 
$K=2\lfloor \frac{3N}{8} \rfloor + k_0$ and $J=2\lfloor \frac{7N}{16} \rfloor + k_0$   
are shown, under the normalization  \par\noindent 
 $ \langle f, \frac{1}{2\pi} (\psi_{k_0,0}+\psi_{k_0,-k_0-1})\rangle_{{\cal H}}=1$. The error of the result with $N+1=500$ is almost invisible there. It is remarkable that the obtained numerical results are very close to  zero for $x<0$ even though the basis functions are not small for $x<0$, though a small oscillation is observed in the result with $N+1=200$. However, the convergence in this case is less rapid than the ODEs without zero points of $p_M(x)$, because of the singularity of the solution at $x=0$ where the  solutions are not $(\mu /2 +1)$-th order differentiable.

To avoid this problem, instead of ODE (\ref{eqn:wtdasslag}), we solve the ODE 
\begin{eqnarray}
u^2g^{\prime\prime}(u)+ug^{\prime}(u)+\left(-c^4u^4+c^2(4\nu +2\mu+2)u^2-\mu^2\right)g(u)=0
\label{eqn:wtdasslag3}
\end{eqnarray}
with a positive constant $c$, which is derived directly from (\ref{eqn:wtdasslag}) by the change of coordinate $x=(cu)^2$ (where $g(u):=f\left((cu)^2\right)$. By this change of coordinate, we can obtain the solutions only for $x\ge 0$, which causes no inconvenience because the true solutions are zero for $x<0$. By this change of coordinate, the accuracy of the numerical result improves drastically, as is shown in Fig. \ref{fig:digits_v_wtdasslag} where the number of significant digits of the ratio $\displaystyle\frac{f(9)}{f(4)}$ are compared between the  ODEs (\ref{eqn:wtdasslag2}) and  (\ref{eqn:wtdasslag3}) with $c=30$.  As is shown in this figure, the number of significant digits empirically increases as a power of $N$ for $N+1>10^3$ in the case of ODE (\ref{eqn:wtdasslag3}), while it increases approximately proportional to $\log N$ due to the singularity at $x=0$ in the case of (\ref{eqn:wtdasslag2}). The bad behavior for ODE (\ref{eqn:wtdasslag2}) is theoretically deduced also from the fact that the order of coefficient $f_n$ in the expansion $f(x)=\sum_n f_n e_n(x)$ is an inverse power of $N$ for a function $f(x)$ with a singularity of this type, in a similar way to the case of Fourier series. The above change of coordinate eliminates this singularity, and it is successful in improving the accuracy to a great extent, up to several hundred digits, as is shown in Figure \ref{fig:digits_v_wtdasslag}. Moreover, the number of significant digits increases almost in a power of $N$ when $N+1\ge 2000$. This implies that the amount of required calculations increases almost in a polynomial order of the number of required significant digits, from the reasons shown in Subsection \ref{subsec:nc} below.

\subsection{Theoretical comparison with existing methods}
\label{subsec:nc}
In this subsection, we compare theoretically the order of the amount of calculations required for a very high accuracy, between some typical existing methods with arbitrary-precision arithmetic and the proposed method. The comparison is made with the following three methods:
\begin{description}
\item
(a) Runge-Kutta methods with arbitrary-precision arithmetic.
\item
(b) Finite element methods with arbitrary-precision arithmetic.
\item
(c) Petrov-Galerkin method with arbitrary-precision arithmetic using the same globally smooth basis functions as this paper. 
\end{description}

In the following, we compare them for the case where we require $Q$ significant digits for the ratio $f(x_1)/f(x_0)$ between the values of a solution function at two points $x_0$ and $x_1$. 

\begin{description}

\item[Order of amount of required calculations in the proposed method]
The amount of required calculations in the proposed method is almost a power of $Q$ when $Q$ is very large, because it requires about $O(N^3 (\log N)^2)$ as was explained in Section \ref{sec:cp} (after Theorem \ref{thm:opt_ratio_gen}) and many numerical resuts show that $Q$ is almost proportional to a power of $N$ when $N$ is very large. Empirically the amount of required calculations is from $O(Q^3)$ to $O(Q^5)$ in many cases. (Moreover, from the discussion in Section \ref{sec:cp}, we may redice it to be from $O(Q^2)$ to $O(Q^4)$ by means of the modification proposed after  Theorem \ref{thm:opt_ratio_gen}.)

\item[Comparison with (a) Runge-Kutta methods]
As is shown in the followong, Method (a) requires more than an exponential order of $Q$. As is known well, the discretization error of the Runge-Kutta methods is proportional to a power $h^p$ $(p\ge 4)$ of the step size $h$ for the discretization of the coordinate. (For example, in the common Runge-Kutta method, proportional to $h^4$.) This implies that the step size $h$ should satisfy the inequality $C_0 h^p < 10^{-Q}$ with a positive constant $C_0$. This inequality implies the inequality $\displaystyle \log h < - \frac{\log 10}{p} Q + C_1$ with another constant $C_1$. Since the numbers $n_s$ of the steps between $x_0$ and $x_1$ is $\displaystyle \frac{|x_0 -x_1|}{h}$, the number of required steps $n_s$ should satisfy the inequality $\displaystyle \log n_s > \frac{\log 10}{p} Q + C_2$ with another constant $C_2$. Hence, $\displaystyle n_s > C_3 \, (10)^{\frac{Q}{p}}$ with a positive constant $C_3$. Moreover, the number of required  digits for the working precision is larger as $Q$ is larger. (At least, it should be larger than $Q$.) \, These fact implies that the amount of required calculations is at least $\displaystyle O(Q^2 \cdot (10)^{\frac{Q}{p}})$. Can you imagine how huge is the amount when $Q$ is several hundreds or several thousands? Hence, the proposed method requires a much smaller amount of calculation than the Runge-Kutta methods, for a very high accuracy.

\item[Comparison with (b) finite element methods]
As is shown in the following, also Method (b) requires more than an exponential order of $Q$. As is known well, the error due to finite-dimensional approximation in the finite element methods is proportional to a power of the support size of the finite elements, at least. For example, when the basis functions are piece-wise polynomials of degree $q$ with support size $d$, the error is approximately proportional to $d^{\, q+1}$ at least, which is easily shown by the Taylor expansions of the true solution $f(x)$. Since the required dimension of the subspace is proportional to $\displaystyle \frac{|x_0 -x_1|}{d}$ and the matrix is band-diagonal, from a discussion very similar to the above Runge-Kutta cases, the amount of required calculations is at least $\displaystyle O(Q^2 \cdot (10)^{\frac{Q}{q+1}})$. This will be huge when $Q$ is very large. Hence, the proposed method requires a much smaller amount of calculation than the finite element  methods, for a very high accuracy.

\item[Comparison with (c) Galerkin methods using globally smooth basis functions]
For Method (c), since the direct order comparison is difficult, here we only point out that it requires a very large amount of calculations. Let $W$ be the number of required digits for the working precision, and $N+1$ be the dimension of the subspace. Then, since the matrix is band-diagonal, the amount of required calculations is $O(W^2 N)$. (If we use fast multiplication algorithms, the Karatsuba algorithm for example, the order may decrease to some extent. However, in our numerical examples, we did not use such algorithms. If we use such algorithms, we can diminish the order to the same extent as that. Therefore, for simplicity, here we compare the order without such algorithms.) \,  

In Method (c), it is very difficult to calculate the eigenvector in a high accuracy because of a heavy `cancelling' due to round-off errors, by the following reason, even if the exact eigenvalue is known. In this method, the matrix is band-diagonal with band width $2\ell _0+1$~\cite{paper1}. When we calculate the elements $f_n$ ($n=0,1,\ldots ,N)$ of the solution vector $\vec{f}$, with unknown intial values $f_0, f_1, \ldots f_{\ell_0-1}$, we should determine these initial values so that the linear equations given by the bottom $\ell _0$ rows of the matrix can be satisfied. This problem can be reduced a system of $\ell _0$ inhomogeneous simultaneous linear equation represented by a $\ell_0\times \ell _0$ matrix.
 
However, this $\ell_0\times \ell _0$ matrix is usually very close to a singular matrix of rank 1, because of the most diverging component contained in the halfway of the calculations of $f_n$ ($n=\ell_0, \ell_0+1, \ldots N-\ell_0$). In other words, with whatever initial values, the vector obtained in the halfway of the calculation are almost parallel to this diverging component, because this diverging component is dominant. Moreover, that $\ell_0\times \ell _0$ matrix is more close to a singular matrix of rank 1, as the dimension $N+1$ is larger. 

Though the order of this approach is difficult to estimate because it depends on the differential equation, anyway we should choose $W$ so that $\displaystyle\frac{W}{(\ell_0-1)Q}$ can be much larger than a monotonously increasing function of $N$ not smaller than $1$, in order to avoid the above mentioned `cancelling'. The total order estimation of this case and the comparison with the proposed method is one of future problem. 

Even if we use other methods for the calculation of the eigenvector, the Gaussian elimination or the diagonalization for example, the basic mathematical structure is the same as the above, and the problem due to the closeness to a linear dependence arises there.

Even if we use the Galerkin methods using another type of globally smooth basis functions, the Hermite functions for example, the basic circumstance is still similar to this. 

\end{description}

Thus, the proposed method requires a much smaller amount of calculations than Runge-Kutta and finite element methods, in order to calculate the values of the solution function $f(x)$ at a finite number of points in an extraordinarily high accuracy. This is the reason why we have many nurerical results by ordinary personal computers where $Q$ is so large as several hundreds or several thousands.

\section{Conclusions}

We have explained how to realize the integer-type algorithm proposed in~\cite{paper1} for linear higher-order differential equations, and proved that the realization proposed there satisfies the conditions given in {\bf C7} which are required for the convergence of numerical results to the true `general solutions' in 
${\cal H}$ 
%$C^M(\mathbb{R} )\cap {\cal H}$  
of the differential equations. 

We have proved that the method based on quasi-orthogonalization can obtain vectors in the quasi-optimal set $(\sigma_{K,N}^{(\Omega )})^{-1}[0,c \underline{\sigma_{K,N}^{(\Omega )}} ]$ with  fixed $c$. Moreover, we have provided some proofs about   the extension for the case with $\dim V\cap \ell ^2(\mathbb{Z}^+)\ge 2$ and the %termination 
halting  
of the procedures.

In addition, we have given a theoretical upper bound for the errors as a function only of the worst truncation error, in terms of one unknown parameter $\Delta _K$ (the worst truncation error of the exact true solutions).   

Numerical results have shown the precision of the proposed method. In addition, we have given an example which attains the exact ratios among the coefficients $f_n$ in the expansion $f(x)=\sum_n f_ne_n(x)$.

In the near future, we will compare the norm of actual errors in numerical results and the upper %limit 
bound  
proposed in this paper, and investigate the tightness of this upper 
%limit
bound.  
Remaining problems include improvements of the algorithm to reduce the number of calculations and better choices of the bilinear form $Q(\vec{f},\vec{g})$, and the orthogonality parameters $h$ and $g$.  

\section*{Acknowledgments}
MH was partially supported by MEXT through a Grant-in-Aid for Scientific Research in the Priority Area "Deepening and Expansion of Statistical Mechanical Informatics (DEX-SMI)", No. 18079014 and
a MEXT Grant-in-Aid for Young Scientists (A) No. 20686026. 
The Center for Quantum Technologies is funded by the Singapore Ministry of Education
and the National Research Foundation as part of the Research Centres of Excellence
programme.

\appendix

\section{Other orthogonality-like relations for %the wavepackets 
$\psi_{k,n} $ }
\label{app:oort}

With respect to the usual $L^2$-inner product  
\begin{eqnarray*}
(f, \,g) := \int_{-\infty}^\infty f(x) \, \overline{g(x)} \, dx \,\, ,
\end{eqnarray*}
the orthogonality-like relation 
\begin{eqnarray*}
(\psi_{k,\, n} , \psi_{k,\, n'} ) = 
\left\{ \begin{array}{@{\,}ll} \displaystyle 
\frac{\pi}{4^k} \,\, (-1)^{n-n'} \frac{(2k)!}{(k+n-n')!\,\, (k+n'-n)!}\,\, &  (\,|n-n'|\le k\,) \\ \\ \displaystyle 
0 & ({\rm otherwise}) 
\end{array} \right. 
\end{eqnarray*}
holds. 
With respect to another inner product 
\begin{eqnarray*}
\langle f, \,g\rangle _{|D|^{-k}} := \int_{-\infty}^\infty ({\cal F}f)(y) \,\overline{({\cal F}g)(y)} \,\, \frac{dy}{|y|^k} 
\end{eqnarray*}
(where ${\cal F}$ denotes the operator of Fourier transformation), another type of \par\noindent   orthogonality-like relation 
\begin{eqnarray*}
 ^\forall n, \, n'\in \mathbb{Z}^+, \,\,\,\,\,\,\,\, \langle \psi_{k,\, n} , \psi_{k,\, n'}\rangle _{|D|^{-k}} = \frac{\pi}{16^k} \,\frac{n!}{(n+2k)!} \,\, \delta_{n n'} 
\end{eqnarray*}
holds. This relation is derived from the orthogonality of the number states associated with the algebra $\mathfrak{su}(1,1)$ ~\cite{SaHa}. 
Here note that $\langle f, \, g\rangle _{(0)} = \langle f, \, g\rangle _{|D|^{-0}} = (f, \, g) $.

\section{On the shape of the wavepackets $\psi_{k,\, n}(x)$}
\label{app:shape}

The wavepackets are complex-valued `wavy' functions. By the scale change $x \to x/\sqrt{k}$, the `envelope' function 
\begin{eqnarray}
|\psi_{k,n} (x)| = (x^2 +1)^{-\frac{k+1}{2}}
\label{eqn:envelope}\end{eqnarray}
of the wavepacket $\psi_{k,\, n}^(x)$ is proportional to the probability density function of the (Student) $t$-distribution with the degree of freedom $k$ (used in statistics), which tends to the standard  Gaussian function as $k\to\infty$. Since the `foot of the mountain' of the probability density function of the $t$-distribution is thicker when $k$ is smaller, the `localization' of the envelope of the wavepacket is better as $k$ becomes larger when a comparison is made with respect to the normalized `width' (or standard deviation). 
The `phase' function of the wavepacket is 
\begin{eqnarray}
\arg \, \psi_{k,\, n} (x) = (k+2n+1) \,\, ( \arctan x -\frac{\pi}{2} ) \, .
\label{eqn:phase}\end{eqnarray} 
This function is almost linear when $|x|$ is not too large (because then $\arctan x\approx x$), and we can approximate the wavepacket by a sinusoidal wavepacket with the above envelope function, because the envelope function is sufficiently small where the phase function deviates considerably from the linear function $ (k+2n+1) \, ( x -\frac{\pi}{2} )$. This approximate picture is very good especially in the cases with large $k$ though it is poor for the cases with small $k$, because the `foot of the mountain' of the envelope vanishes very rapidly for large $k$. In fact, for sufficiently large $k$, 
\begin{eqnarray*}
\psi_{k,n} (x) \approx (- i)^{k+2n+1} \,\,\,  e^{ i(k+2n+1)x} \,\,\, (x^2 +1)^{-\frac{k+1}{2}} \\ {\rm or} \,\,\,\,\,\,\,\,\,\,\,\, \psi_{k,n} (x) \approx {\textstyle \frac{\sqrt{k \pi} \, \Gamma (\frac{k}{2})}{\Gamma ({\frac{k+1}{2})}}} \,\,\, e^{i(k+2n+1)x} \,\,\, e^{-\frac{k}{2}x^2} \,\,\, .
\end{eqnarray*}
For this, we can prove easily that the function $\displaystyle \, \Xi_k \, e^{- i \xi_{k,n}x} \, \psi_{k,n} (\frac{x}{\sqrt{k}}) \,\,$ with  \par\noindent $\displaystyle \Xi_k :=  \frac{\Gamma ({\frac{k+1}{2})}}{\pi \sqrt{2k} \, \Gamma (\frac{k}{2})}$ and $\displaystyle \xi_{k,n}:=\frac{k+2n+1}{\sqrt{k}} \, $ converges to the standard Gaussian function \par\noindent $\displaystyle\frac{1}{\sqrt{2\pi}} \, e^{-\frac{1}{2}x^2}$ as $k \to \infty$ for the $L^2$-norm, by means  of the upper and lower bounds of the Stirling formula, because the properties (\ref{eqn:envelope}) and (\ref{eqn:phase}) lead us to   
\begin{eqnarray}
 \label{eqn:zzzzzz}\hspace{1cm}
&& \Xi_k \, e^{- i \, \xi_{k,n}x} \, \psi_{k,n} (\frac{x}{\sqrt{k}}) - \frac{1}{\sqrt{2\pi}} \, e^{-\frac{1}{2}x^2}
\\&&  = \frac{1}{\sqrt{2\pi}} \, \left(  e^{- i \, (k+2n+1)  \bigl(\frac{x}{\sqrt{k}} \, - \, \arctan \frac{x}{\sqrt k} \bigr)} \, 
\sqrt{\frac{2\pi \,\, \Xi_k^2 \,\, e^{x^2}}{\displaystyle\bigl(1+\frac{x^2}{k}\bigr)^{k+1}}} \,\,  - \, \, 1 \, \right) \,
e^{-\frac{1}{2}x^2} . \nonumber 
\end{eqnarray}

This picture of the wavepackets, `almost-sinusoidal' oscillations with a spindle-shaped envelope, which are similar to Gaussian-weighted (complex) sinusoidal \par\noindent wavepackets, is `natural' and useful in many  applications, and it is very convenient for the interpretation of the basis systems used in this paper (for this, the relationship between these wavepackets and the Fourier series is shown in Appendix B of ~\cite{paper2}).

\end{document}